\documentclass[small]{article}

\usepackage{amsmath,amstext,amsthm,amsfonts}
\usepackage{makeidx}
\usepackage{amssymb,latexsym}
\usepackage{txfonts, pxfonts}
\usepackage[english]{babel}

\numberwithin{equation}{section}
\theoremstyle{plain}
\newtheorem{theorem}{Theorem}[section]
\newtheorem{proposition}[theorem]{Proposition}
\newtheorem{lemma}[theorem]{Lemma}
\newtheorem{conjecture}[theorem]{Conjecture}
\newtheorem{corollary}[theorem]{Corollary}

\theoremstyle{definition}
\newtheorem{definition}[theorem]{Definition}

\newtheorem*{notation}{\textbf{Notation}}

\newtheorem*{convention}{\textbf{Convention}}
\newtheorem{remark}[theorem]{Remark}

\def\e{\epsilon}

\def\N{\mathbb{N}}
\def\R{\mathbb{R}}
\def\Rn{{\mathbb{R}}^n_+}
\def\d{\partial}
\def\dst{d\sigma_{g(t)}}
\def\dvt{dv_{g(t)}}
\def\dv{dv}
\def\ds{d\sigma}
\def\D{\Delta}

\def\a{\alpha}
\def\b{\beta}
\def\l{\lambda}
\def\g{g_{\nu}}
\def\cesc{R_{g(t)}}
\def\cmedia{H}
\def\cm{H_{g(t)}}
\def\cmbar{\overline{H}_{g(t)}}
\def\cmzbar{\overline{H}_{g_0}}
\def\cminfbar{\overline{H}_{\infty}}
\def\cmz{H_{g_0}}
\def\cmv{H_{\g}}
\def\cmvbar{\overline{H}_{\g}}

\def\u{u_{\nu}}
\def\uinf{u_{\infty}}
\def\uinu{\bar{u}_{(x_{i,\nu},\e_{i,\nu})}}
\def\ujnu{\bar{u}_{(x_{j,\nu},\e_{j,\nu})}}

\def\uknu{\bar{u}_{(x_{k,\nu},\e_{k,\nu})}}
\def\uk{\bar{u}_{(x_{k},\e_{k})}}

\def\w{w_{\nu}}
\def\crit {\frac{2n}{n-2}}
\def\critbordo{\frac{2(n-1)}{n-2}}

\def\conjbordo{\frac{2(n-1)}{n}}\def\ba{\begin{align}}
\def\ea{\end{align}}
\def\bp{\begin{proof}}
\def\ep{\end{proof}}
\def\Q{Q(B^n,\d B)}
\def\QM{Q(M,\d M)}
\def\s{\sigma}

\def\uum{\bar{u}_{(x_1,\e_1)}}
\def\udois{\bar{u}_{(x_2,\e_2)}}
\def\ubar{\bar{U}_{(x_0,\e)}}
\def\func:u{\bar{u}_{(x_0,\e)}}
\def\U{U_{\epsilon}}

\begin{document}

\title{Convergence of scalar-flat metrics on manifolds with boundary under a Yamabe-type flow}

\author{\textsc{S\'ergio Almaraz}~\footnote{Supported by FAPERJ, CAPES and CNPq (Brazil).}}
\date{}

\maketitle

\begin{abstract}
We study a conformal flow for compact Riemannian manifolds of dimension greater than two with boundary. Convergence to a scalar-flat metric with constant mean curvature on the boundary is established in dimensions up to seven, and in any dimensions if the manifold is spin or if it satisfies a generic condition. 
\end{abstract}

\tableofcontents


\section{Introduction}\label{sec:intr}
Let $M^n$ be a closed manifold with dimension $n\geq 3$. In order to solve the Yamabe problem (see \cite{yamabe}), R. Hamilton introduced the Yamabe flow, which evolves Riemannian metrics on $M$ according to the equation
$$
\frac{\d}{\d t}g(t)=-(R_{g(t)}-\overline{R}_{g(t)})g(t)\,,
$$ 
where $R_g$ denotes the scalar curvature of the metric $g$ and $\overline{R}_g$ stands for the average $\displaystyle{\left(\int_M \dv_g\right)^{-1}}\int_M R_g\dv_g$. Here, $\dv_g$ is the volume form of $(M,g)$. Although the Yamabe problem was solved using a different approach in \cite{aubin1,schoen1,trudinger}, the Yamabe flow is a natural geometric deformation to metrics of constant scalar curvature. The convergence of the Yamabe flow on closed manifolds was studied in \cite{chow, struwe-flow, ye}. This question was completed solved in \cite{brendle-flow, brendle-invent}, where the author makes use of the positive mass theorem.

In this work, we study the convergence of a similar flow on compact n-dimensional manifolds with boundary, when $n\geq 3$. For those manifolds, J. Escobar raised the question of existence of conformal scalar-flat metrics on $M$ which have the boundary as a constant mean curvature hypersurface. This problem was studied in \cite{ahmedou, escobar3,escobar4, marques-weyl,marques-4-5, almaraz2, chen}. (The question of existence of conformal metrics with constant scalar curvature and minimal boundary was studied in \cite{brendle-chen,escobar2}; see also \cite{ambrosetti-li-malchiodi,han-li}.)
 
Let $(M^n,g_0)$ be a compact Riemannian manifold with boundary $\d M$ and dimension $n\geq 3$.
We consider the following conformal invariant defined in \cite{escobar3}:
\ba
Q(M,\d M)&=\inf_{g\in [g_0]}
\frac{\int_M R_g\dv_g+2\int_{\d M}\cmedia_g\ds_g}{\left(\int_{\d M}\ds_g\right)^{\frac{n-2}{n-1}}}\notag
\\
&=\inf_{\{\:u\in C^1(\bar{M}), u\nequiv 0 \:\text{on}\: \d M\}}
\frac{\int_M \left(\frac{4(n-1)}{n-2}|d u|_{g_0}^2+R_{g_0}u^2\right)dv_{g_0}+\int_{\d M}2\cmedia_{g_0}u^2d\sigma_{g_0}}
{\left(\int_{\d M}|u|^\frac{2(n-1)}{n-2}d\sigma_{g_0}\right)^{\frac{n-2}{n-1}}}\,,\notag
\end{align}
where $\cmedia_g$ and  $\ds_g$ denote respectively the trace of the 2nd fundamental form and the volume form of $\d M$, with respect to the metric $g$, and $[g_0]$ stands for the conformal class of the metric $g_0$. Although we always have  $Q(M, \d M)\leq Q(B^n, \d B)$, where $B^n$ is the closed unit ball in $\R^n$,  we may have $Q(M, \d M)=-\infty$ (see \cite{escobar3'}).

Conformal scalar-flat metrics in compact manifolds with boundary can be easily obtained under the hypothesis that  $Q(M,\d M)>-\infty$  (which is the case when the scalar curvature is non-negative). To that end, we can use, as the conformal factor, the first eigenfunction of a linear eigenvalue problem (see  \cite[Proposition 1.4]{escobar3}).

We are interested in a formulation of a Yamabe-type flow for compact scalar-flat manifolds with boundary proposed by S. Brendle in \cite{brendle-boundary}. This flow evolves a conformal family of metrics $g(t)$, $t\geq 0$, according to the equations
\begin{equation}\label{eq:evol}
\begin{cases}
\cesc=0\,,&\text{in}\:M\,,
\\
\frac{\d}{\d t}g(t)=-2(\cm-\cmbar)g(t)\,,&\text{on}\:\d M\,,
\end{cases}
\end{equation}
where, $\overline{\cmedia}_g$ stands for the average $\displaystyle{\left(\int_{\d M} \ds_g\right)^{-1}}\int_{\d M} \cmedia_g\ds_g $. (We refer the reader to Section \ref{sec:prelim} for the formulation in terms of the conformal factor.)

Brendle proved short-time existence of a unique solution to (\ref{eq:evol}) for a given initial metric and the following long-time result: 
\begin{theorem}[\cite{brendle-boundary}]\label{brendle:boundary:thm}
Suppose that:

(i) $Q(M,\d M)\leq 0$, \: or

(ii) $Q(M,\d M)>0$, $M$ is locally conformally flat with umbilic boundary, and the boundary of the universal cover of $M$ is connected.

Then, for every initial scalar-flat metric $g(0)$ on M, the flow \eqref{eq:evol} exists for all time $t\geq 0$ and converges to a scalar-flat metric with constant mean curvature on the boundary.
\end{theorem}

Inspired by the ideas in \cite{brendle-flow, brendle-invent}, we handle the remaining cases of this problem.
Define
\begin{equation*}
\mathcal{Z}=
\{x_0\in\d M\,;\:
\limsup_{x\to x_0}d_{g_0}(x,x_0)^{2-d}|W_{g_0}(x)|
=\limsup_{x\to x_0}d_{g_0}(x,x_0)^{1-d}|\pi_{g_0}(x)|=0\}\,,
\end{equation*}
where $W_{g_0}$ denotes the Weyl tensor of $M$, $\pi_{g_0}$ the trace-free second fundamental form of $\d M$, and $d=\Big{[}\frac{n-2}{2}\Big{]}$.
Our first result is the following:
\begin{theorem}\label{first:thm}
Suppose that $(M^n,g_0)$ is not conformally diffeomorphic to the unit ball $B^n$ and satisfies $Q(M,\d M)>0$. 
If 

(a) $\mathcal{Z}= \emptyset$, \: or

(b) $n\leq 7$, \: or

(c) $M$ is spin,
\\
then, for any initial scalar-flat metric $g(0)$, the flow \eqref{eq:evol} exists for all time $t\geq 0$ and converges to a scalar-flat metric with constant mean curvature on the boundary.
\end{theorem}

Since Euclidean domains are spin, the following is an immediate consequence of Theorems \ref{brendle:boundary:thm} and \ref{first:thm}:

\begin{corollary}

If $M\subset \R^n$ is a compact domain with smooth boundary, then the flow \eqref{eq:evol}, starting with any scalar-flat metric, exists for all time $t\geq 0$ and converges to a scalar-flat metric with constant mean curvature on the boundary.
\end{corollary}

Condition (a) in Theorem \ref{first:thm} is particularly satisfied if  the trace-free second fundamental form of $\d M$ is nonzero everywhere.
In dimensions $n\geq 4$, the set of metrics which satisfy this latter condition  is an open and dense subset of the space of all Riemannian metrics on $M$. This hypothesis was used in \cite{almaraz1} to prove compactness of the set of solutions to the Yamabe problem on manifolds with boundary.

Conditions (b) and (c) allow us to make use of a positive mass theorem for manifolds with a non-compact boundary, very recently proved in \cite{almaraz-barbosa-lima}. 

Before stating our main result, from which Theorem \ref{first:thm} follows, we will discuss this positive mass theorem and the concept of mass for those manifolds.

Let $(N, g)$ be a Riemannian manifold with non-compact  boundary $\d N$. 
\begin{definition}\label{def:asym}
We say that $N$ is {\it{asymptotically flat}} with order $p>0$, if there is a compact set $K\subset N$ and a diffeomorphism $f:N\backslash K\to \Rn\backslash \overline{B^+_1(0)}$ such that, in the coordinate chart defined by $f$ (called {\it  asymptotic coordinates} of $M$), we have
$$
|g_{ab}(y)-\delta_{ab}|+|y||g_{ab,c}(y)|+|y|^2|g_{ab,cd}(y)|=O(|y|^{-p})\,,
\:\:\:\:\text{as}\:\:|y|\to\infty\,,
$$
where $a,b,c,d=1,...,n$.
\end{definition}

Here, $\Rn=\{(y_1,...,y_n)\in\R^n\,;\:y_n\geq 0\}$ and $\overline{B^+_1(0)}=\{y\in\Rn\,;\:|y|\leq 1\}$.

Suppose that $N^n$, with dimension $n\geq 3$,  is asymptotically flat with order $p>\frac{n-2}{2}$. Let $(y_1,...,y_n)$ be the  asymptotic coordinates induced by the diffeomorphism $f$ as above. We also assume that $R_g$ is integrable on $N$, and $\cmedia_g$  is integrable on $\d N$. Then the limit
\begin{align}\label{def:mass}
m(g)&:=
\\
&\lim_{R\to\infty}\left\{
\sum_{a,b=1}^{n}\int_{y\in\Rn,\, |y|=R}(g_{ab,b}-g_{bb,a})\frac{y_a}{|y|}\,\ds_{R}
+\sum_{i=1}^{n-1}\int_{y\in\d\Rn,\, |y|=R}g_{ni}\frac{y_i}{|y|}\,\ds_{R}\right\}\notag
\end{align}
exists, and we call it the {\it mass} of $(M, g)$ . Moreover, $m(g)$ is a geometric invariant in the sense that it does not depend on the asymptotic coordinates. (This definition of mass was presented to me by F. Marques.)
\begin{conjecture}[Positive mass]\label{pmt}
If $R_g, \cmedia_g\geq 0$, then we have $m(g)\geq 0$ and the equality holds if and only if $N$ is isometric to $\Rn$.
\end{conjecture}
In \cite{almaraz-barbosa-lima}, this conjecture is reduced to the case of manifolds without boundary, known in the spin case for any dimensions (\cite{witten}) and for $n\leq 7$ in general (\cite{schoen2,schoen-yau}),  so we have the following result:
\begin{theorem}[\cite{almaraz-barbosa-lima}]
Conjecture \ref{pmt} holds true if $n\leq 7$ or if $N$ is spin.
\end{theorem}
\begin{remark}
Special cases of this conjecture were previously obtained by S. Raulot in \cite{raulot} and by J. Escobar in the appendix of \cite{escobar2}.
\end{remark}

The asymptotically flat manifolds  used in this paper are obtained as the generalized stereographic projections of the compact Riemannian  manifold with boundary $(M,g_0)$. Those stereographic projections are performed around points $x_0\in \d M$ by means of the Green functions $G_{x_0}$, with singularity at $x_0$, obtained in Appendix B. After choosing a new background metric $g_{x_0}\in [g_0]$ with better coordinates expansion around $x_0$ (see Section \ref{sub:sec:deftestfunct}), we consider the asymptotically flat manifold $(M\backslash \{x_0\}, \bar{g}_{x_0})$, where $\bar{g}_{x_0}= G_{x_0}^{\frac{4}{n-2}}g_{x_0}$ satisfies $R_{\bar {g}_{x_0}}\equiv 0$ and $\cmedia_{\bar {g}_{x_0}}\equiv 0$. If $x_0\in \mathcal{Z}$, according to Proposition \ref{propo19}, this manifold has asymptotic order $p>\frac{n-2}{2}$, so Conjecture \ref{pmt} claims that $m(\bar{g}_{x_0})>0$ unless $M$ is conformally equivalent to the unit ball.

Our main result, which implies Theorem \ref{first:thm}, is the following:

\begin{theorem}\label{main:thm}
Suppose that $(M^n,g_0)$ is not conformally diffeomorphic to the unit ball $B^n$ and satisfies $Q(M,\d M)>0$. 

If $m(\bar{g}_{x_0})>0$ for all $x_0\in \mathcal{Z}$, then, for any initial scalar-flat metric $g(0)$, the flow \eqref{eq:evol} exists for all time $t\geq 0$ and converges to a scalar-flat metric with constant mean curvature on the boundary.
\end{theorem}

The proof of Theorem \ref{main:thm} follows the arguments in \cite{brendle-flow}. An essential step is the construction of a family of test functions on $M$, whose energies are uniformly bounded by the Sobolev quotient $\Q$. 
This construction is inspired by the test functions introduced by S. Brendle in \cite{brendle-invent} for the case of closed manifolds. 
The functions we use here were obtained in \cite{chen} in the case of umbilic boundary, where S. Chen addresses the existence of solutions to the Yamabe problem for manifolds with boundary, using an approach similar to the one in \cite{brendle-chen}. In the present work, we extend those functions to the case when the boundary does not need to be umbilic.    

Another crucial result used in the proof of our main theorem is the result in \cite{almaraz4}, which is a modification of a compactness theorem due to M. Struwe in \cite{struwe}; see also Chapter 3 of \cite{druet-hebey-robert} and \cite{cao-noussair-yan,chabrowski-girao,pierotti-terracini}.

This paper is organized as follows. In Section \ref{sec:prelim}, we establish some preliminaries and prove the long-time existence of the flow. In Section
\ref{sec:testfunc}, we construct the necessary test functions by modifying the arguments in \cite{chen}. In Section \ref{sec:blowup}, we make use of the compactness theorem in \cite{almaraz4} to carry out a blow-up analysis using the test functions. In Section \ref{sec:mainthm}, firstly we use the blow-up analysis to prove a result which is analogous to Proposition 3.3 of \cite{brendle-flow}. Then we use this result to prove the main theorem by estimating the solution to the flow uniformly in $t\geq 0$.  In Appendix A, we establish some elliptic estimates. In Appendix B, we construct the Green function used in this work and prove some of its properties.

\section{Preliminary results and long-time existence}\label{sec:prelim}
\begin{notation}
In the rest of this paper, $M^n$ will denote a compact manifold of dimension $n\geq 3$ with boundary $\d M$,  and $g_0$ will denote a background Riemannian metric on $M$. 
We will denote by $B_r(x)$ (resp. $D_r(x)$) the metric ball in $M$ (resp. $\d M$) of radius $r$ with center $x\in M$ (resp. $x\in \d M$). 

For any Riemannian metric $g$ on $M$, $\eta_g$ will denote the inward unit normal vector to $\d M$ respect to $g$ and $\Delta_g$ the Laplace-Beltrami operator,.

If $z_0\in \Rn$, we set $B_r^+(z_0)=\{z\in\Rn\,;\:|z-z_0|< r\}$,
$$
\d^+B^+_{r}(z_0)=\d B^+_{r}(z_0)\cap\Rn\,,
\:\:\:\:\text{and}\:\:\:\:\:
\d 'B^+_{r}(z_0)=B^+_{r}(z_0)\cap \d\Rn\,.
$$
Finally, for any $z=(z_1,...,z_{n-1}, z_n)$ we set $\bar{z}=(z_1,...,z_{n-1},0)\in \d\Rn\cong\R^{n-1}$. 
\end{notation} 
\begin{convention}
We assume that $(M,g_0)$ satisfies $Q(M,\d M)>0$. According to \cite[Proposition 1.4]{escobar3}, we can also assume that $R_{g_0}\equiv 0$ and $\cmz>0$, after a conformal change of the metric. Multiplying $g_0$ by a positive constant, we can suppose that $\int_{\d M}\ds_{g_0}=1$.

The Sobolev spaces $H^p(M)$ and $L^p(M)$ are defined with respect to the metric $g_0$, and $H^p(\d M)$ and $L^p(\d M)$ with respect to the induced metric on $\d M$.

We will adopt the summation convention whenever confusion is not possible, and use indices $a,b,c,d=1,...,n$, and $i,j,k,l=1,...,n-1$.
\end{convention}

If $g=u^{\frac{4}{n-2}}g_0$ for some  positive smooth function $u$ on $M$, we know that 
\begin{equation}\label{eq:R:H}
\begin{cases}
\displaystyle R_g=u^{-\frac{n+2}{n-2}}\left(-\frac{4(n-1)}{n-2}\Delta_{g_0}u+R_{g_0}\right)\,,&\text{in}\:M\,,
\\
\displaystyle \cmedia_g=u^{-\frac{n}{n-2}}\left(-\frac{2(n-1)}{n-2}\frac{\d}{\d \eta_{g_0}}u+\cmedia_{g_0} u\right)\,,&\text{on}\:\d M\,,
\end{cases}
\end{equation}
and the operators $L_{g}=\frac{4(n-1)}{n-2}\Delta_{g}u-R_{g}$ and $B_g=\frac{2(n-1)}{n-2}\frac{\d}{\d \eta_{g}}u-\cmedia_{g}$ satisfy
\begin{equation}\label{propr:L}
L_{u^{\frac{4}{n-2}}g_0}(u^{-1}\zeta)=u^{-\frac{n+2}{n-2}}L_{g_0}\zeta,
\end{equation}
\begin{equation}\label{propr:B}
B_{u^{\frac{4}{n-2}}g_0}(u^{-1}\zeta)=u^{-\frac{n}{n-2}}B_{g_0}\zeta\,,
\end{equation}
for any smooth function $\zeta$.

If $u(t)=u(\cdot,t)$ is a 1-parameter family of positive smooth functions on $M$ and $g(t)=u(t)^{\frac{4}{n-2}}g_0$ with $R_{g_0}\equiv 0$, then (\ref{eq:evol}) can be written as
\begin{equation}\label{eq:evol:u}
\begin{cases}
\Delta_{g_0}u(t)=0\,,&\text{in}\:M\,,
\\
\displaystyle\frac{\d}{\d t}u(t)=-\frac{n-2}{2}(\cm-\cmbar)u(t)\,,&\text{on}\:\d M\,.
\end{cases}
\end{equation}
The second equation of (\ref{eq:evol:u}) can also be written as 
\ba
\frac{\d}{\d t}u(t)&=(n-1)u(t)^{-\frac{2}{n-2}}\frac{\d}{\d \eta_{g_0}}u(t)
-\frac{n-2}{2}\cmedia_{g_0}u(t)^{1-\frac{2}{n-2}}\notag
\\
&\hspace{1cm}-\frac{n-2}{2}u(t)\int_{\d M}\left(\frac{2(n-1)}{n-2}\frac{\d}{\d \eta_{g_0}}u(t)-\cmedia_{g_0}u(t)\right)u(t)\ds_{g_0}\,.\notag
\end{align}
Recall that short-time existence of solutions to the equations (\ref{eq:evol:u}) was obtained in \cite{brendle-boundary}. 
Hence, those equations have a solution $u(t)$ defined for all $t$ in the maximal interval $[0,T_{max})$.

According to \cite[Lemma 3.8]{brendle-boundary}, the function $\cm$ on $\d M$ can be extended to a smooth function on $M$, also denoted by  $\cm$,  satisfying
\begin{equation}\label{eq:evol:H}
\begin{cases}
\Delta_{g(t)}\cm=0\,,&\text{in}\:M\,,
\\
\displaystyle\frac{\d}{\d t}\cm=(n-1)\frac{\d}{\d\eta_{g(t)}}\cm+\cm(\cm-\cmbar)\,,&\text{on}\:\d M\,.
\end{cases}
\end{equation}

Hence, the evolution equations for the volume form $\ds_{g(t)}$ of $\d M$ and for $\cmbar$ are given by
\begin{equation}\label{eq:evol:vol}
\frac{d}{dt}\ds_{g(t)}=-(n-1)(\cm-\cmbar)\ds_{g(t)}
\end{equation}
and 
\begin{equation}\label{eq:evol:Hbar}
\frac{d}{dt}\cmbar=-(n-2)\int_{\d M}(\cm-\cmbar)^2\ds_{g(t)}\,.
\end{equation}
In particular, we can assume that
\begin{equation}\label{bd:area:norm}
\int_{\d M}\ds_{g(t)}=1\,,
\:\:\:\:\:\text{for all}\:t\in [0, T_{max})\,,
\end{equation}
and  we see that $\cmbar$ is decreasing. 

The next proposition is a direct application of the maximum principle to the equations (\ref{eq:evol:H}).
\begin{proposition}\label{Propo2.1} We have 
$$
\inf_{\d M}\cm \geq \min\{\inf_{\d M}\cmedia_{g(0)},0\}\,,
\quad\quad\quad\text{for all}\:t\in [0,T_{max}).
$$  
\end{proposition}

Set 
$$
\s=1-\min\{0,\inf_{\d M}\cmedia_{g(0)}\}=\max\{\sup_{\d M}(1-\cmedia_{g(0)}),1\}\,.
$$ 
By Proposition \ref{Propo2.1}, we have $\cm+\sigma\geq 1$ for all $t\in [0,T_{max})$.

In order to prove that $T_{max}=\infty$, we will prove uniform estimates for $u(t)$ on $[0,T)$, if $T$ is finite.
\begin{proposition}\label{Propo2.4}
Let $0<T\leq T_{max}$. If $T<\infty$, then there exist $C(T),c(T)>0$ such that
\begin{equation}\label{Propo2.4:1}
\sup_{M}u(t)\leq C(T)\:\:\:\:\text{and}\:\:\:\:\inf_{M}u(t)\geq c(T)\,,\:\:\:\:\:\text{for all}\:\:t\in[0,T)\,. 
\end{equation}
\end{proposition}
\bp
It follows from the evolution equations (\ref{eq:evol:u}) and (\ref{eq:evol:Hbar}),  and from the inequality $\cm+\sigma\geq 1$ that
$$
\frac{\d }{\d t}\log u(t)=-\frac{n-2}{2}(\cm-\cmbar)\leq \frac{n-2}{2}(\overline{\cmedia}_{g(0)}+\sigma)\,,
\:\:\:\:\text{on}\:\d M\,.
$$
Since $T<\infty$, there exists $C(T)>0$ such that $\sup_{\d M}u(t)\leq C(T)$ for all $t\in[0,T)$, and the first estimate of (\ref{Propo2.4:1}) follows from the maximum principle.

In order to prove the second one, first we will prove that there exists $c(T)>0$ such that 
\begin{equation}\label{Propo2.4:2}
\|u(t)\|_{L^{\crit}(M)}\geq c(T)\,,
\:\:\:\:\:\text{for all}\:t\in[0,T)\,.
\end{equation} 

Suppose by contradiction this is not true. Then there exists a sequence $\{t_j\}_{j=1}^{\infty}\subset [0,T)$ such that $u_j=u(t_j)\to 0$ in $L^{\crit}(M)$ as $j\to\infty$. Using (\ref{eq:R:H}), (\ref{eq:evol:Hbar}), and the boundary area normalization (\ref{bd:area:norm}), we see that 
\begin{equation}\label{Propo2.4:3}
\int_{M}\critbordo|du(t)|^2_{g_0}\dv_{g_0}+\int_{\d M}\cmz u(t)^2\ds_{g_0}=\overline{\cmedia}_{g(t)}\leq \overline{\cmedia}_{g(0)}\,,\:\:\:\:\text{for all}\:t\geq 0.
\end{equation}
Hence, there exists $u_0\in H^1(M)$ such that, up to a subsequence, $u_j\rightharpoonup u_0$ in $H^1(M)$. By the Sobolev embedding theorems, we can also assume that $u_j\to u_0$ in $L^2(M)$ and, at the same time, $u_j\to u_0$ in $L^2(\d M)$. Since we are assuming $u_j\to 0$ in $L^{\crit}(M)$, we see that $u_0\equiv 0$ a.e., and thus $u_j\to 0$ in $L^2(\d M)$. Since $\sup_{\d M}u_j\leq C(T)$, it follows from interpolation that $u_j\to 0$ in $L^{\critbordo}(\d M)$. This contradicts the boundary area normalization and proves the estimate (\ref{Propo2.4:2}).
 
We set $P=\cmz+\sigma C(T)^{\frac{2}{n-1}}$
and observe that, for all $t\in [0,T)$,
\ba
-\critbordo\frac{\d u(t)}{\d\eta_{g_0}}+Pu(t)
&\geq -\critbordo\frac{\d u(t)}{\d\eta_{g_0}}+\cmz u(t)+\sigma u(t)^{\frac{n}{n-2}}\notag
\\
&=(\cm+\sigma)u(t)^{\frac{n}{n-2}}
\geq 0\,.\notag
\end{align}
Then it follows from Proposition \ref{PropoA.2} that there exists $c(T)>0$ such that
\begin{equation}\notag
\big(\inf_{M} u(t)\big)^{\frac{n-2}{2n}}\big(\sup_{M} u(t)\big)^{\frac{n+2}{2n}}
\geq c(T)\left(\int_{M}u(t)^{\crit}\dv_{g_0}\right)^{\frac{n-2}{2n}}\notag
\end{equation}
for all $t\in [0,T)$. Then the second estimate of (\ref{Propo2.4:1}) easily follows using the fact that $\sup_{M}u(t)\leq C(T)$.
\ep
Now we proceed as in \cite[p.642]{brendle-boundary} to conclude that, if $T$ is finite, all higher order derivatives of $u$ are uniformly bounded on $[0,T)$. This implies that $u(t)$ is defined for all $t\geq 0$.

\begin{notation}
We define 
\begin{equation}\label{eq:def:cminfbar}
\cminfbar=\lim_{t\to\infty}\cmbar
\end{equation}
and observe that $\cminfbar\geq \frac{1}{2}\QM>0$.
\end{notation}

Next we establish some auxiliary results to be used in the rest of the paper.
\begin{lemma}\label{Lemma2.2} For any $p>2$ we have
\ba
\frac{d}{dt}&\int_{\d M}(\cm+\s)^{p-1}\dst=\notag
\\
&-\frac{4(n-1)(p-2)}{p-1}\int_{M}\left|d(\cm+\s)^{\frac{p-1}{2}}\right|_{g(t)}^2\dvt\notag
\\
&-(n-p)\int_{\d M}\left\{(\cm+\s)^{p-1}-(\cmbar+\s)^{p-1}\right\}(\cm-\cmbar)\dst\notag
\\
&-(p-1)\int_{\d M}\s\left\{(\cm+\s)^{p-2}-(\cmbar+\s)^{p-2}\right\}(\cm-\cmbar)\dst\,.\notag
\end{align}
\end{lemma}
\bp
This lemma is a direct computation using the equations (\ref{eq:evol:H}) and (\ref{eq:evol:vol}).
\ep
\begin{lemma}\label{Lemma2.3}
For any $p>n-1$ there exists $C>0$ such that
\begin{align}\label{Lemma2.3:1}
\frac{d}{dt}\int_{\d M}|\cm-\cmbar|^p\dst
\leq
&\,C\int_{\d M}|\cm-\cmbar|^p\dst
\\
&+C\left\{\int_{\d M}|\cm-\cmbar|^p\dst\right\}^{\frac{p+2-n}{p+1-n}}\notag
\end{align}
for all $t$. 
\end{lemma}
\bp
From the evolution equations (\ref{eq:evol:H}), (\ref{eq:evol:vol}), and (\ref{eq:evol:Hbar}), we obtain
\ba
&\frac{d}{dt}\int_{\d M}|\cm-\cmbar|^p\dst\notag
\\
&=p(n-1)\int_{\d M}|\cm-\cmbar|^{p-2}(\cm-\cmbar)\frac{\d\cm}{\d \eta_{g(t)}}\ds_{g(t)}\notag
\\
&\hspace{0.5cm}+p\int_{\d M}|\cm-\cmbar|^{p}\cm\ds_{g(t)}\notag
\\
&\hspace{0.5cm}-(n-1)\int_{\d M}|\cm-\cmbar|^{p}(\cm-\cmbar)\ds_{g(t)}\notag
\\
&\hspace{0.5cm}+p(n-2)\int_{\d M}|\cm-\cmbar|^{p-2}(\cm-\cmbar)\ds_{g(t)}
\int_{\d M}(\cm-\cmbar)^2\ds_{g(t)}\,.\notag
\end{align}
Using the identity 
$$
p\int_M|f|^{p-2}f\Delta_g f\dv_g
+\frac{4(p-1)}{p}\int_M\big{|}d|f|^{\frac{p}{2}}\big{|}_g^2\dv_g
=-p\int_{\d M}|f|^{p-2}f\frac{\d f}{\d \eta_g}\ds_g\,,
$$
we can write
\ba
&\frac{d}{dt}\int_{\d M}|\cm-\cmbar|^p\dst\notag
\\
&=-\frac{(p-1)(n-2)}{p}\Big\{
\int_M\frac{4(n-1)}{n-2}\big{|}d|\cm-\cmbar|^{\frac{p}{2}}\big{|}^2\dv_{g(t)}\notag
\\
&\hspace{3cm}+\int_{\d M}2\cm|\cm-\cmbar|^{p}\ds_{g(t)}
\Big\}\notag
\\
&\hspace{0.5cm}+\left(\frac{2(p-1)(n-2)}{p}+p+1-n\right)\int_{\d M}|\cm-\cmbar|^{p}(\cm-\cmbar)\ds_{g(t)}\notag
\\
&\hspace{0.5cm}+\left(\frac{2(p-1)(n-2)}{p}+p\right)\int_{\d M}\cmbar|\cm-\cmbar|^{p}\ds_{g(t)}\notag
\\
&\hspace{0.5cm}+p(n-2)\int_{\d M}|\cm-\cmbar|^{p-2}(\cm-\cmbar)\ds_{g(t)}
\int_{\d M}(\cm-\cmbar)^2\ds_{g(t)}\,.\notag
\end{align}
Since $p>n-1$ and $\cmbar$ is nonincreasing, using H\"{o}lder's inequality and $\int_{\d M}\dst=1$, we obtain
\ba\label{Lemma2.3:2}
&\frac{d}{dt}\int_{\d M}|\cm-\cmbar|^p\dst
\\
&\hspace{0.5cm}\leq -\frac{(p-1)(n-2)}{p}Q(M,\d M)\left\{\int_{\d M}|\cm-\cmbar|^{\frac{p(n-1)}{n-2}}\ds_{g(t)}\right\}^{\frac{n-2}{n-1}}\notag
\\
&\hspace{1cm}+\left(\frac{2(p-1)(n-2)}{p}+p+1-n\right)\int_{\d M}|\cm-\cmbar|^{p+1}\ds_{g(t)}\notag
\\
&\hspace{1cm}+\left(\frac{2(p-1)(n-2)}{p}+p\right)\int_{\d M}\cmzbar|\cm-\cmbar|^{p}\ds_{g(t)}\notag
\\
&\hspace{1cm}+p(n-2)\left\{\int_{\d M}|\cm-\cmbar|^{p}\ds_{g(t)}\right\}^{\frac{p+1}{p}}\,.\notag
\end{align}

Applying the Young's inequality $AB\leq \a A^{\frac{1}{\a}}+(1-\a)A^{\frac{1}{1-\a}}$ to the interpolation inequality 
$\|f\|_{L^{p+1}(\d M)}^{p+1}\leq \|f\|_{L^{\frac{p(n-1)}{n-2}}(\d M)}^{\a p}\|f\|_{L^{p}(\d M)}^{1+(1-\a)p}$ with $\a=\frac{n-1}{p}<1$, we obtain
\ba
\int_{\d M}|\cm-\cmbar|^{p+1}\ds_{g(t)}
\leq
&\,\delta\left\{\int_{\d M}|\cm-\cmbar|^{\frac{p(n-1)}{n-2}}\ds_{g(t)}\right\}^{\frac{n-2}{n-1}}\notag
\\
&+\delta^{-\frac{\a}{1-\a}}\left\{\int_{\d M}|\cm-\cmbar|^{p}\ds_{g(t)}\right\}^{\frac{p+2-n}{p+1-n}}\notag
\end{align}
for any $0<\delta<1$. Choosing $\delta$ small, we substitute this last inequality in (\ref{Lemma2.3:2}) and apply again Young's inequality to obtain the estimate (\ref{Lemma2.3:1}).
\ep
\begin{proposition}\label{Propo3.1}
Fix $n-1<p<n$. Then 
$$\lim_{t\to\infty}\int_{\d M}|\cm-\cmbar|^p\ds_{g(t)}=0\,.$$
\end{proposition}
\bp
Since $p>n-1\geq 2$, it follows from Lemma \ref{Lemma2.2} that 
\ba
\frac{d}{dt}\int_{\d M}&(\cm+\s)^{p-1}\dst\notag
\\
&\leq
-(n-p)\int_{\d M}\left\{(\cm+\s)^{p-1}-(\cmbar+\s)^{p-1}\right\}(\cm-\cmbar)\,\dst\,.\notag
\end{align}
One can also check that
$$
\left\{(\cm+\s)^{p-1}-(\cmbar+\s)^{p-1}\right\}(\cm-\cmbar)\geq c|\cm-\cmbar|^p\,.
$$ 
Hence, for $p<n$ we have 
$$
\frac{d}{dt}\int_{\d M}(\cm+\s)^{p-1}\dst
\leq
-c\int_{\d M}|\cm-\cmbar|^p\dst\,.
$$
Integrating, we obtain
$$
\int_0^{\infty}\int_{\d M}|\cm-\cmbar|^p\dst\,dt
\leq 
c^{-1}\int_{\d M}(\cmz+\s)^{p-1}\ds_{g_0}\,,
$$
which implies
$$
\liminf_{t\to\infty}\int_{\d M}|\cm-\cmbar|^p\dst=0\,.
$$
On the other hand, since $p>n-1$, we can apply Lemma \ref{Lemma2.3} to conclude the proof.
\ep
\begin{corollary}\label{Corol3.2}
For any $1<p<n$ we have 
$$\lim_{t\to\infty}\int_{\d M}|\cm-\cminfbar|^p\ds_{g(t)}=0\,.$$
\end{corollary}


\section{The test function}\label{sec:testfunc}
In this section, we construct a test function to be used in our subsequent blow-up analysis. 
Since our construction follows the same steps of \cite{chen}, we only point out the necessary modifications.


\subsection{The auxiliary function $\phi$ and some algebraic preliminaries}\label{subsec:algebraic}
First we fix some notations. If $\e>0$, we define
\begin{equation}\label{eq:def:U}
\U(y)=\left(\frac{\e}{(\e+y_n)^2+|\bar{y}|^2}\right)^{\frac{n-2}{2}}
\:\:\:\:\text{for}\:\:y\in\Rn\,.
\end{equation}
It is well known that the $\U$ satisfy
\begin{align}\label{eq:Ue}
\begin{cases}
\Delta\U=0\,,&\text{in}\:\Rn\,,
\\
\d_n\U+(n-2)\U^{\frac{n}{n-2}}=0\,,&\text{on}\:\d\Rn\,,
\end{cases}
\end{align}
and
\begin{equation}\label{eq:U:Q}
4(n-1)\left(\int_{\d\Rn}\U(y)^{\critbordo}dy\right)^{\frac{1}{n-1}}=\Q\,.
\end{equation}

In this section, $\mathcal{H}$ will denote a symmetric trace-free 2-tensor on $\Rn$ with components $\mathcal{H}_{ab}$, $a,b=1,...,n$, satisfying
\begin{equation}\label{propr:H}
\begin{cases}
\mathcal{H}_{ab}(0)=0\,,&\text{for}\: a,b=1,...,n\,,
\\
\mathcal{H}_{an}(x)=0\,,&\text{forl}\:x\in\Rn,\: a=1,...,n\,,
\\
\d_k\mathcal{H}_{ij}(0)=0\,,&\text{for}\: i,j,k=1,...,n-1\,,
\\
\sum_{j=1}^{n-1}x_j\mathcal{H}_{ij}(x)=0\,,&\text{for}\:x\in\d\Rn,\: i=1,...,n-1\,.
\end{cases}
\end{equation}
We will also assume that those components are of the form
\begin{equation}\label{forma:H}
\mathcal{H}_{ab}(x)=\sum_{|\a|=1}^{d}h_{ab,\a}x^{\a}
\:\:\:\:\:\:\text{for}\:x\in\Rn\,,
\end{equation}
where $d=\left[\frac{n-2}{2}\right]$ and each $\a$ stands for a multi-index. Obviously, the constants $h_{ab,\a}\in\R$ satisfy $h_{an,\a}=0$ for any $\a$, and $h_{ab,\a}=0$ for any $\a\neq (0,...,0,1)$ with $|\a|=1$, where  $a,b=1,...,n$.

Let $\eta:\R\to\R$ be a non-negative smooth function such that
$\eta|_{[0,4/3]}\equiv 1$ and $\eta|_{[5/3,\infty)}\equiv 0$. If $\rho>0$, we define 
\begin{equation}\label{def:eta}
\eta_{\rho}(x)=\eta\left(\frac{|x|}{\rho}\right)
\:\:\:\:\text{for}\:x\in \Rn\,.
\end{equation} 
Notice that $\d_n\eta_{\rho}=0$ on $\d\Rn$.

Let $V=V(\e, \rho, \mathcal{H})$ be the smooth vector field on $\Rn$ obtained in  \cite[Proposition 12]{chen}, which satisfies
\begin{align}\label{eq:V}
\begin{cases}
\sum_{b=1}^{n}\d_b\left\{\U^{\crit}(\eta_{\rho}\mathcal{H}_{ab}-\d_aV_b-\d_bV_a+\frac{2}{n}(\text{div} V)\delta_{ab})\right\}=0\,,&\text{in}\:\Rn\,,
\\
\d_nV_i=V_n=0\,,&\text{on}\:\d\Rn\,,
\end{cases}
\end{align}
for $a=1,...,n$, and $i=1,...,n-1$, and 
\begin{equation}\label{est:V}
|\d^{\b}V(x)|\leq C(n,|\b|)\sum_{i,j=1}^{n-1}\sum_{|\a|=1}^{d}|h_{ij,\a}|(\e+|x|)^{|\a|+1-|\b|}
\end{equation} 
for any multi-index $\b$. Here,
\begin{equation}
\delta_{ab}=
\begin{cases}\notag
1\,,&\text{if}\:a=b\,,
\\
0\,,&\text{if}\:a\neq b\,.
\end{cases}
\end{equation}
We define symmetric trace-free 2-tensors $S$ and $T$ on $\Rn$ by 
$$
S_{ab}=\d_aV_b+\d_bV_a-\frac{2}{n}(\text{div}V)\delta_{ab}\quad\quad\text{and}\quad\quad T=\mathcal{H}-S\,.
$$
Observe that $T_{in}=S_{in}=0$ on $\d\Rn$ for $i=1,...,n-1$.
It follows from (\ref{eq:V}) that $T$ satisfies
$$
\U\d_bT_{ab}+\frac{2n}{n-2}\d_b\U T_{ab}=0\,,
\:\:\:\:\text{in}\:\:B^+_{\rho}(0)\,,
\:\:\:\:\text{for}\:\:a=1,...,n\,.
$$
(Recall that we are adopting the summation convention.) In particular,
\begin{equation}\label{eq:U:T:2}
\frac{n-2}{4(n-1)}\U\d_a\d_bT_{ab}+\d_a(\d_b\U T_{ab})=0\,,
\:\:\:\:\text{in}\:\:B^+_{\rho}(0)\,,
\end{equation}
where we have used the identity 
$\U\d_a\d_b\U-\frac{n}{n-2}\d_a\U\d_b\U=-\frac{1}{n-2}|d\U|^2\delta_{ab}$ in $\Rn$ for all $a,b=1,...,n$.

Next we define the auxiliary function $\phi=\phi_{\e,\rho,\mathcal{H}}$ 
by
\begin{equation}\label{eq:def:phi}
\phi=\d_a\U V_a+\frac{n-2}{2n}\U\text{div}V\,.
\end{equation}
By  \cite[Propositions 1 and 5]{chen} and equation (\ref{eq:U:T:2}),  we have
\begin{align}
\begin{cases}\notag
\Delta \phi=\frac{n-2}{4(n-1)}\U \d_b\d_a \mathcal{H}_{ab}+\d_b(\d_a\U \mathcal{H}_{ab}),&\text{in}\:B^+_{\rho}(0)\,,
\\
\d_n \phi-\frac{n}{n-2}\U^{-1}\d_n\U\phi=-\frac{1}{2(n-1)}\d_n\U S_{nn}, &\text{on}\:\d\Rn\,.
\end{cases}
\end{align}

Observe that if $n=3$ then $d=0$, in which case $\mathcal{H}\equiv 0$ and $\phi\equiv 0$.

\begin{convention} 
In the rest of Section \ref{subsec:algebraic} we will assume that $n\geq 4$.
\end{convention}

We define algebraic Schouten tensor and algebraic Weyl tensor by   
$$A_{ac}=\d_c\d_e\mathcal{H}_{ae}+\d_a\d_e\mathcal{H}_{ce}-\d_e\d_e\mathcal{H}_{ac}-\frac{1}{n-1}\d_e\d_f\mathcal{H}_{ef}\delta_{ac}$$
and
\ba
Z_{abcd}&=\d_b\d_d\mathcal{H}_{ac}-\d_b\d_c\mathcal{H}_{ad}+\d_a\d_c\mathcal{H}_{db}-\d_a\d_d\mathcal{H}_{bc}\notag
\\
&\hspace{0.5cm}+\frac{1}{n-2}\left(A_{ac}\delta_{bd}-A_{ad}\delta_{bc}
+A_{bd}\delta_{ac}-A_{bc}\delta_{db}\right)\,.\notag
\end{align}
We also set
\ba
Q_{ab,c}=\U \d_cT_{ab}&-\frac{2}{n-2}\d_a\U T_{bc}-\frac{2}{n-2}\d_b\U T_{ac}\notag
\\
&+\frac{2}{n-2}\d_d\U T_{ad}\delta_{bc}+\frac{2}{n-2}\d_d\U T_{bd}\delta_{ac}\,.\notag
\end{align}
\begin{lemma}\label{lemma1}
If the tensor $\mathcal{H}$ satisfies
\ba
\begin{cases}\notag
Z_{abcd}=0,&\text{in}\:\Rn\,,
\\
\d_n\mathcal{H}_{ij}=0,&\text{on}\:\d\Rn\,,
\end{cases}
\end{align}
then $\mathcal{H}=0$ in $\Rn$.
\end{lemma}
\bp
Observe that the hypothesis $\d_n\mathcal{H}_{ij}=0$ on $\d\Rn$ implies that $h_{ij,\a}=0$ for $\a=(0,...,0,1)$. In this case, the expression (\ref{forma:H}) can be written as 
$$
\mathcal{H}_{ab}(x)=\sum_{|\a|=2}^{d}h_{ab,\a}x^{\a}\,.
$$
Now the result is just Proposition 2.3 of \cite{brendle-chen}.
\ep
\begin{proposition}\label{propo2}
Set $U_r=B_{r/4}(0,...,0,\frac{3r}{2})\subset \Rn$. Then there exists $C=C(n)>0$ such that
$$
\sum_{i,j=1}^{n-1}\sum_{|\a|=1}^{d}|h_{ij,\a}|^2r^{2|\a|-4+n}
\leq C\int_{U_r}Z_{abcd}Z_{abcd}
+Cr^{-1}\int_{\d 'B^+_{\frac{5r}{3}}(0)\backslash\d 'B^+_{\frac{4r}{3}}(0)}\d_n \mathcal{H}_{ij}\d_n \mathcal{H}_{ij}\,,
$$
for all $r>0$.
\end{proposition}
\bp
If $r=1$,  observe that the square roots of both sides of the inequality are norms in $\mathcal{H}$, due to Lemma \ref{lemma1}. The general case follows by scaling. 
\ep
\begin{lemma}\label{lemma3}
There exists $C=C(n)>0$ such that
$$
\e^{n-2}r^{6-2n}\int_{U_r}Z_{abcd}Z_{abcd}
\leq 
\frac{C}{\theta}\int_{B^+_{2r}(0)\backslash B^+_{r}(0)}Q_{ab,c}Q_{ab,c}
+\theta\e^{n-2}\sum_{i,j=1}^{n-1}\sum_{|\a|=1}^{d}|h_{ij,\a}|^2r^{2|\a|+2-n}
$$
for all $0<\theta<1$ and all $r\geq \e$.
\end{lemma}
\bp
This follows from the third formula in the proof of Proposition 7 in \cite{brendle-chen}, by means of Young's inequality. Observe that, in our calculations, we are using the range $1\leq|\a|\leq d$ in the summation formulas, instead of the range $2\leq|\a|\leq d$ used in \cite{brendle-chen}. 
\ep
\begin{lemma}\label{lemma4}
There exists $C=C(n)>0$ such that
\ba
\e^{n-2}r^{5-2n}&\int_{\d 'B^+_{\frac{5r}{3}}(0)\backslash\d 'B^+_{\frac{4r}{3}}(0)}\d_n \mathcal{H}_{ij}\d_n \mathcal{H}_{ij}\notag
\\
&\leq 
\theta\e^{n-2}\sum_{i,j=1}^{n-1}\sum_{|\a|=1}^{d}|h_{ij,\a}|^2r^{2|\a|+2-n}\notag
\\
&\hspace{0.5cm}+C\int_{\d 'B^+_{2r}(0)\backslash\d 'B^+_{r}(0)}(-\d_n\U)\U(S_{nn})^2
+\frac{C}{\theta}\int_{B^+_{2r}(0)\backslash B^+_{r}(0)}Q_{ij,n}Q_{ij,n}\notag
\end{align}
for all $0<\theta<1$ and all $r\geq \e$. 
\end{lemma}
\bp
Let $\chiup:\R\to\R$ be a non-negative smooth function such that
$\chiup(t)=1$ for $t\in [4/3,5/3]$ and $\chiup(t)=0$ for $t\notin [1,2]$. For $r>0$ and $x\in \Rn$ we define $\chiup_r(x)=\chiup (|x|/r)$.
It follows from $\d_nS_{ij}=\frac{2n}{(n-1)(n-2)}\d_n\U\U^{-1}S_{nn}\delta_{ij}$, on $\d\Rn$, (see the proof of Proposition 5 in \cite{chen}) that, on $\d\Rn$, we have
\ba
\U\d_n\U(S_{nn})^2
&=\frac{(n-1)(n-2)^2}{4n^2}\U^{3}(\d_n\U)^{-1}\d_nS_{ij}\d_nS_{ij}\notag
\\
&=\frac{(n-1)(n-2)^2}{4n^2}\U^{3}(\d_n\U)^{-1}(\d_n\mathcal{H}_{ij}-\d_nT_{ij})(\d_n\mathcal{H}_{ij}-\d_nT_{ij})\notag\,.
\end{align}
Using the fact that $\frac{1}{2}a^2\leq (a-b)^2+b^2$ for any $a,b\in\R$, we obtain
\ba\label{lemma4:1}
\int_{\d\Rn}\U(-\d_n\U)(S_{nn})^2\chiup_r&+\frac{(n-1)(n-2)^2}{4n^2}\int_{\d\Rn}\U^{3}(-\d_n\U)^{-1}\d_nT_{ij}\d_nT_{ij}\chiup_r\notag
\\
&\geq
\frac{(n-1)(n-2)^2}{8n^2}\int_{\d\Rn}\U^{3}(-\d_n\U)^{-1}\d_n\mathcal{H}_{ij}\d_n\mathcal{H}_{ij}\chiup_r\notag
\\
&\geq C^{-1}\e^{n-2}r^{5-2n}\int_{\d'B^+_{\frac{5r}{3}}(0)\backslash \d'B^+_{\frac{4r}{3}}(0)}\d_n\mathcal{H}_{ij}\d_n\mathcal{H}_{ij}\,,
\end{align}
where $C=C(n)>0$.
Since $\U \d_nT_{ij}=Q_{ij,n}$ on $\d\Rn$, integration by parts gives
\ba\label{lemma4:2}
\int_{\d\Rn}&\U^{3}(-\d_n\U)^{-1}\d_nT_{ij}\d_nT_{ij}\chiup_r
=\int_{\d\Rn}\U(-\d_n\U)^{-1}Q_{ij,n}Q_{ij,n}\chiup_r
\\
&=\int_{\Rn}\d_n\big{\{}\U(\d_n\U)^{-1}Q_{ij,n}Q_{ij,n}\chiup_r\big{\}}\notag
\\
&=\int_{\Rn}Q_{ij,n}Q_{ij,n}\chiup_r
-\int_{\Rn}\U(\d_n\U)^{-2}(\d_n\d_n\U)Q_{ij,n}Q_{ij,n}\chiup_r\notag
\\
&\hspace{0.5cm}+\int_{\Rn}2\U(\d_n\U)^{-1}\d_nQ_{ij,n}Q_{ij,n}\chiup_r
+\int_{\Rn}\U(\d_n\U)^{-1}Q_{ij,n}Q_{ij,n}\d_n\chiup_r\,.\notag
\end{align}
Estimating the terms on the right-hand side of (\ref{lemma4:2}) and using H\"{o}lder's and Young's inequalities, we obtain
\ba\label{lemma4:3}
\int_{\d\Rn}\U^{3}(-\d_n\U)^{-1}&\d_nT_{ij}\d_nT_{ij}\chiup_r
\\
&\leq C\left\{\e^{n-2}\sum_{i,j=1}^{n-1}\sum_{|\a|=1}^{d}|h_{ij,\a}|^2r^{2|\a|+2-n}\right\}^{\frac{1}{2}}
\cdot\left\{\int_{B^+_{2r}(0)\backslash B^+_{r}(0)}Q_{ij,n}Q_{ij,n}\right\}^{\frac{1}{2}}\notag
\\
&\leq \theta\e^{n-2}\sum_{i,j=1}^{n-1}\sum_{|\a|=1}^{d}|h_{ij,\a}|^2r^{2|\a|+2-n}
+\frac{C}{\theta}\int_{B^+_{2r}(0)\backslash B^+_{r}(0)}Q_{ij,n}Q_{ij,n}\,.\notag
\end{align}
Now the result follows from the estimates (\ref{lemma4:1}) and (\ref{lemma4:3}).
\ep
\begin{proposition}\label{propo5}
There exists $\l=\l(n)>0$ such that
\ba
\l\e^{n-2}\sum_{i,j=1}^{n-1}&\sum_{|\a|=1}^{d}|h_{ij,\a}|^2\int_{B^+_{\rho}(0)}(\e+|x|)^{2|\a|+2-2n}dx\notag
\\
&\leq \frac{1}{4}\int_{B^+_{\rho}(0)}Q_{ab,c}Q_{ab,c}dx\notag
-\frac{n^2}{2(n-1)(n-2)}\int_{\d 'B^+_{\rho}(0)}\d_n\U\U(S_{nn})^2dx\notag
\end{align}
for all $\rho\geq 2\e$.
\end{proposition}
\bp
It follows from Proposition \ref{propo2}, Lemma \ref{lemma3}, and Lemma \ref{lemma4} that
\ba
\e^{n-2}\sum_{i,j=1}^{n-1}\sum_{|\a|=1}^{d}&|h_{ij,\a}|^2r^{2|\a|+2-n}\notag
\\
\leq
&\,C\int_{\d 'B^+_{2r}(0)\backslash\d 'B^+_{r}(0)}(-\d_n\U)\U(S_{nn})^2
+C\int_{B^+_{2r}(0)\backslash B^+_{r}(0)}Q_{ab,c}Q_{ab,c}\notag
\end{align}
for all $r\geq \e$. Now the assertion follows.
\ep


\subsection{Defining the test function $\bar{u}_{(x,\e)}$ and estimating its Sobolev quotient}\label{sub:sec:deftestfunct}

\begin{definition}
Fix $x_0\in\d M$ and geodesic normal coordinates for $\d M$ centered at $x_0$. 
Let $(x_1,...,x_{n-1})$ be the coordinates of  $x\in\d M$ and $\nu(x)$ be the inward unit vector normal to $\d M$ at $x$. 
For small $x_n\geq 0$, the point $\exp_{x}(x_n\nu(x))\in M$ is said to have {\it{Fermi coordinates}} $(x_1,...,x_n)$ (centered at $x_0$). 
\end{definition}
For small $\rho>0$, the Fermi coordinates centered at $x_0$ define a smooth map $\psi_{x_0}:B_{\rho}^+(0)\subset \Rn\to M$. We will sometimes omit the symbols $\psi_{x_0}$ in order to simplify our notations, identifying $\psi_{x_0}(x)\in M$ with $x\in B^+_{\rho}(0)$. In those coordinates, we have the properties $g_{ab}(0)=\delta_{ab}$ and $g_{nb}(x)=\delta_{nb}$, for any $x\in B^+_{\rho}(0)$ and $a,b=1,...,n$. If we write $g=\exp(h)$, where $\exp$ denotes the matrix exponential, then the symmetric 2-tensor $h$ satisfies the following properties:
\begin{equation}\notag
\begin{cases}
h_{ab}(0)=0\,,&\text{for}\: a,b=1,...,n\,,
\\
h_{an}(x)=0\,,&\text{for}\:x\in B^+_{\rho}(0),\: a=1,...,n\,,
\\
\d_kh_{ij}(0)=0\,,&\text{for}\: i,j,k=1,...,n-1\,,
\\
\sum_{j=1}^{n-1}x_jh_{ij}(x)=0\,,&\text{for}\:x\in \d 'B^+_{\rho}(0),\: i=1,...,n-1\,.
\end{cases}
\end{equation}  
The last two properties follow from the fact that Fermi coordinates are normal on the boundary.

According to  \cite[Proposition 3.1]{marques-weyl}, for each $x_0\in\d M$ we can find a conformal metric $g_{x_0}=f_{x_0}^{\frac{4}{n-2}}g_0$, with $f_{x_0}(x_0)=1$, and Fermi coordinates centered at $x_0$ such that $\text{det}(g_{x_0})(x)=1+O(|x|^{2d+2})$. In particular, 
if we write $g_{x_0}=\exp(h_{x_0})$, we have $\text{tr}(h_{x_0})(x)=O(|x|^{2d+2})$. Moreover, $\cmedia_{g_{x_0}}$, the trace of the second fundamental form of $\d M$, satisfies 
\begin{equation}\label{est:H}
\cmedia_{g_{x_0}}(x)=-\frac{1}{2}g^{ij}\d_ng_{ij}(x)=-\frac{1}{2}\d_n(\log \text{det} (g_{x_0}))(x)=O(|x|^{2d+1})\,.
\end{equation}

Since $M$ is compact, we can fix a small $\rho$ such that $1/2\leq f_{x_0}\leq 3/2$ for any $x_0\in \d M$.
\begin{notation}
In order to simplify our notations, in the coordinates above, we will write $g_{ab}$ and $g^{ab}$ instead of $(g_{x_0})_{ab}$ and $(g_{x_0})^{ab}$ respectively, and $h_{ab}$ instead of $(h_{x_0})_{ab}$.
\end{notation}

In this section, we denote by 
\begin{equation}\notag 
\mathcal{H}_{ab}(x)=\sum_{1\leq |\a|\leq d}h_{ab,\a}x^{\a}
\end{equation}
the Taylor expansion of order $d=\big[\frac{n-2}{2}\big]$ associated with the function $h_{ab}(x)$. Thus, $h_{ab}(x)=\mathcal{H}_{ab}(x)+O(|x|^{d+1})$.  Observe that $\mathcal{H}$ is a symmetric trace-free 2-tensor on $\Rn$, which satisfies the properties (\ref{propr:H}) and has the form (\ref{forma:H}). 
Then we can use the function $\phi=\phi_{\e,\rho,\mathcal{H}}$ 
(see formula (\ref{eq:def:phi})) and the results obtained in Section \ref{subsec:algebraic}. 

Let us assume $Q(M,\d M)>0$.  Recall the definitions of $\U$ in (\ref{eq:def:U}), $\eta_{\rho}$ in (\ref{def:eta}), and $\cminfbar$ in (\ref{eq:def:cminfbar}). Define
\ba\label{def:test:func}
\ubar(x)
=&\left(\frac{2(n-1)}{\cminfbar}\right)^{\frac{n-2}{2}}\eta_{\rho}(\psi_{x_0}^{-1}(x))\big(\U(\psi_{x_0}^{-1}(x))+\phi(\psi_{x_0}^{-1}(x))\big)
\\
&\hspace{0.1cm}+\left(\frac{2(n-1)}{\cminfbar}\right)^{\frac{n-2}{2}}\e^{\frac{n-2}{2}}\big(1-\eta_{\rho}(\psi_{x_0}^{-1}(x))\big)G_{x_0}(x)\,,\notag
\end{align}
if $x\in\psi_{x_0}(B^+_{2\rho}(0))$, and 
$$
\ubar(x)=G_{x_0}(x)\,,
\:\:\:\:\text{otherwise}\,. 
$$ 
Here, $G_{x_0}$ is the Green's function of the conformal Laplacian $L_{g_{x_0}}=\Delta_{g_{x_0}}-\frac{n-2}{4(n-1)}R_{g_{x_0}}$, with pole at $x_0\in\d M$, satisfying the boundary condition 
\begin{equation}\label{eq:G:bordo}
\frac{\d}{\d\eta_{g_{x_0}}}G_{x_0}-\frac{n-2}{2(n-1)}\cmedia_{g_{x_0}}G_{x_0}=0
\end{equation}
and the normalization $\lim_{|y|\to 0}|y|^{n-2}G_{x_0}(\psi_{x_0}(y))=1$. This function, obtained in Proposition \ref{green:point}, satisfies 
\ba\label{estim:G}
|G_{x_0}(\psi_{x_0}(y))-|y|^{2-n}|
&\leq
C\sum_{i,j=1}^{n-1}\sum_{|\a|=1}^{d}|h_{ij,\a}||y|^{|\a|+2-n}+C|y|^{d+3-n}\,,
\\
\left|\frac{\d}{\d y_b}(G_{x_0}(\psi_{x_0}(y))-|y|^{2-n})\right|
&\leq
C\sum_{i,j=1}^{n-1}\sum_{|\a|=1}^{d}|h_{ij,\a}||y|^{|\a|+1-n}+C|y|^{d+2-n}\,,
\end{align}
for all $b=1,...,n$.

By the estimate (\ref{est:V}), $\phi$ satisfies $|\phi(y)|\leq C\e^{\frac{n-2}{2}}\sum_{i,j=1}^{n-1}\sum_{|\a|=1}^{d}|h_{ij,\a}|(\e+|y|)^{|\a|+2-n}$, for all $y\in\Rn$, and 
$$
\left|\d_n\phi(y)+n\U^{\frac{2}{n-2}}\phi(y)\right|
\leq C\e^{\frac{n}{2}}\sum_{i,j=1}^{n-1}\sum_{|\a|=1}^{d}|h_{ij,\a}|(\e+|y|)^{|\a|-n}\,,
$$
for all $y\in\d\Rn$.

We define the test function 
\begin{equation}\label{def:test:func:u}
\func:u=f_{x_0}\ubar\,.
\end{equation}

Our main result in this section is the following estimate for the energy of $\func:u$:
\begin{proposition}\label{Propo:energy:test}
Under the hypothesis of Theorem \ref{main:thm}, there exists $\e_0>0$, depending only on $(M,g_0)$, such that 
\ba
&\frac{\int_M\frac{4(n-1)}{n-2}|d \func:u|_{g_0}^2\dv_{g_0}+\int_{\d M}2\cmedia_{g_0}\func:u^2\ds_{g_0}}{\left(\int_{\d M}\func:u^{\frac{2(n-1)}{n-2}}\ds_{g_0}\right)^{\frac{n-2}{n-1}}}\notag
\\
&\hspace{0.5cm}=\frac{\int_M\left\{\frac{4(n-1)}{n-2}|d \ubar|_{g_{x_0}}^2
+R_{g_{x_0}}\ubar^2\right\}\dv_{g_{x_0}}
+\int_{\d M}2\cmedia_{g_{x_0}}\ubar^2\ds_{g_{x_0}}}{\left(\int_{\d M}\ubar^{\frac{2(n-1)}{n-2}}\ds_{g_{x_0}}\right)^{\frac{n-2}{n-1}}}\notag
\\
&\hspace{0.5cm}\leq \Q\notag
\end{align}
for all $x_0\in\d M$ and $\e\in(0,\e_0)$.
\end{proposition}
\begin{convention} 
In the rest of Section \ref{sec:testfunc}, we will use the normalization $\cminfbar=2(n-1)$, without loss of generality.
\end{convention}

Let $\l$ be the constant obtained in Proposition \ref{propo5}.
\begin{proposition}\label{propo6}
There exist $C=C(n,g_0)$ and $\rho_0=\rho_0(n,g_0)$ such that
\ba
&\int_{B^+_{\rho}(0)}\left\{\frac{4(n-1)}{n-2}|d(\U+\phi)|_{g_{x_0}}^2+R_{g_{x_0}}(\U+\phi)^2\right\}dx
+\int_{\d 'B^+_{\rho}(0)}2\cmedia_{g_{x_0}}(\U+\phi)^2dx\notag
\\
&\hspace{1cm}\leq
4(n-1)\int_{\d 'B^+_{\rho}(0)}\U^{\frac{2}{n-2}}
\left\{\U^2+2\U\phi+\frac{n}{n-2}\phi^2-\frac{n-2}{8(n-1)^2}\U^2|S_{nn}|^2\right\}dx\notag
\\
&\hspace{1.5cm}+\int_{\d^+B^+_{\rho}(0)}\left\{\frac{4(n-1)}{n-2}\U\d_a\U+\U^2\d_bh_{ab}-\d_b\U^2h_{ab}\right\}
\frac{x_a}{|x|}\ds_{\rho}\notag
\\
&\hspace{1.5cm}-\frac{\l}{2}\sum_{i,j=1}^{n-1}\sum_{|\a|=1}^{d}|h_{ij,\a}|^2\e^{n-2}\int_{B^+_{\rho}(0)}(\e+|x|)^{2|\a|+2-2n}dx\notag
\\
&\hspace{1.5cm}+C\sum_{i,j=1}^{n-1}\sum_{|\a|=1}^{d}|h_{ij,\a}|\e^{n-2}\rho^{|\a|+2-n}
+C\e^{n-2}\rho^{2d+4-n}\notag
\end{align}
for all $0<2\e\leq \rho\leq\rho_0$.
\end{proposition}
\bp
It follows from  \cite[Proposition 11]{brendle-invent} that the scalar curvature satisfies
\begin{equation}\label{estim:R:1}
|R_{g_{x_0}}-\d_a\d_b\mathcal{H}_{ab}|\leq C\sum_{i,j=1}^{n-1}\sum_{|\a|=1}^{d}|h_{ij,\a}||x|^{|\a|-1}+C|x|^{d-1}
\end{equation}
and
\begin{align}\label{estim:R:2}
&\big|R_{g_{x_0}}-\d_a\d_bh_{ab}+\d_b(\mathcal{H}_{ab}\d_c\mathcal{H}_{ac})-\frac{1}{2}\d_b\mathcal{H}_{ab}\d_c\mathcal{H}_{ac}+\frac{1}{4}\d_c\mathcal{H}_{ab}\d_c\mathcal{H}_{ab}\big|
\\
&\hspace{1cm}\leq C\sum_{i,j=1}^{n-1}\sum_{|\a|=1}^{d}|h_{ij,\a}|^2|x|^{2|\a|-1}
+C\sum_{i,j=1}^{n-1}\sum_{|\a|=1}^{d}|h_{ij,\a}||x|^{|\a|+d-1}+C|x|^{2d}\,.\notag
\end{align}
We point out that, although these estimates are a little weaker than those in \cite[Proposition 3]{chen}, they are enough to prove our result.   

Following the steps  in \cite[Proposition 7]{chen} we obtain
\ba
\int_{B^+_{\rho}(0)}&\left\{\frac{4(n-1)}{n-2}|d(\U+\phi)|_{g_{x_0}}^2+R_{g_{x_0}}(\U+\phi)^2\right\}dx
+\int_{\d 'B^+_{\rho}(0)}2\cmedia_{g_{x_0}}(\U+\phi)^2dx\notag
\\
&\hspace{0.5cm}\leq
-\frac{4(n-1)}{n-2}\int_{\d 'B^+_{\rho}(0)}
\left\{\U\d_n\U+2\d_n\U\phi+\frac{n}{n-2}\U^{-1}\d_n\U\phi^2\right\}dx\notag
\\
&\hspace{1cm}+\frac{n+2}{2(n-2)}\int_{\d 'B^+_{\rho}(0)}\U\d_n\U(S_{nn})^2dx
-\frac{1}{4}\int_{B^+_{\rho}(0)}Q_{ab,c}Q_{ab,c}dx\notag
\\
&\hspace{1cm}+\int_{\d^+B^+_{\rho}(0)}\left\{\frac{4(n-1)}{n-2}\U\d_a\U+\U^2\d_bh_{ab}-\d_b\U^2h_{ab}\right\}
\frac{x_a}{|x|}\ds_{\rho}\notag
\\
&\hspace{1cm}+\frac{\l}{2}\sum_{i,j=1}^{n-1}\sum_{|\a|=1}^{d}|h_{ij,\a}|^2\e^{n-2}\int_{B^+_{\rho}(0)}(\e+|x|)^{2|\a|+2-2n}dx\notag
\\
&\hspace{1cm}+C\sum_{i,j=1}^{n-1}\sum_{|\a|=1}^{d}|h_{ij,\a}|\e^{n-2}\rho^{|\a|+2-n}
+C\e^{n-2}\rho^{2d+4-n}\,.\notag
\end{align}
Now the assertion follows from Proposition \ref{propo5} and the second equation of (\ref{eq:Ue}).
\ep

As in \cite{chen} (see also \cite{brendle-invent, brendle-chen}), we define the flux integral
\ba\label{def:I}
\mathcal{I}(x_0,\rho)
=\frac{4(n-1)}{n-2}&\int_{\d^+B^+_{\rho}(0)}(|x|^{2-n}\d_aG_{x_0}-\d_a|x|^{2-n}G_{x_0})\frac{x_a}{|x|}\ds_{\rho}
\\
-&\int_{\d^+B^+_{\rho}(0)}|x|^{2-2n}(|x|^{2}\d_bh_{ab}-2nx_bh_{ab})\frac{x_a}{|x|}\ds_{\rho}\,,\notag
\end{align}
for $\rho>0$ sufficiently small. 
\begin{proposition}\label{propo7}
There exists $\rho_0=\rho_0(n,g_0)$ such that
\ba
&\int_M\left\{\frac{4(n-1)}{n-2}|d\ubar|_{g_{x_0}}^2+R_{g_{x_0}}\ubar^2\right\}\dv_{g_{x_0}}
+\int_{\d M}2\cmedia_{g_{x_0}}\ubar^2\ds_{g_{x_0}}\notag
\\
&\hspace{2cm}\leq
\Q\left\{\int_{\d M}\ubar^{\critbordo}\ds_{g_{x_0}}\right\}^{\frac{n-2}{n-1}}
-\e^{n-2}\mathcal{I}(x_0,\rho)\notag
\\
&\hspace{2.5cm}-\frac{\l}{4}\sum_{i,j=1}^{n-1}\sum_{|\a|=1}^{d}|h_{ij,\a}|^2\e^{n-2}
\int_{B^+_{\rho}(0)}(\e+|x|)^{2|\a|+2-2n}dx\notag
\\
&\hspace{2.5cm}+C\sum_{i,j=1}^{n-1}\sum_{|\a|=1}^{d}|h_{ij,\a}|\e^{n-2}\rho^{|\a|+2-n}
+C\e^{n-2}\rho^{2d+4-n}
+C\rho^{1-n}\e^{n-1}\notag
\end{align}
for all $0<2\e\leq \rho\leq\rho_0$. 
\end{proposition}
\bp
Once we have proved Proposition \ref{propo6}, our proof  is analogous to the one in \cite[Proposition 9]{chen}. A necessary step is the estimate
\ba\label{est:propo8}
4(n-1)&\int_{\d 'B^+_{\rho}(0)}\U^{\frac{2}{n-2}}\left(\U^2+2\U\phi+\frac{n}{n-2}\phi^2-\frac{n-2}{8(n-1)^2}\U^2S_{nn}^2\right)dx
\\
&\leq \Q\left(\int_{\d 'B^+_{\rho}(0)}(\U+\phi)^{\critbordo}dx\right)^{\frac{n-2}{n-1}}
+\sum_{i,j=1}^{n-1}\sum_{|\a|=1}^{d}|h_{ij,\a}|\rho^{|\a|+1-n}\e^{n-1}\notag
\\
&\hspace{1cm}+C\sum_{i,j=1}^{n-1}\sum_{|\a|=1}^{d}|h_{ij,\a}|^2\e^{n-1}\rho\int_{\d 'B^+_{\rho}(0)}(\e+|x|)^{2|\a|+2-2n}dx\notag
\end{align}
for all $0<2\e\leq \rho\leq\rho_0$ and $\rho_0$ sufficiently small. This inequality is slightly different from the one in \cite[Proposition 8]{chen}, since the Taylor expansion (\ref{forma:H}) for $\mathcal{H}_{ab}$ includes terms of order $|\a|=1$. However, the estimate (\ref{est:propo8}) is enough to prove our assertion. Also  observe that we are assuming a different boundary condition for the Green's function $G_{x_0}$ (see (\ref{eq:G:bordo})) which differ from the one in \cite{chen} by the term $\frac{n-2}{2(n-1)}\cmedia_{g_{x_0}}G_{x_0}$. However, this term is easily estimated using (\ref{est:H}) and (\ref{estim:G}).
\ep
\begin{corollary}\label{corol8}
There exist $\rho_0,\: \theta,\: C_0>0$, depending only on $(M,g_0)$, such that
\ba
&\int_M\left\{\frac{4(n-1)}{n-2}|d\ubar|_{g_{x_0}}^2+R_{g_{x_0}}\ubar^2\right\}\dv_{g_{x_0}}
+\int_{\d M}2\cmedia_{g_{x_0}}\ubar^2\ds_{g_{x_0}}\notag
\\
&\hspace{2cm}\leq
\Q\left\{\int_{\d M}\ubar^{\critbordo}\ds_{g_{x_0}}\right\}^{\frac{n-2}{n-1}}
-\e^{n-2}\mathcal{I}(x_0,\rho)\notag
\\
&\hspace{2.5cm}-\theta\e^{n-2}\int_{B^+_{\rho}(0)}|W_{g_0}(x)|^2(\e+|x|)^{6-2n}dx\notag
\\
&\hspace{2.5cm}-\theta\e^{n-2}\int_{\d 'B^+_{\rho}(0)}|\pi_{g_0}(x)|^2(\e+|x|)^{5-2n}dx\notag
\\
&\hspace{2.5cm}+C_0\e^{n-2}\rho^{2d+4-n}
+C_0\left(\frac{\e}{\rho}\right)^{n-2}\frac{1}{\log(\rho/\e)}\notag
\end{align}
for all $0<2\e\leq \rho\leq\rho_0$. Here, we denote by $W_{g_0}$ the Weyl tensor of $(M,g_0)$ and by $\pi_{g_0}$ the trace-free 2nd fundamental form of $\d M$.
\end{corollary}
\bp
By Young's inequality, given $C>0$ there exists $C'>0$ such that
$$
C|h_{ij,\a}|\e^{n-2}\rho^{|\a|+2-n}
\leq
\frac{\l}{8}|h_{ij,\a}|^2\e^{n-2}\int_{B^+_{\rho}}(\e+|x|)^{2|\a|+2-2n}dx
+C'\left(\frac{\e}{\rho}\right)^{2n-4-2|\a|}\,,
$$
for $|\a|<\frac{n-2}{2}$, and 
$$
C|h_{ij,\a}|\e^{n-2}\rho^{|\a|+2-n}
\leq
\frac{\l}{8}|h_{ij,\a}|^2\e^{n-2}\int_{B^+_{\rho}}(\e+|x|)^{2|\a|+2-2n}dx
+C'\left(\frac{\e}{\rho}\right)^{n-2}\frac{1}{\log(\rho/\e)}\,,
$$
for $|\a|=\frac{n-2}{2}$. Then, according to Proposition \ref{propo7}, we have
\ba\label{corol8:1}
&\int_M\left\{\frac{4(n-1)}{n-2}|d\ubar|_{g_{x_0}}^2+R_{g_{x_0}}\ubar^2\right\}\dv_{g_{x_0}}
+\int_{\d M}2\cmedia_{g_{x_0}}\ubar^2\ds_{g_{x_0}}
\\
&\hspace{2cm}\leq
\Q\left\{\int_{\d M}\ubar^{\critbordo}\ds_{g_{x_0}}\right\}^{\frac{n-2}{n-1}}
-\e^{n-2}\mathcal{I}(x_0,\rho)\notag
\\
&\hspace{2.5cm}-\frac{\l}{8}\sum_{i,j=1}^{n-1}\sum_{|\a|=1}^{d}|h_{ij,\a}|^2\e^{n-2}\int_{B^+_{\rho}}(\e+|x|)^{2|\a|+2-2n}dx\notag
\\
&\hspace{2.5cm}+C\e^{n-2}\rho^{2d+4-n}
+C\left(\frac{\e}{\rho}\right)^{n-2}\frac{1}{\log(\rho/\e)}\notag
\end{align}

On the other hand, we have the pointwise estimates
$$
|W_{g_0}(x)|=|W_{g_{x_0}}(x)|
\leq C|\d^2h(x)|+C|\d h(x)|
\leq C\sum_{1\leq|\a|\leq d}\sum_{i,j=1}^{n-1}|h_{ij,\a}||x|^{|\a|-2}+C|x|^{d-1}
$$
and
$$
|\pi_{g_0}(x)|=|\pi_{g_{x_0}}(x)|
\leq C|\d h(x)|
\leq C\sum_{1\leq|\a|\leq d}\sum_{i,j=1}^{n-1}|h_{ij,\a}||x|^{|\a|-1}+C|x|^{d}\,.
$$
Hence,
\ba\label{corol8:2}
\int_{B^+_{\rho}(0)}&|W_{g_{x_0}}(x)|^2(\e+|x|)^{6-2n}dx
+\int_{\d 'B^+_{\rho}(0)}|\pi_{g_{x_0}}(x)|^2(\e+|x|)^{5-2n}dx
\\
&\leq
C\sum_{1\leq|\a|\leq d}\sum_{i,j=1}^{n-1}|h_{ij,\a}|^2\int_{B^+_{\rho}(0)}(\e+|x|)^{2|\a|+2-2n}dx\notag
\\
&\hspace{0.5cm}+C\sum_{1\leq|\a|\leq d}\sum_{i,j=1}^{n-1}|h_{ij,\a}|^2\int_{\d 'B^+_{\rho}(0)}(\e+|x|)^{2|\a|+3-2n}dx
+C\rho^{2d+4-n}\notag
\\
&\leq
C'\sum_{1\leq|\a|\leq d}\sum_{i,j=1}^{n-1}|h_{ij,\a}|^2\int_{B^+_{\rho}(0)}(\e+|x|)^{2|\a|+2-2n}dx+C\rho^{2d+4-n}\,.\notag
\end{align}

Now the result follows from the estimates (\ref{corol8:1}) and (\ref{corol8:2}).
\ep
Recall that we denote by $\mathcal{Z}$ the set of all points $x_0\in\d M$ such that
$$
\limsup_{x\to x_0}d_{g_0}(x,x_0)^{2-d}|W_{g_0}(x)|
=\limsup_{x\to x_0}d_{g_0}(x,x_0)^{1-d}|\pi_{g_0}(x)|=0\,.
$$
\begin{proposition}\label{propo18}
The functions $\mathcal{I}(x_0,\rho)$ converge uniformly to a continuous function $I:\mathcal{Z}\to\R$ as $\rho\to 0$.
\end{proposition}
\bp
We will prove that there exists $C>0$ such that
\begin{equation}\label{propo18:0}
\sup_{x_0\in\mathcal{Z}}|\mathcal{I}(x_0,\rho)-\mathcal{I}(x_0,\tilde{\rho})|\leq C\rho^{2d+4-n}\,,
\:\:\:\:\text{for all}\:\:0<\tilde{\rho}<\rho\,.
\end{equation}

Our proof follows the same steps of \cite[Proposition 18]{brendle-invent}. However, our computations are slightly different because here we cannot assume $x^a h_{ab}(x)=0$, since this property is a consequence of the use of normal coordinates in  \cite{brendle-invent}.

Fix $x_0\in\mathcal{Z}$ and consider Fermi coordinates $\psi_{x_0}:B^+_{2\rho}(0)\to M$ as in the beginning of Section \ref{sub:sec:deftestfunct}. We will write $B^+_{\rho}=B^+_{\rho}(0)$ and $B^+_{\tilde{\rho}}=B^+_{\tilde{\rho}}(0)$ for short. Integrating by parts, we see that
\ba\label{propo18:1}
\mathcal{I}(x_0,\rho)-\mathcal{I}(x_0,\tilde{\rho})&=
\frac{4(n-1)}{n-2}\int_{B^+_{\rho}\backslash B^+_{\tilde{\rho}}}|x|^{2-n}\Delta G_{x_0} dx
\\
&\hspace{1cm}-\int_{\d^+B^+_{\rho}}\left\{|x|^{3-2n}x_i\d_jh_{ij}-2n|x|^{1-2n}x_ix_jh_{ij})\right\}\ds_{\rho}\notag
\\
&\hspace{1cm}+\int_{\d^+B^+_{\tilde{\rho}}}\left\{|x|^{3-2n}x_i\d_jh_{ij}-2n|x|^{1-2n}x_ix_jh_{ij})\right\}\ds_{\tilde{\rho}}\notag
\\
&\hspace{1cm}+O(\rho^{2d+4-n})\,.\notag
\end{align} 
Here, $\Delta$ stands for the Euclidean Laplacian and we have used (\ref{est:H}).

Since $x_0\in\mathcal{Z}$, we have $g_{ij}(x)=\delta_{ij}+O(|x|^{d+1})$ and $G_{x_0}(x)=|x|^{2-n}+O(|x|^{d+3-n})$. 
Then
\ba\label{propo18:2}
\int_{B^+_{\rho}\backslash B^+_{\tilde{\rho}}}|x|^{2-n}\Delta G_{x_0}\, dx
&=-\int_{B^+_{\rho}\backslash B^+_{\tilde{\rho}}}|x|^{2-n}(L_{g_{x_0}}-\Delta)|x|^{2-n} dx
\\
&\hspace{1cm}-\int_{B^+_{\rho}\backslash B^+_{\tilde{\rho}}}|x|^{2-n}(L_{g_{x_0}}-\Delta) (G_{x_0}-|x|^{2-n})\, dx\notag
\\
&=-\int_{B^+_{\rho}\backslash B^+_{\tilde{\rho}}}|x|^{2-n}(L_{g_{x_0}}-\Delta)|x|^{2-n} dx
+O(\rho^{2d+4-n})\,.\notag
\end{align}

Using  $g^{ij}(x)=\delta_{ij}-h_{ij}(x)+O(|x|^{2d+2})$, $\text{tr} (h) (x)=O(|x|^{2d+2})$, $\det(g_{x_0})(x)=1+O(|x|^{2d+2})$, and (\ref{estim:R:2}), we obtain
\ba
(L_{g_{x_0}}-\Delta)|x|^{2-n}
&=g^{ij}\d_i\d_j|x|^{2-n}+\d_ig^{ij}\d_j|x|^{2-n}
+\frac{1}{2}\frac{\d_i\det(g_{x_0})}{\det(g_{x_0})}g^{ij}\d_j|x|^{2-n}
\\
&\hspace{1cm}-\frac{n-2}{4(n-1)}R_{g_{x_0}}|x|^{2-n}\notag
\\
&=-n(n-2)|x|^{-2-n}x_ix_jh_{ij}+(n-2)|x|^{-n}x_j\d_ih_{ij}\notag
\\
&\hspace{1cm}-\frac{n-2}{4(n-1)}|x|^{2-n}\d_i\d_jh_{ij}+O(|x|^{2+2d-n})\,.\notag
\end{align}
Hence,
\ba\label{propo18:3}
-\int_{B^+_{\rho}\backslash B^+_{\tilde{\rho}}}&|x|^{2-n}(L_{g_{x_0}}-\Delta)|x|^{2-n}dx
\\
&=n(n-2)\int_{B^+_{\rho}\backslash B^+_{\tilde{\rho}}}|x|^{-2n}x_ix_jh_{ij}\,dx
-(n-2)\int_{B^+_{\rho}\backslash B^+_{\tilde{\rho}}}|x|^{2-2n}x_j\d_ih_{ij}\,dx\notag
\\
&\hspace{1cm}+\frac{n-2}{4(n-1)}\int_{B^+_{\rho}\backslash B^+_{\tilde{\rho}}}|x|^{4-2n}\d_i\d_jh_{ij}\,dx
+O(\rho^{4+2d-n})\,.\notag
\end{align}
Integrating by parts, we obtain
\ba\label{propo18:4}
&-(n-2)\int_{B^+_{\rho}\backslash B^+_{\tilde{\rho}}}|x|^{2-2n}x_j\d_ih_{ij}\,dx
\\
&=-2(n-1)(n-2)\int_{B^+_{\rho}\backslash B^+_{\tilde{\rho}}}|x|^{-2n}x_ix_jh_{ij}\,dx
+(n-2)\int_{B^+_{\rho}\backslash B^+_{\tilde{\rho}}}|x|^{2-2n}\delta_{ij}h_{ij}\,dx\notag
\\
&\hspace{0.3cm}-(n-2)\int_{\d^+B^+_{\rho}}|x|^{1-2n}x_ix_jh_{ij}\,\ds_{\rho}
+(n-2)\int_{\d^+B^+_{\tilde{\rho}}}|x|^{1-2n}x_ix_jh_{ij}\,\ds_{\tilde{\rho}}\notag
\end{align}
and
\ba\label{propo18:5}
&\frac{n-2}{4(n-1)}\int_{B^+_{\rho}\backslash B^+_{\tilde{\rho}}}|x|^{4-2n}\d_i\d_jh_{ij}\,dx
\\
&=\frac{(n-2)^2}{2(n-1)}\int_{B^+_{\rho}\backslash B^+_{\tilde{\rho}}}|x|^{2-2n}x_i\d_jh_{ij}\,dx\notag
\\
&\hspace{0.5cm}+\frac{n-2}{4(n-1)}\int_{\d^+B^+_{\rho}}|x|^{3-2n}x_i\d_jh_{ij}\,\ds_{\rho}
-\frac{n-2}{4(n-1)}\int_{\d^+B^+_{\tilde{\rho}}}|x|^{3-2n}x_i\d_jh_{ij}\,\ds_{\tilde{\rho}}\notag
\\
&=(n-2)^2\int_{B^+_{\rho}\backslash B^+_{\tilde{\rho}}}|x|^{-2n}x_ix_jh_{ij}\,dx
-\frac{(n-2)^2}{2(n-1)}\int_{B^+_{\rho}\backslash B^+_{\tilde{\rho}}}|x|^{2-2n}\delta_{ij}h_{ij}\,dx\notag
\\
&\hspace{0.5cm}+\frac{(n-2)^2}{2(n-1)}\int_{\d^+B^+_{\rho}}|x|^{1-2n}x_ix_jh_{ij}\,\ds_{\rho}
-\frac{(n-2)^2}{2(n-1)}\int_{\d^+B^+_{\tilde{\rho}}}|x|^{1-2n}x_ix_jh_{ij}\,\ds_{\tilde{\rho}}\notag
\\
&\hspace{0.5cm}+\frac{n-2}{4(n-1)}\int_{\d^+B^+_{\rho}}|x|^{3-2n}x_i\d_jh_{ij}\,\ds_{\rho}
-\frac{n-2}{4(n-1)}\int_{\d^+B^+_{\tilde{\rho}}}|x|^{3-2n}x_i\d_jh_{ij}\,\ds_{\tilde{\rho}}\,.\notag
\end{align}

Substituting (\ref{propo18:4}) and (\ref{propo18:5}) in (\ref{propo18:3}), the coefficients of $\int_{B^+_{\rho}\backslash B^+_{\tilde{\rho}}}|x|^{-2n}x_ix_jh_{ij}\,dx$ cancel out and we obtain
\ba
-\int_{B^+_{\rho}\backslash B^+_{\tilde{\rho}}}&|x|^{2-n}(L_{g_{x_0}}-\Delta)|x|^{2-n}dx\notag
\\
=&-\frac{n(n-2)}{2(n-1)}\left\{\int_{\d^+B^+_{\rho}}|x|^{1-2n}x_ix_jh_{ij}\,\ds_{\rho}
-\int_{\d^+B^+_{\tilde{\rho}}}|x|^{1-2n}x_ix_jh_{ij}\,\ds_{\tilde{\rho}}\right\}\notag
\\
&+\frac{n-2}{4(n-1)}\left\{\int_{\d^+B^+_{\rho}}|x|^{3-2n}x_i\d_jh_{ij}\,\ds_{\rho}
-\int_{\d^+B^+_{\tilde{\rho}}}|x|^{3-2n}x_i\d_jh_{ij}\,\ds_{\tilde{\rho}}\right\}\notag
\\
&+O(\rho^{2d+4-n})\,,\notag
\end{align}
where we used again that $\text{tr}(h)(x)=O(|x|^{2d+2})$\,. Hence, we have
\ba\label{propo18:6}
\frac{4(n-1)}{n-2}&\int_{B^+_{\rho}\backslash B^+_{\tilde{\rho}}}|x|^{2-n}\Delta G_{x_0}\,dx
\\
=&-2n\left\{\int_{\d^+B^+_{\rho}}|x|^{1-2n}x_ix_jh_{ij}\,\ds_{\rho}
-\int_{\d^+B^+_{\tilde{\rho}}}|x|^{1-2n}x_ix_jh_{ij}\,\ds_{\tilde{\rho}}\right\}\notag
\\
&+\left\{\int_{\d^+B^+_{\rho}}|x|^{3-2n}x_i\d_jh_{ij}\,\ds_{\rho}
-\int_{\d^+B^+_{\tilde{\rho}}}|x|^{3-2n}x_i\d_jh_{ij}\,\ds_{\tilde{\rho}}\right\}\notag
\\
&+O(\rho^{2d+4-n})\,.\notag
\end{align}
 Now the assertion follows from (\ref{propo18:1}) and (\ref{propo18:6}). 
\ep

The following proposition relates $\mathcal{I}(x_0)$ with the mass defined by (\ref{def:mass}).
\begin{proposition}\label{propo19} 
Let $x_0\in\mathcal{Z}$ and consider inverted coordinates $y=x/|x|^2$, where $x=(x_1,...,x_n)$ are Fermi coordinates centered at $x_0$. If we define the metric $\bar{g}=G_{x_0}^{\frac{4}{n-2}}g_{x_0}$ on $M\backslash\{x_0\}$, then the following statements hold:

{\bf{(i)}} $(M\backslash\{x_0\},\bar{g})$ is an asymptotically flat manifold with order  $p>\frac{n-2}{2}$ (in the sense of Definition \ref{def:asym}), and  satisfies $R_{\bar{g}}\equiv 0$ and $\cmedia_{\bar{g}}\equiv 0$.

{\bf{(ii)}} We have 
$$
\mathcal{I}(x_0)=\lim_{R\to\infty}
\left\{\int_{\d^+B^+_R(0)}\frac{y_a}{|y|}\frac{\d}{\d y_b}
\bar{g}\left(\frac{\d}{\d y_a},\frac{\d}{\d y_b}\right)\ds_{R}
-\int_{\d^+B^+_R(0)}\frac{y_a}{|y|}\frac{\d}{\d y_a}
\bar{g}\left(\frac{\d}{\d y_b},\frac{\d}{\d y_b}\right)\ds_{R}
\right\}\,.
$$

In particular, $\mathcal{I}(x_0)$ is the mass $m(\bar{g})$ of $(M\backslash\{x_0\},\bar{g})$.
\end{proposition}
\bp
The item (i) follows from the fact that $\bar{g}\left(\frac{\d}{\d y_a},\frac{\d}{\d y_b}\right)=\delta_{ab}+O(|y|^{-d-1})$ and the definition of $G_{x_0}$.
In order to prove (ii), we can mimic the proof in \cite[Proposition 4.3]{brendle-chen} to obtain
\ba
\int_{\d^+B^+_{\rho^{-1}}(0)}\frac{y_a}{|y|}\frac{\d}{\d y_b}
\bar{g}\left(\frac{\d}{\d y_a},\frac{\d}{\d y_b}\right)\ds_{\rho^{-1}}
&-\int_{\d^+B^+_{\rho^{-1}}(0)}\frac{y_a}{|y|}\frac{\d}{\d y_a}
\bar{g}\left(\frac{\d}{\d y_b},\frac{\d}{\d y_b}\right)\ds_{\rho^{-1}}\notag
\\
&=\mathcal{I}(x_0,\rho)+O(\rho^{2d+4-n})\,,\notag
\end{align}
where we used (\ref{est:H}).

The last statement of Proposition \ref{propo19} follows from the fact that 
$$
\bar{g}\left(\frac{\d}{\d y_i}, \frac{\d}{\d y_n}\right)=0\,,
\:\:\:\:\:\text{for}\: i=1,...,n-1\,,
\:\:\:\text{if}\:y_n=0\,.
$$
\ep

We are now able to prove Proposition \ref{Propo:energy:test}.
\bp[Proof of Proposition \ref{Propo:energy:test}]
Assume that  $\mathcal{I}(x_0)>0$ for all $x_0\in\mathcal{Z}$. Since $\mathcal{Z}\subset \d M$ is compact and $\mathcal{I}$ is continuous on $\mathcal{Z}$, we know that $\inf_{x_0\in\mathcal{Z}}\mathcal{I}(x_0)>0$. 

By Proposition \ref{propo18},  $\sup_{x_0\in\mathcal{Z}}|\mathcal{I}(x_0,\rho)-\mathcal{I}(x_0)|\to 0$ as $\rho\to 0$. Hence, we can find $\rho\in (0,\rho_0]$ such that
$$
\inf_{x_0\in\mathcal{Z}}\mathcal{I}(x_0,\rho)>C_0\rho^{2d+4-n}\,.
$$
Here, $\rho_0$ and $C_0$ are the constants appearing in Corollary \ref{corol8}.
By continuity, there exists an open subset $\Omega\subset \d M$, containing $\mathcal{Z}$, such that
\begin{equation}\label{Propo:energy:test:1}
\inf_{x_0\in\Omega}\mathcal{I}(x_0,\rho)>C_0\rho^{2d+4-n}\,.
\end{equation}
(If $\mathcal{Z}=\emptyset$ we set $\Omega=\emptyset$.)

Observe that $\mathcal{Z}=\d M$ if $n=3$. If $n\geq 4$, we will prove that 
\ba\label{Propo:energy:test:2}
\int_{B_{\rho}(x_0)}&|W_g(x)|^2d_{g_0}(x,x_0)^{6-2n}\dv_{g_0}
\\
&+\int_{D_{\rho}(x_0)}|\pi_g(x)|^2d_{g_0}(x,x_0)^{5-2n}\ds_{g_0}
=\infty\,,
\:\:\:\:\:\text{for all}\:x_0\in \d M\backslash \Omega\,.\notag
\end{align}

Since $\mathcal{Z}\subset\Omega$, the equation (\ref{Propo:energy:test:2}) holds for any $n\geq 6$ by the definition of $\mathcal{Z}$. If $n=4,5$, then $d=1$. In this case,
$$
\limsup_{x\to x_0}d_{g_0}(x,x_0)^{2-d}|W_{g_0}(x)|=0\,,
\:\:\:\:\:\:\text{for all}\:x_0\in\d M\,.
$$
Hence, $x_0\in\mathcal{Z}$ if and only if $\limsup_{x\to x_0}d_{g_0}(x,x_0)^{1-d}|\pi_{g_0}(x)|=0$.
Thus,
$$
\int_{D_{\rho}(x_0)}|\pi_g(x)|^2d_{g_0}(x,x_0)^{5-2n}\ds_{g_0}
=\infty\,,
\:\:\:\:\:\text{for all}\:x_0\in \d M\backslash \Omega\,,\notag
$$
and (\ref{Propo:energy:test:2}) holds.

Since $\d M\backslash\Omega$ is compact, it follows from Dini's theorem that 
\ba\label{Propo:energy:test:3}
\inf_{x_0\in\d M\backslash\Omega}&\int_{B_{\rho}(x_0)}|W_g(x)|^2(\e+d_{g_0}(x,x_0))^{6-2n}\dv_{g_0}
\\
&+\inf_{x_0\in\d M\backslash\Omega}\int_{D_{\rho}(x_0)}|\pi_g(x)|^2(\e+d_{g_0}(x,x_0))^{5-2n}\ds_{g_0}
\to \infty\,,
\:\:\:\:\:\text{as}\:\e\to 0\,.\notag
\end{align}

By the identities (\ref{Propo:energy:test:1}) and (\ref{Propo:energy:test:3}), we can choose $\e_0\in (0,\rho/2]$ such that
\ba
\inf_{x_0\in \Omega}\mathcal{I}(x_0,\rho)
&+\inf_{x_0\in\d M\backslash\Omega}\theta\int_{B_{\rho}(x_0)}|W_g(x)|^2(\e+d_{g_0}(x,x_0))^{6-2n}\dv_{g_0}\notag
\\
&+\inf_{x_0\in\d M\backslash\Omega}\theta\int_{D_{\rho}(x_0)}|\pi_g(x)|^2(\e+d_{g_0}(x,x_0))^{5-2n}\ds_{g_0}\notag
\\
&\hspace{3cm}>\:C_0\rho^{2d+4-n}+C_0\rho^{2-n}\frac{1}{\log(\rho/\e)}\,,\notag
\end{align}
for all $\e\in(0,\e_0]$. Now the assertion follows from Corollary \ref{corol8}.
\ep


\subsection{Further estimates}
In this section, we prove some results to be used in the next section. We use the same notations  of Section \ref{sub:sec:deftestfunct}.
Since $M$ is compact, we can assume that $\frac{1}{2}d_{g_0}(x_0,x)\leq d_{g_{x_0}}(x_0,x)\leq 2d_{g_0}(x_0,x)$ and $\frac{1}{2}d_{g_0}(x_0,x)\leq |\psi_{x_0}^{-1}(x)|\leq 2d_{g_0}(x_0,x)$ for all $x\in \psi_{x_0}(B^+_{2\rho}(0))$ and $x_0\in \d M$.
\begin{proposition}\label{PropoB.1}
If $2\e\leq \rho$, then we have
\ba
&\left|\frac{4(n-1)}{n-2}\Delta_{g_{x_0}}\ubar(x)-R_{g_{x_0}}\ubar(x)\right|\notag
\\
&\hspace{2cm}\leq C\left(\frac{\e}{\e^2+d_{g_{x_0}}(x,x_0)^2}\right)^{\frac{n-2}{2}}1_{\psi_{x_0}(B^+_{\rho}(0))}(x)
\\
&\hspace{2.5cm}+C\left\{
\e^{\frac{n}{2}}\rho^{-1-n}
+\e^{\frac{n-2}{2}}\rho^{1-n}
\right\}1_{\psi_{x_0}(B^+_{2\rho}(0)\backslash B^+_{\rho}(0))}(x)\notag
\end{align}
for all $x\in M$, and
\ba
&\left|\frac{2(n-1)}{n-2}\frac{\d}{\d\eta_{g_{x_0}}}\ubar(x)-\cmedia_{g_{x_0}}\ubar(x)+\cminfbar\ubar^{\frac{n}{n-2}}(x)\right|\notag
\\
&\hspace{2cm}\leq C\left(\frac{\e}{\e^2+d_{g_{x_0}}(x,x_0)^2}\right)^{\frac{n-2}{2}}1_{\psi_{x_0}(\d'B^+_{\rho}(0))}(x)\notag
\\
&\hspace{3cm}+C\left(\frac{\e}{\e^2+d_{g_{x_0}}(x,x_0)^2}\right)^{\frac{n}{2}}1_{\d M\backslash \psi_{x_0}(\d'B^+_{\rho}(0))}(x)\notag
\end{align}
for all $x\in \d M$.
\end{proposition}
\bp
In order to simplify our notations, we identify points $\psi_{x_0}(x)\in M$ with $x\in B^+_{2\rho}(0)$, omitting the symbol $\psi_{x_0}$. In particular, we identify $x_0\in \d M$ with $0\in B^+_{2\rho}(0)$. Recall that we are assuming $\cminfbar=2(n-1)$.

By the definition of $\ubar$,
\ba
\Delta_{g_{x_0}}&\ubar-\frac{n-2}{4(n-1)}R_{g_{x_0}}\ubar\notag
\\
&=\Delta_{g_{x_0}}\eta_{\rho}\cdot(\U+\phi_{x_0}-\e^{\frac{n-2}{2}}|x|^{2-n})\notag
-\Delta_{g_{x_0}}\eta_{\rho}\cdot\e^{\frac{n-2}{2}}(G_{x_0}-|x|^{2-n})\notag
\\
&\hspace{0.5cm}+2<d\eta_{\rho},d(\U+\phi_{x_0}-\e^{\frac{n-2}{2}}|x|^{2-n})>_{g_{x_0}}\notag
\\
&\hspace{0.5cm}-2\e^{\frac{n-2}{2}}<d\eta_{\rho},d(G_{x_0}-|x|^{2-n})>_{g_{x_0}}\notag
\\
&\hspace{0.5cm}+\eta_{\rho}\cdot\left(\Delta_{g_{x_0}}\U-\frac{n-2}{4(n-1)}R_{g_{x_0}}\U+\Delta\phi_{x_0}\right)\notag
\\
&\hspace{0.5cm}+\eta_{\rho}\cdot\left((\Delta_{g_{x_0}}-\Delta)\phi_{x_0}-\frac{n-2}{4(n-1)}R_{g_{x_0}}\phi_{x_0}\right)\notag
\end{align}
and
\ba
\critbordo&\frac{\d}{\d\eta_{g_{x_0}}}\ubar-\cmedia_{g_{x_0}}\ubar+2(n-1)\ubar^{\frac{n}{n-2}}\notag
\\
&=\critbordo\d_n\left(\eta_{\rho}(\U+\phi_{x_0})+\e^{\frac{n-2}{2}}(1-\eta_{\rho})G_{x_0}\right)\notag
\\
&\hspace{0.5cm}-\cmedia_{g_{x_0}}\left(\eta_{\rho}(\U+\phi_{x_0})+\e^{\frac{n-2}{2}}(1-\eta_{\rho})G_{x_0}\right)\notag
\\
&\hspace{0.5cm}+2(n-1)\left(\eta_{\rho}(\U+\phi_{x_0})+\e^{\frac{n-2}{2}}(1-\eta_{\rho})G_{x_0}\right)^{\frac{n}{n-2}}\notag
\\
&=\eta_{\rho}\cdot \left\{\critbordo\d_n(\U+\phi_{x_0})+2(n-1)(\U+\phi_{x_0})^{\frac{n}{n-2}}\right\}\notag
\\
&\hspace{0.5cm}-\cmedia_{g_{x_0}}\eta_{\rho}(\U+\phi_{x_0})\notag
\\
&\hspace{0.5cm}+2(n-1)\left\{\left(\eta_{\rho}(\U+\phi_{x_0})+\e^{\frac{n-2}{2}}(1-\eta_{\rho})G_{x_0}\right)^{\frac{n}{n-2}}
-\eta_{\rho}(\U+\phi_{x_0})^{\frac{n}{n-2}}\right\}\,.\notag
\end{align}
Now the result easily follows.
\ep
\begin{lemma}\label{LemmaB.4}
We have
$$
\int_{\d M}\uum\udois^{\frac{n}{n-2}}\ds_{g_0}
\leq C\left(\frac{\e_1\e_2}{\e_2^2+d_{g_0}(x_1,x_2)^2}\right)^{\frac{n-2}{2}}\,.
$$
\end{lemma}
\bp
As in \cite[Lemma B.4]{brendle-flow}, one can prove that

\ba
\int_{\d M}&\left(\frac{\e_1}{\e_1^2+d_{g_0}(x,x_1)^2}\right)^{\frac{n-2}{2}}
\left(\frac{\e_2}{\e_2^2+d_{g_0}(x,x_2)^2}\right)^{\frac{n}{2}}\ds_{g_0}\notag
\\
&\leq C\left(\frac{\e_1\e_2}{\e_2^2+d_{g_0}(x_1,x_2)^2}\right)^{\frac{n-2}{2}}\,.\notag
\end{align}
From this the assertion follows.
\ep
\begin{lemma}\label{LemmaB.5}
For all $\e_1, \e_2\leq \rho^2\leq 1/4$,
$$
\int_M\uum|\Delta_{g_0}\udois|\dv_{g_0}
\leq C\rho\left(\frac{\e_1\e_2}{\e_2^2+d_{g_0}(x_1,x_2)^2}\right)^{\frac{n-2}{2}}
$$
and
\ba
\int_{\d M}&\left|\frac{2(n-1)}{n-2}\frac{\d}{\d\eta_{g_0}}\udois-\cmz\udois+\cminfbar\udois^{\frac{n}{n-2}}\right|\uum\ds_{g_0}\notag
\\
&\hspace{3cm}\leq
C\left(\rho+\frac{\e_2}{\rho}\right)\left(\frac{\e_1\e_2}{\e_2^2+d_{g_0}(x_1,x_2)^2}\right)^{\frac{n-2}{2}}\,.\notag
\end{align}
\end{lemma}
\bp
It follows from Proposition \ref{PropoB.1} that
\ba
&\left|\frac{4(n-1)}{n-2}\Delta_{g_{x_2}}\bar{U}_{(x_2,\e_2)}(x)-R_{g_{x_2}}\bar{U}_{(x_2,\e_2)}(x)\right|\notag
\\
&\hspace{2cm}\leq C\rho^{-1}\left(\frac{\e_2}{\e_2^2+d_{g_0}(x,x_2)^2}\right)^{\frac{n-2}{2}}1_{\{d_{g_0}(y,x_2)\leq 4\rho\}}(x)\notag
\end{align}
for all $x\in M$, and
\ba
&\left|\frac{2(n-1)}{n-2}\frac{\d}{\d\eta_{g_{x_2}}}\bar{U}_{(x_2,\e_2)}(x)-\cmedia_{g_{x_2}}\bar{U}_{(x_2,\e_2)}(x)+\cminfbar\bar{U}_{(x_2,\e_2)}^{\frac{n}{n-2}}(x)\right|\notag
\\
&\hspace{2cm}\leq C\left(\frac{\e_2}{\e_2^2+d_{g_0}(x,x_2)^2}\right)^{\frac{n-2}{2}}1_{\{d_{g_0}(y,x_2)\leq 4\rho\}\cap \d M}(x)\notag
\\
&\hspace{3cm}+C\left(\frac{\e_2}{\e_2^2+d_{g_0}(x,x_2)^2}\right)^{\frac{n}{2}}1_{\{d_{g_0}(y,x_2)\geq \rho/2\}\cap \d M}(x)\notag
\end{align}
for all $x\in \d M$. 

Proceeding as in \cite[Lemma B.5]{brendle-flow}, we can show that
\ba
\int_{\{d_{g_0}(y,x_2)\leq 4\rho\}}&\left(\frac{\e_1}{\e_1^2+d_{g_0}(x,x_1)^2}\right)^{\frac{n-2}{2}}
\left(\frac{\e_2}{\e_2^2+d_{g_0}(x,x_2)^2}\right)^{\frac{n-2}{2}}\dv_{g_0}(y)\notag
\\
&\leq C\rho^2\left(\frac{\e_1\e_2}{\e_2^2+d_{g_0}(x_1,x_2)^2}\right)^{\frac{n-2}{2}}\,,\notag
\end{align}
\ba
\int_{\{d_{g_0}(y,x_2)\leq 4\rho\}\cap \d M}&\left(\frac{\e_1}{\e_1^2+d_{g_0}(x,x_1)^2}\right)^{\frac{n-2}{2}}
\left(\frac{\e_2}{\e_2^2+d_{g_0}(x,x_2)^2}\right)^{\frac{n-2}{2}}\ds_{g_0}(y)\notag
\\
&\leq C\rho\left(\frac{\e_1\e_2}{\e_2^2+d_{g_0}(x_1,x_2)^2}\right)^{\frac{n-2}{2}}\,,\notag
\end{align}
and
\ba
\int_{\{d_{g_0}(y,x_2)\geq \rho/2\}\cap \d M}&\left(\frac{\e_1}{\e_1^2+d_{g_0}(x,x_1)^2}\right)^{\frac{n-2}{2}}
\left(\frac{\e_2}{\e_2^2+d_{g_0}(x,x_2)^2}\right)^{\frac{n}{2}}\ds_{g_0}(y)\notag
\\
&\leq C\frac{\e_2}{\rho}\left(\frac{\e_1\e_2}{\e_2^2+d_{g_0}(x_1,x_2)^2}\right)^{\frac{n-2}{2}}\,.\notag
\end{align}

Now the assertion follows.
\ep


\section{Blow-up analysis}\label{sec:blowup}
In this section, we carry out the blow-up analysis for sequences of solutions to the equations (\ref{eq:evol:u}) that will be necessary for the proof of Theorem \ref{main:thm}.

Let $u(t)$, $t\geq 0$, be the solution of (\ref{eq:evol:u}) obtained in Section \ref{sec:prelim}, and let
$\{t_{\nu}\}_{\nu=1}^{\infty}$ be a sequence satisfying $\lim_{\nu\to\infty}t_{\nu}=\infty$. We set $\u=u(t_{\nu})$ and $\g=g(t_{\nu})=\u^{\frac{4}{n-2}}g_0$. Then 
$$
\int_{\d M}\u^{\critbordo}\ds_{g_0}=\int_{\d M}\ds_{\g}=1\,,\:\:\:\:\text{for all}\:\nu\,.
$$
It follows from Corollary \ref{Corol3.2} that
$$ 
\int_{\d M}
\left|\frac{2(n-1)}{n-2}\frac{\d\u}{\eta_{g_0}}-\cmedia_{g_0}\u+\cminfbar\u^{\frac{n}{n-2}}\right|^{\conjbordo}
\ds_{g_{0}}
=\int_{\d M}|\cmv-\cminfbar|^{\conjbordo}\ds_{\g}
\to 0
$$
as $\nu\to\infty$.
\begin{proposition}\label{Propo4.1}
After passing to a subsequence, there exist an integer $m\geq 0$, a smooth function $\uinf\geq 0$, and a sequence of $m$-tuplets  $\{(x_{k,\nu}^*,\e_{k,\nu}^*)_{1\leq k\leq m}\}_{\nu=1}^{\infty}$,  such that:

(i) The function $\uinf$ satisfies
\begin{equation}\notag
\begin{cases}
\D_{g_0}\uinf=0\,,&\text{in}\:M\,,
\\
\frac{2(n-1)}{n-2}\frac{\d}{\d\eta_{g_0}}\uinf-\cmedia_{g_0}\uinf+\cminfbar\uinf^{\frac{n}{n-2}}=0\,,&\text{on}\:\d M\,.
\end{cases} 
\end{equation}

(ii) For all $i\neq j$,
$$
\lim_{\nu\to\infty}\left\{\frac{\e_{i,\nu}^*}{\e_{j,\nu}^*}+\frac{\e_{j,\nu}^*}{\e_{i,\nu}^*}
+\frac{d_{g_0}(x_{i,\nu}^*,x_{j,\nu}^*)^2}{\e_{i,\nu}^*\e_{j,\nu}^*}\right\}=\infty\,.
$$

(iii) We have
$$
\lim_{\nu\to\infty}
\big{\|}\u-\uinf-\sum_{k=1}^{m}\bar{u}_{(x_{k,\nu}^*,\e_{k,\nu}^*)}\big{\|}_{H^1(M)}= 0\,,
$$
where the functions $\bar{u}_{(x_{k,\nu}^*,\e_{k,\nu}^*)}$ were defined by equation (\ref{def:test:func:u}).
\end{proposition}
\bp
This is the content of \cite{almaraz4}. Observe that, although functions $\bar{u}_{(x_{k,\nu}^*,\e_{k,\nu}^*)}$ differ from the ones used in \cite{almaraz4}, it's easy to check that their difference converge to zero in $H^1(M)$.
The regularity of $\uinf$ was established by P. Cherrier in \cite{cherrier}.
\ep
\begin{proposition}\label{Propo4.2}
If $\uinf(x)=0$ for some $x\in M$, then $\uinf\equiv 0$.
\end{proposition}
\bp
This is just a consequence of the maximum principle.
\ep

Define the functionals
$$
E(u)=\frac{\frac{4(n-1)}{n-2}\int_M|du|_{g_0}^2\dv_{g_0}+2\int_{\d M}\cmedia_{g_0}u^2\ds_{g_0}}{\left(\int_{\d M}u^{\frac{2(n-1)}{n-2}}\ds_{g_0}\right)^{\frac{n-2}{n-1}}}
$$
and
$$
F(u)=\frac{\frac{4(n-1)}{n-2}\int_M|du|_{g_0}^2\dv_{g_0}+2\int_{\d M}\cmedia_{g_0}u^2\ds_{g_0}}{\int_{\d M}u^{\frac{2(n-1)}{n-2}}\ds_{g_0}}\,.
$$
Observe that $\cminfbar=\frac{1}{2}F(\uinf)$. Hence,
\ba
1&=\lim_{\nu\to\infty}\int_{\d M}u_{\nu}^{\critbordo}\ds_{g_0}=\lim_{\nu\to\infty}\left\{\int_{\d M}\uinf^{\critbordo}\ds_{g_0}
+\sum_{k=1}^{m}\int_{\d M}\bar{u}_{(x^*_{k,\nu},\e^*_{k,\nu})}^{\critbordo}\ds_{g_0}\right\}\notag
\\
&=\left(\frac{E(\uinf)}{2\cminfbar}\right)^{n-1}+m\left(\frac{Q(B^n,\d B)}{2\cminfbar}\right)^{n-1}\,.\notag
\end{align}
Thus, 
\begin{equation}\label{eq:cminbar}
\cminfbar=\frac{1}{2}\left(E(\uinf)^{n-1}+mQ(B^n,\d B)^{n-1}\right)^{\frac{1}{n-1}}\,.
\end{equation}

\subsection{The case $\uinf\equiv 0$}

We set 
\ba
\mathcal{A}_{\nu}=\Big{\{}(x_k,\e_k,\a_k)_{k=1,...,m}\in &(\d M\times\R_+\times\R_+)^m\,,\:\text{such that}\notag
\\
&d_{g_0}(x_{k}, x^*_{k,\nu})\leq \e^*_{k,\nu}\,,\:\frac{1}{2}\leq\frac{\e_k}{\e^*_{k,\nu}}\leq 2
\,,\:\frac{1}{2}\leq\a_k\leq 2\Big{\}}\,.\notag
\end{align}
For each $\nu$, we can choose a triplet $(x_{k,\nu},\e_{k,\nu},\a_{k,\nu})_{k=1,...,m}\in\mathcal{A}_{\nu}$ such that
\ba
&\int_M\frac{2(n-1)}{n-2}\big| d(u_{\nu}-\sum_{k=1}^{m}\a_{k,\nu}\uknu)\big|_{g_0}^2\dv_{g_0}
+\int_{\d M}\cmedia_{g_0}\big(u_{\nu}-\sum_{k=1}^{m}\a_{k,\nu}\uknu\big)^2\ds_{g_0}\notag
\\
&\hspace{0.5cm}\leq\int_M\frac{2(n-1)}{n-2}\big{|}d(u_{\nu}-\sum_{k=1}^{m}\a_{k}\uk)\big{|}_{g_0}^2\dv_{g_0}
+\int_{\d M}\cmedia_{g_0}\big(u_{\nu}-\sum_{k=1}^{m}\a_{k}\uk\big)^2\ds_{g_0}\notag
\end{align}
for all $(x_k,\e_k,\a_k)_{k=1,...,m}\in\mathcal{A}_{\nu}$.

\vspace{0.2cm}
The proof of the next two propositions are the same of Propositions 5.1 and 5.2 in \cite{brendle-flow}:
\begin{proposition}\label{Propo5.1}
We have: 

(i) For all $i\neq j$, 
$$
\lim_{\nu\to\infty}\left\{\frac{\e_{i,\nu}}{\e_{j,\nu}}+\frac{\e_{j,\nu}}{\e_{i,\nu}}
+\frac{d_{g_0}(x_{i,\nu},x_{j,\nu})^2}{\e_{i,\nu}\e_{j,\nu}}\right\}=\infty\,.
$$

(ii) We have
$$
\lim_{\nu\to\infty}
\big{\|}\u-\sum_{k=1}^{m}\a_{k,\nu}\uknu\big{\|}_{H^1(M)}= 0\,.
$$
\end{proposition}
\begin{proposition}\label{Propo5.2}
We have 
$$
d_{g_0}(x_{k,\nu}, x^*_{k,\nu})\leq o(1)\e^*_{k,\nu}\,,\:\:\:
\frac{\e_{k,\nu}}{\e^*_{k,\nu}}=1+o(1)\,,\:\:\:
\text{and}\:\:\:\a_{k,\nu}=1+o(1)\,,
$$
for all $k=1,...,m$. In particular, $(x_{k, \nu},\e_{k, \nu},\a_{k, \nu})_{k=1,...,m}$ is an interior point of $\mathcal{A}_{\nu}$ for $\nu$ sufficiently large.
\end{proposition}
\begin{convention}
Assume that $\e_{i,\nu}\leq \e_{j,\nu}$ for all $i\leq j$, without loss of generality.
\end{convention}
\begin{notation}
We write $u_{\nu}=v_{\nu}+w_{\nu}$, where
\begin{equation}\label{def:v:w:1}
v_{\nu}=\sum_{k=1}^{m}\a_{k,\nu}\uknu\:\:\:\:\text{and}\:\:\:\:
w_{\nu}=u_{\nu}-\sum_{k=1}^{m}\a_{k,\nu}\uknu\,.
\end{equation}
\end{notation}
Observe that by Proposition \ref{Propo5.1} we have
\begin{equation}\label{propr:w:1}
\int_M\frac{2(n-1)}{n-2}|dw_{\nu}|_{g_0}^2\dv_{g_0}
+\int_{\d M}\cmedia_{g_0}w_{\nu}^2\ds_{g_0}=o(1)\,.
\end{equation}
\begin{proposition}\label{Propo5.3}
Let $\psi_{k,\nu}:B^+_{2\rho}(0)\to M$ be Fermi coordinates centered at $x_{k,\nu}$. If we set 
$$
C_{\nu}=\left(\int_{\d M}|\w|^{\critbordo}\ds_{g_0}\right)^{\frac{n-2}{2(n-1)}}
+\left(\int_{M}|\w|^{\crit}\dv_{g_0}\right)^{\frac{n-2}{2n}}\,,
$$ 
then for all $k=1,...,m$, we have:

(i) $\displaystyle{\big{|} \int_{\d M} \uknu^{\frac{n}{n-2}}\,\w\,\ds_{g_0}\big{|}
\leq o(1)C_{\nu}}$\,.

(ii) $\displaystyle{\big{|} \int_{\psi_{k,\nu}(\d 'B^+_{2\rho}(0))} \uknu^{\frac{n}{n-2}}
\frac{\e_{k,\nu}^2-|\psi_{k,\nu}^{-1}(x)|^2}{\e_{k,\nu}^2+|\psi_{k,\nu}^{-1}(x)|^2}\,\w\,\ds_{g_0}\big{|}
\leq o(1)C_{\nu}}$\,.

(iii) $\displaystyle{\big{|} \int_{\psi_{k,\nu}(\d 'B^+_{2\rho}(0))} \uknu^{\frac{n}{n-2}}
\frac{\e_{k,\nu}\psi_{k,\nu}^{-1}(x)}{\e_{k,\nu}^2+|\psi_{k,\nu}^{-1}(x)|^2}\,\w\,\ds_{g_0}\big{|}
\leq o(1)C_{\nu}}$\,.
\end{proposition}
\bp
\noindent 
(i) It follows from the definition of $(x_{k,\nu},\e_{k,\nu},\a_{k,\nu})_{k=1,...,m}$ that
$$
\int_M\frac{2(n-1)}{n-2}<d\uknu,dw_{\nu}>_{g_0}\dv_{g_0}
+\int_{\d M}\cmedia_{g_0}\uknu w_{\nu}\,\ds_{g_0}=0\,.
$$
Integrating by parts, 
\ba
\int_M&\frac{2(n-1)}{n-2}\Delta_{g_0}\uknu w_{\nu}\,\dv_{g_0}\notag
\\
&+\int_{\d M}\left\{\frac{2(n-1)}{n-2}\frac{\d}{\d\eta_{g_0}}\uknu
-\cmedia_{g_0}\uknu \right\}w_{\nu}\,\ds_{g_0}=0\,,\notag
\end{align}
which implies
\ba\label{Propo5.3:1}
\int_{\d M}&\cminfbar \uknu^{\frac{n}{n-2}}w_{\nu}\,\ds_{g_0}
\\
&=\int_{M}\frac{2(n-1)}{n-2}\Delta_{g_0}\uknu w_{\nu}\,\dv_{g_0}\notag
\\
&+\int_{\d M}\left\{\frac{2(n-1)}{n-2}\frac{\d}{\d\eta_{g_0}}\uknu
-\cmedia_{g_0}\uknu +\cminfbar \uknu^{\frac{n}{n-2}}\right\}w_{\nu}\,\ds_{g_0}\,.\notag
\end{align}
Then, using Proposition \ref{PropoB.1} and a conformal change of the metric, we can prove that
$$
\left|\int_{\d M}\uknu^{\frac{n}{n-2}}w_{\nu}\ds_{g_0}\right|
\leq o(1)\left(\|w_{\nu}\|_{L^{\frac{2n}{n-2}}(M)}+\|w_{\nu}\|_{L^{\frac{2(n-1)}{n-2}}(\d M)}\right)
$$
for $k=1,...,m$.

\bigskip\noindent 
(ii) Let us set $\tilde{\psi}_{k,\nu}=\frac{\d}{\d\e}\Big{|}_{\e=\e_{k,\nu}}\bar{u}_{(x_{k,\nu},\e)}$. Similarly to (\ref{Propo5.3:1}) we obtain
\ba
\int_{\d M}&\frac{n}{n-2}\cminfbar \uknu^{\frac{2}{n-2}}\tilde{\psi}_{k,\nu} w_{\nu}\,\ds_{g_0}\notag
\\
&=\int_{M}\frac{2(n-1)}{n-2}\Delta_{g_0}\tilde{\psi}_{k,\nu} w_{\nu}\,\dv_{g_0}\notag
\\
&\hspace{0.3cm}+\int_{\d M}\left\{\frac{2(n-1)}{n-2}\frac{\d}{\d\eta_{g_0}}\tilde{\psi}_{k,\nu}
-\cmedia_{g_0}\tilde{\psi}_{k,\nu} +\frac{n}{n-2}\cminfbar \uknu^{\frac{2}{n-2}}\tilde{\psi}_{k,\nu}\right\}w_{\nu}\,\ds_{g_0}\,.\notag
\end{align}
Using the estimate (\ref{est:V}) we observe that $\tilde{\psi}_{k,\nu}$ satisfies 
$$
\e_{k,\nu}\uknu^{-1}\tilde{\psi}_{k,\nu}=\frac{n-2}{2}\frac{|x|^2-\e_{k,\nu}^2}{(\e_{k,\nu}+x_n)^2+|x|^2}
+O\big{(}(\e_{k,\nu}+|\bar{x|})\big{)}\,,
\:\:\:\:\text{in}\:\:B^+_{\rho}(0)\,.
$$
Now the result follows as in the item (i), and the item (iii) follows similarly.
\ep

\begin{proposition}\label{Propo5.4}
There exists $c>0$ such that 
\ba
\frac{n}{n-2}\cminfbar\int_{\d M}\sum_{k=1}^{m}&\uknu^{\frac{2}{n-2}}\w^2\,\ds_{g_0}\notag
\\
&\leq
(1-c)\left\{\int_{M}\frac{2(n-1)}{n-2}|d\w|_{g_0}^2\dv_{g_0}
+\int_{\d M}\cmedia_{g_0}\w^2\,\ds_{g_0}\right\}\notag
\end{align}
for all $\nu$ sufficiently large.
\end{proposition}
\bp
Suppose by contradiction this is not true. Upon rescaling we can find a sequence $\{\tilde{w}_{\nu}\}$ satisfying
\begin{equation}\label{Propo5.4:1}
\int_M \frac{2(n-1)}{n-2}|d\tilde{w}_{\nu}|_{g_0}^2\dv_{g_0}
+\int_{\d M}\cmedia_{g_0}\tilde{w}_{\nu}^2\,\ds_{g_0}=1
\end{equation}
and
\begin{equation}\label{Propo5.4:2}
\lim_{\nu\to\infty}\frac{n}{n-2}\cminfbar\int_{\d M}\sum_{k=1}^{m}\uknu^{\frac{2}{n-2}}\tilde{w}_{\nu}^2\,\ds_{g_0}
\geq 1\,.
\end{equation}
Observe that the identity (\ref{Propo5.4:1}) implies
\begin{equation}\label{Propo5.4:w:bordo}
\int_{\d M}|\tilde{w}_{\nu}|^{\critbordo}\ds_{g_0}\leq \QM^{-\frac{n-1}{n-2}}
\end{equation}
and
\begin{equation}\label{Propo5.4:w:int}
\int_{M}|\tilde{w}_{\nu}|^{\crit}\dv_{g_0}\leq Q(M)^{-\frac{n}{n-2}}\,,
\end{equation}
where $Q(M)$ is the conformal invariant defined in \cite{escobar2}, which has the same sign of $\QM$ (see \cite[Proposition 1.2]{escobar3}).

In view of Proposition \ref{Propo5.1}, we can choose a sequence $\{N_{\nu}\}$, such that $N_{\nu}\to\infty$, $N_{\nu}\e_{k,\nu}\to 0$ for all $k=1,...,m$, and
$$
\frac{\e_{j,\nu}+d_{g_0}(x_{i,\nu},x_{j,\nu})}{N_{\nu}\e_{i,\nu}}\to\infty
\:\:\:\:\:\text{for all}\:\:i< j\,.
$$  
Set $\Omega_{j,\nu}=B_{N_{\nu}\e_{j,\nu}}(x_{j,\nu})\backslash \bigcup_{i=1}^{j-1}B_{N_{\nu}\e_{i,\nu}}(x_{i,\nu})$ for $1\leq j\leq m$. It follows from (\ref{Propo5.4:1}) and (\ref{Propo5.4:2}) that there exists $1\leq k\leq m$ such that
$$
\lim_{\nu\to\infty}\int_{\d M}\uknu^{\frac{2}{n-2}}\tilde{w}_{\nu}^2\,\ds_{g_0}>0
$$
and
\ba
\lim_{\nu\to\infty}&\left\{\int_{\Omega_{k,\nu}} \frac{2(n-1)}{n-2}|d\tilde{w}_{\nu}|_{g_0}^2\dv_{g_0}
+\int_{\Omega_{k,\nu}\cap\d\Rn}\cmedia_{g_0}\tilde{w}_{\nu}^2\,\ds_{g_0}\right\}\notag
\\
&\leq 
\lim_{\nu\to\infty}\frac{n}{n-2}\cminfbar\int_{\d M}\uknu^{\frac{2}{n-2}}\tilde{w}_{\nu}^2\,\ds_{g_0}\,.\notag
\end{align}

For each $\nu$, let $\psi_{\nu}:B^+_{\rho}(0)\subset \Rn\to M$ be Fermi coordinates centered at $x_{k,\nu}$. We set
$$
\hat{w}_{\nu}(y)=\e_{k,\nu}^{\frac{n-2}{2}}\tilde{w}_{\nu}(\psi_{\nu}(\e_{k,\nu}y))\,.
$$
Then
\begin{equation}\label{Propo5.4:3'}
\lim_{\nu\to\infty}\int_{\{y\in\Rn,\:|y|\leq N_{\nu}\}}\frac{2(n-1)}{n-2}|d\hat{w}_{\nu}(y)|^2\,dy
\leq 1\,,
\end{equation}
\begin{equation}\label{Propo5.4:2'}
\lim_{\nu\to\infty}\int_{\{y\in\d\Rn,\:|y|\leq N_{\nu}\}}|\hat{w}_{\nu}(y)|^{\critbordo}\,dy
\leq \QM^{-\frac{n-1}{n-2}}\,,
\end{equation}
and
\begin{equation}\notag
\lim_{\nu\to\infty}\int_{\{y\in\Rn,\:|y|\leq N_{\nu}\}}|\hat{w}_{\nu}(y)|^{\crit}\,dy
\leq Q(M)^{-\frac{n}{n-2}}\,.
\end{equation}
Thus, we can assume that $\hat{w}_{\nu}\rightharpoonup \hat{w}$ in $H_{loc}^1(\Rn)$ for some $\hat{w}$ satisfying
\begin{equation}\label{Propo5.4:3}
\int_{\d\Rn}\frac{1}{1+|y|^2}\hat{w}^2(y)\,dy>0
\end{equation}
and
\begin{equation}\label{Propo5.4:4}
\int_{\Rn}|d\hat{w}(y)|^2dy
\leq n\int_{\d\Rn}\frac{1}{1+|y|^2}\hat{w}^2(y)\,dy\,.
\end{equation}
Moreover, Proposition \ref{Propo5.3}, together with the inequalities (\ref{Propo5.4:w:bordo}) and (\ref{Propo5.4:w:int}), implies that
\begin{equation}\label{Propo5.4:5}
\int_{\d\Rn}\left(\frac{1}{1+|y|^2}\right)^{\frac{n}{2}}\hat{w}(y)\,dy=0\,,
\end{equation}
\begin{equation}\label{Propo5.4:6}
\int_{\d\Rn}\left(\frac{1}{1+|y|^2}\right)^{\frac{n}{2}}\frac{1-|y|^2}{1+|y|^2}\hat{w}(y)\,dy=0\,,
\end{equation}
\begin{equation}\label{Propo5.4:7}
\int_{\d\Rn}\left(\frac{1}{1+|y|^2}\right)^{\frac{n}{2}}\frac{y_j}{1+|y|^2}\hat{w}(y)\,dy=0\,,
\:\:\:\:j=1,..,n-1\,,
\end{equation}
where $y=(y_1,...,y_{n-1},0)$.

Let $B_{1/2}$ be the Euclidean ball in $\R^{n}$ of radius $1/2$ with center $(0,...,0,-1/2)$. 
We set $\mathcal{C}=\{w\in H^1(B_{1/2});\:\int_{\d B_{1/2}}w\,\ds=0\}$. Observe that 
$$
\inf_{0\neq w\in\mathcal{C}}\frac{\int_{B_{1/2}}|dw|^2dy}{\int_{\d B_{1/2}}w^2\ds}=2\,,
$$
and this infimum is realized only by the coordinate functions $z_1,...,z_n$ of $\R^n$, taken with center $(0,...,0,-1/2)$, restricted to $B_{1/2}$. 

The ball $B_{1/2}$ is conformally equivalent to the half-space $\Rn$ by means of the inversion 
$F:\mathbb{R}_+^n\to B_{1/2}\backslash\{ (0,...,0,-1)\}$ given by
$$F(y_1,...y_n)=\frac{(y_1,...,y_{n-1},y_n+1)}{y_1^2+...+y_{n-1}^2+(y_n+1)^2}+(0,...,0,-1)\,.$$
An easy calculation shows that $F$ is a conformal map and $F^*g_{eucl} = U_1^{\frac{4}{n-2}}g_{eucl}$ in $\Rn$, where $g_{eucl}$ is the Euclidean metric, and 
$$
z_j\circ F(y)=\frac{y_j}{(y_n+1)^2+|\bar{y}|^2}\,,
\:\:\:\:\:
z_n\circ F(y)=\frac{1}{2}\frac{1-|y|^2}{(y_n+1)^2+|\bar{y}|^2}\,.
$$
(See (\ref{eq:def:U}) for the definition of $U_1$.)

Using the estimate (\ref{Propo5.4:3'}), we can easily check that 
$$
w=(U_1^{-1}\hat{w})\circ F^{-1}\in H^1(B_{1/2})\,.
$$
In view of the identities \eqref{Propo5.4:5}, \eqref{Propo5.4:6} and \eqref{Propo5.4:7}, $w$ is $L^2(\d B_{1/2})$-orthogonal to the functions $1,z_1,...,z_n$. Since by (\ref{Propo5.4:3}) we have $w\nequiv 0$ in $B_{1/2}$, this function satisfies 
$$
\int_{B_{1/2}}|dw|^2 dy > 2\int_{\d B_{1/2}}w^2 \ds\,,
$$
which corresponds to 
\begin{equation}\notag
\int_{\Rn}|d\hat{w}(y)|^2dy
- n\int_{\d\Rn}\frac{1}{1+|y|^2}\hat{w}^2(y)dy>0\,.
\end{equation} 
This contradicts the inequality (\ref{Propo5.4:4}).
\ep
\begin{corollary}\label{Corol5.5}
There exists $c>0$ such that 
$$
\frac{n}{n-2}\cminfbar\int_{\d M}v_{\nu}^{\frac{2}{n-2}}w_{\nu}^2\,\ds_{g_0}
\leq
(1-c)\left\{\int_{M}\frac{2(n-1)}{n-2}|d\w|_{g_0}^2\dv_{g_0}
+\int_{\d M}\cmedia_{g_0}\w^2\,\ds_{g_0}\right\}
$$
for all $\nu$ sufficiently large.
\end{corollary}
\bp
By the definition of $v_{\nu}$ (equation (\ref{def:v:w:1})), we have
$$
\lim_{\nu\to\infty}\int_{\d M}\big|v_{\nu}^{\frac{2}{n-2}}-\sum_{k=1}^{m}\bar{u}_{(x_{k,\nu}, \e_{k,\nu})}^{\frac{2}{n-2}}\big|^{n-1}\ds_{g_0}=0\,.
$$
Hence, the assertion follows from Proposition \ref{Propo5.4}.
\ep

The next proposition is similar to Proposition 5.6 of \cite{brendle-flow} and we will just outline its proof.
\begin{proposition}\label{Propo5.6}
For all $\nu$ sufficiently large, we have $E(v_{\nu})\leq\left\{\sum_{k=1}^{m}E(\uk)^{n-1}\right\}^{\frac{1}{n-1}}$.
\end{proposition}
\bp
Observe that, given $i<j$, there exist $C,c>0$  such that
$$
\uinu(x)^{\frac{n}{n-2}}\ujnu(x)
\geq c\left(\frac{\e_{i,\nu}\e_{j,\nu}}{\e_{j,\nu}^2+d_{g_0}(x_{i,\nu},x_{j,\nu})^2}\right)^{\frac{n-2}{2}}\e_{i,\nu}^{1-n}
$$
and
$$
\uinu(x)\ujnu(x)^{\frac{n}{n-2}}
\leq C\left(\frac{\e_{i,\nu}\e_{j,\nu}}{\e_{j,\nu}^2+d_{g_0}(x_{i,\nu},x_{j,\nu})^2}\right)^{\frac{n}{2}}\e_{i,\nu}^{1-n}\,,
$$
for all $x\in\d M$ such that $d_{g_0}(x,x_{i,\nu})\leq \e_{i,\nu}$, and $\nu$ sufficiently large. 

Proceeding as in \cite{brendle-flow}, we obtain
\ba\label{Propo5.6:4}
\frac{1}{2}E(v_{\nu})&\left(\int_{\d M}v_{\nu}^{\critbordo}\ds_{g_0}\right)^{\frac{n-2}{n-1}}
\\
&\leq
\frac{1}{2}\left(\sum_{k=1}^{m}E(\uknu)^{n-1}\right)^{\frac{1}{n-1}}\left(\int_{\d M}v_{\nu}^{\critbordo}\ds_{g_0}\right)^{\frac{n-2}{n-1}}\notag
\\
&\hspace{0.5cm}-\sum_{i<j}2\a_{i,\nu}\a_{j,\nu}\int_M\critbordo\Delta_{g_0}\ujnu\uinu\dv_{g_0}\notag
\\
&\hspace{0.5cm}-\sum_{i<j}2\a_{i,\nu}\a_{j,\nu}\int_{\d M}\Big(\critbordo\frac{\d}{\d\eta_{g_0}}\ujnu-\cmz\ujnu\notag
\\
&\hspace{5cm}+\frac{1}{2}F(\ujnu)\ujnu^{\frac{n}{n-2}}\Big)\uinu\ds_{g_0}\notag
\\
&\hspace{0.5cm}-c\sum_{i<j}\left(\frac{\e_{i,\nu}\e_{j,\nu}}{\e_{j,\nu}^2+d_{g_0}(x_{i,\nu},x_{j,\nu})^2}\right)^{\frac{n-2}{2}}\,.\notag
\end{align}
It follows from Lemmas \ref{LemmaB.4} and \ref{LemmaB.5} that
\ba\label{Propo5.6:5}
&\int_M\critbordo\left|\Delta_{g_0}\ujnu\uinu\right|\dv_{g_0}
\\
&\hspace{0.5cm}+\int_{\d M}\left|\critbordo\frac{\d}{\d\eta_{g_0}}\ujnu-\cmz\ujnu+\frac{1}{2}F(\ujnu)\ujnu^{\frac{n}{n-2}}\right|\uinu\ds_{g_0}\notag
\\
&\hspace{2cm}\leq
C\left(\rho+\frac{\e_{j,\nu}}{\rho}\right)
\left(\frac{\e_{i,\nu}\e_{j,\nu}}{\e_{j,\nu}^2+d_{g_0}(x_{i,\nu},x_{j,\nu})^2}\right)^{\frac{n-2}{2}}
+o(1)\left(\frac{\e_{i,\nu}\e_{j,\nu}}{\e_{j,\nu}^2+d_{g_0}(x_{i,\nu},x_{j,\nu})^2}\right)^{\frac{n-2}{2}}\,,\notag
\end{align}
where we used that $\lim_{\nu\to\infty}\big|\frac{1}{2}F(\ujnu)-\cminfbar\big|=0$, for all $j=1,...,m$. Now the assertion follows from the estimates (\ref{Propo5.6:4}) and (\ref{Propo5.6:5}), choosing $\rho$ small and $\nu$ large.
\ep
\begin{corollary}\label{Corol5.7}
Under the hypothesis of Theorem \ref{main:thm},  we have
$$
E(v_{\nu})\leq (mQ(B^n,\d B)^{n-1})^{\frac{1}{n-1}}
$$ 
for all $\nu$ sufficiently large.
\end{corollary}
\bp
Using Proposition \ref{Propo:energy:test}, we obtain $E(\uknu)\leq \Q$ for all $k=1,...,m$. Then the result follows from Proposition \ref{Propo5.6}
\ep

\subsection{The case $\uinf >0$}
\begin{proposition}\label{Propo6.1}
There exist sequences $\{\psi_a\}_{a\in\N}\subset C^{\infty}(M)$ and $\{\lambda_a\}_{a\in\N}\subset \R$, with $\lambda_a>0$,  satisfying:

(i) For all $a\in\N$, 
\begin{equation}\notag
\begin{cases} 
\Delta_{g_0}\psi_a=0\,,&\text{in}\:M\,,
\\
\critbordo\frac{\d}{\d\eta_{g_0}}\psi_a
-\cmz\psi_a+\lambda_au_{\infty}^{\frac{2}{n-2}}\psi_a=0\,,&\text{on}\:\d M\,.
\end{cases} 
\end{equation}

(ii) For all $a,b\in\N$,
$$
\int_{\d M}\psi_a\psi_bu_{\infty}^{\frac{2}{n-2}}\ds_{g_0}=
\begin{cases} 
1\,,&\text{if}\:a=b\,,
\\
0\,,&\text{if}\:a\neq b\,.
\end{cases}
$$

(iii) The span of $\{\psi_a\}_{a\in\N}$ is dense in $L^2(\d M)$.

(iv) We have $\lim_{a\to\infty}\lambda_a=\infty$.
\end{proposition}
\bp
Since we are assuming $H_{g_0}>0$, for each $f\in L^2(\d M)$ we can define $T(f)=u$, where $u\in H^{1}(M)$ is the unique solution of 
\begin{equation}\notag
\begin{cases} 
\Delta_{g_0}u=0\,,&\text{in}\:M\,,
\\
\critbordo\frac{\d}{\d\eta_{g_0}}u
-\cmz u=fu_{\infty}^{\frac{2}{n-2}}\,,&\text{on}\:\d M\,.
\end{cases} 
\end{equation}
Since $H^{1}(M)$ is compactly embedded in $L^2(\d M)$, the operator $T:L^2(\d M)\to L^2(\d M)$ is compact. Integrating by parts, we see that $T$ is symmetric with respect to the inner product $(\psi_1,\psi_2)\mapsto\int_{\d M}\psi_1\psi_2u_{\infty}^{\frac{2}{n-2}}\ds_{g_0}$. Then the result follows from the spectral theorem for compact operators.
\ep

Let $A\subset \N$ be a finite set such that $\l_a>\frac{n}{n-2}\cminfbar$ for all $a\notin A$, and define the projection 
$$
\Gamma(f)
=\sum_{a\notin A}\left(\int_{\d M}\psi_af\ds_{g_0}\right)\psi_a\uinf^{\frac{2}{n-2}}
=f-\sum_{a\in A}\left(\int_{\d M}\psi_af\ds_{g_0}\right)\psi_a\uinf^{\frac{2}{n-2}}\,.
$$
\begin{lemma}\label{Lemma6.2}
For any $1\leq p<\infty$ there exists $C>0$ such that
\ba
\|f\|_{L^p(\d M)}
\leq &\,C\,\Big\|\critbordo\frac{\d f}{\d \eta_{g_0}}-\cmz f+\frac{n}{n-2}\cminfbar u_{\infty}^{\frac{2}{n-2}}f\Big\|_{L^p(\d M)}\notag
\\
&+C\sup_{a\in A}\left|\int_{\d M}u_{\infty}^{\frac{2}{n-2}}\psi_a f\ds_{g_0}\right|\notag
\end{align}
for all $f\in C^2( M)$ satisfying $\Delta_{g_0}f=0$ in $M$.
\end{lemma}
\bp
Set $T(f)=\critbordo\frac{\d}{\d \eta_{g_0}}f-\cmz f+\frac{n}{n-2}\cminfbar u_{\infty}^{\frac{2}{n-2}}f$ on $\d M$.
Suppose the result is not true. Then we can find a sequence of harmonic functions $\{f_j\}$ satisfying
$$
1=\|f_j\|_{L^p(\d M)}\geq j\|T(f_j)\|_{L^p(\d M)}+j\sup_{a\in A}\left|\int_{\d M}u_{\infty}^{\frac{2}{n-2}}\psi_a f_j\ds_{g_0}\right|\,.
$$ 
By \cite[Lemma 3.2]{brendle-boundary},
$$
\|f_j\|_{W^{1,p}(\d M)}\leq C\|T(f_j)\|_{L^{p}(\d M)}+C\|f_j\|_{L^{p}(\d M)}\,,
\:\:\:\:\text{if}\:p>1\,,
$$ 
and by Proposition \ref{green:global} and Corollary \ref{green:corol} we have
$$
\|f_j\|_{W^{\frac{1}{2},1}(\d M)}\leq C\|T(f_j)\|_{L^{1}(\d M)}+C\|f_j\|_{L^{1}(\d M)}\,.
$$ 
It follows from compactness that we can find a function $f$ satisfying
$$
\|f\|_{L^p(\d M)}=1\,,
\:\:\:\:\:\:\:\:\:\sup_{a\in A}\left|\int_{\d M}u_{\infty}^{\frac{2}{n-2}}\psi_a f\ds_{g_0}\right|=0\,,
$$
and
$$
\int_{\d M}T(\psi_a)f\ds_{g_0}=0\:\:\:\:\:\:\text{for any}\:a\in \N\,.
$$
Hence,
$$
\left(\l_a-\frac{n}{n-2}\cminfbar\right)\int_{\d M}\psi_afu_{\infty}^{\frac{2}{n-2}}\ds_{g_0}=0\:\:\:\:\:\text{for all}\:a\in \N\,.
$$
In particular, $\int_{\d M}\psi_afu_{\infty}^{\frac{2}{n-2}}\ds_{g_0}=0$ for all $a\notin A$, which implies $f\equiv 0$ on $\d M$. This contradicts  $\|f\|_{L^p(\d M)}=1$.
\ep
\begin{lemma}\label{Lemma6.3} There exists $C>0$ such that
\ba\label{Lemma6.3:i}
\|f\|_{L^{\frac{n}{n-2}}(\d M)}
\leq &\,C\,\Big\|\Gamma\left(\critbordo\frac{\d f}{\d \eta_{g_0}}-\cmz f+\frac{n}{n-2}\cminfbar u_{\infty}^{\frac{2}{n-2}}f\right)\Big\|_{L^{\frac{n(n-1)}{n^2-2n+2}}(\d M)}
\\
&+C\sup_{a\in A}\left|\int_{\d M}u_{\infty}^{\frac{2}{n-2}}\psi_a f\ds_{g_0}\right|\notag
\end{align}
and
\ba\label{Lemma6.3:ii}
\|f\|_{L^1(\d M)}
\leq &\, C\,\Big\|\Gamma\left(\critbordo\frac{\d f}{\d \eta_{g_0}}-\cmz f+\frac{n}{n-2}\cminfbar u_{\infty}^{\frac{2}{n-2}}f\right)\Big\|_{L^1(\d M)}
\\
&+C\sup_{a\in A}\left|\int_{\d M}u_{\infty}^{\frac{2}{n-2}}\psi_a f\ds_{g_0}\right|\,,\notag
\end{align}
for all $f\in C^2(M)$ satisfying $\Delta_{g_0}f=0$ in $M$.
\end{lemma}
\bp
We set $p_0=\frac{n(n-1)}{n^2-2n+2}$ and follow the notation in the proof of Lemma \ref{Lemma6.2}. By \cite[Lemma 3.2]{brendle-boundary},
$$
\|f\|_{W^{1,p_0}(\d M)}\leq C\|T(f)\|_{L^{p_0}(\d M)}+C\|f\|_{L^{p_0}(\d M)}\,.
$$ 
Thus, it follows from Lemma \ref{Lemma6.2} that
\begin{equation}\label{Lemma6.3:1}
\|f\|_{W^{1,p_0}(\d M)}\leq C\|T(f)\|_{L^{p_0}(\d M)}
+C\sup_{a\in A}\left|\int_{\d M}u_{\infty}^{\frac{2}{n-2}}\psi_af\ds_{g_0}\right|\,.
\end{equation}
By the definition of $\Gamma$, we have
\begin{equation}\label{Lemma6.3:2}
T(f)=\Gamma(T(f))
+\sum_{a\in A}\left(\frac{n}{n-2}\cminfbar-\l_a\right)\left\{\int_{\d M}\uinf^{\frac{2}{n-2}}\psi_af\ds_{g_0}\right\}\uinf^{\frac{2}{n-2}}\psi_a\,.
\end{equation}
Hence,
\begin{equation}\label{Lemma6.3:3}
\|T(f)\|_{L^{p_0}(\d M)}\leq \|\Gamma(T(f))\|_{L^{p_0}(\d M)}
+C\sup_{a\in A}\left|\int_{\d M}u_{\infty}^{\frac{2}{n-2}}\psi_af\ds_{g_0}\right|\,.
\end{equation}
Now the estimate (\ref{Lemma6.3:i}) follows from (\ref{Lemma6.3:1}), (\ref{Lemma6.3:3}), and Sobolev inequalities.  

In order to prove (\ref{Lemma6.3:ii}), observe that by Lemma \ref{Lemma6.2} we have
$$
\|f\|_{L^1(\d M)}\leq C\|T(f)\|_{L^1(\d M)}
+C\sup_{a\in A}\left|\int_{\d M}u_{\infty}^{\frac{2}{n-2}}\psi_af\ds_{g_0}\right|\,.
$$  
Now the result follows from  (\ref{Lemma6.3:2}).
\ep
\begin{lemma}\label{Lemma6.4}
There exists $\zetaup>0$ with the following significance: for all $z=(z_1,...,z_a)\in\R^{A}$ with $|z|\leq \zetaup$, there exists a smooth function $\bar{u}_z$ satisfying $\Delta_{g_0}\bar{u}_z=0$ in $M$, 
\begin{equation}\label{Lemma6.4:1}
\int_{\d M}\uinf^{\frac{2}{n-2}}(\bar{u}_z-\uinf)\psi_a\ds_{g_0}=z_a
\:\:\:\:\text{for all}\:a\in A\,,
\end{equation}
and
\begin{equation}\label{Lemma6.4:2}
\Gamma\left(\critbordo\frac{\d \bar{u}_z}{\d\eta_{g_0}}-\cmz\bar{u}_z+\cminfbar\bar{u}_z^{\frac{n}{n-2}}\right)=0\,.
\end{equation}
Moreover, the mapping $z\mapsto \bar{u}_z$ is real analytic. 
\end{lemma}
\bp
This is just an application of the implicit function theorem.
\ep
\begin{lemma}\label{Lemma6.5}
There exists $0<\gamma<1$ such that
$$
\big{|}E(\bar{u}_z)-E(\uinf)\big{|}
\leq
C\sup_{a\in A}\left|
\int_{\d M}\psi_a\left(\critbordo\frac{\d \bar{u}_z}{\d\eta_{g_0}}-\cmz\,\bar{u}_z+\cminfbar\,\bar{u}_z^{\frac{n}{n-2}}\ds_{g_0}
\right)
\right|^{1+\gamma}\,,
$$
if $|z|$ is sufficiently small.
\end{lemma}
\bp
Observe that the function $z\mapsto E(\bar{u}_z)$ is real analytic. According to results of Lojasiewicz (see (2.4) in \cite[p.538]{simon}), there exists $0<\gamma<1$ such that
$$
|E(\bar{u}_z)-E(\uinf)|\leq \sup_{a\in A}\left|\frac{\d}{\d z_a}E(\bar{u}_z)\right|^{1+\gamma}\,,
$$
if $|z|$ is sufficiently small.
Differentiating $E(\bar{u}_z)$, we obtain
\ba\label{Lemma6.5:1}
\left(\int_{\d M}\bar{u}_z^{\critbordo}\right)^{\frac{n-2}{n-1}}&\frac{\d}{\d z_a}E(\bar{u}_z)
\\
&=-2(F(\bar{u}_z)-2\cminfbar)\int_{\d M}\bar{u}_z^{\frac{n}{n-2}}\tilde{\psi}_{a,z}\,\ds_{g_0}\notag
\\
&\hspace{0.3cm}-2\int_{\d M}\left(\frac{4(n-1)}{n-2}\frac{\d \bar{u}_z}{\d\eta_{g_0}}-2\cmz\bar{u}_z+2\cminfbar\bar{u}_z^{\frac{n}{n-2}}\right)\tilde{\psi}_{a,z}\,\ds_{g_0}\,,\notag
\end{align}
where we have set $\tilde{\psi}_{a,z}=\displaystyle{\frac{\d\bar{u}_z}{\d z_a}}$ for $a\in A$. Differentiating  (\ref{Lemma6.4:1}), we obtain
\ba\label{Lemma6.5:2}
\int_{\d M}\uinf^{\frac{2}{n-2}}\tilde{\psi}_{a,z}\psi_b\,\ds_{g_0}
=
\begin{cases}
1\,, &\text{if}\: a=b\,,
\\
0\,, &\text{if}\:a\neq b\,,
\end{cases}
\end{align}
for all $b\in A$.

Integrating by parts and using the identity (\ref{Lemma6.4:2}), we see that
\ba\label{Lemma6.5:3}
(F(\bar{u}_z)-2\cminfbar)&\int_{\d M}\bar{u}_z^{\critbordo}\ds_{g_0}
\\
&=-\int_{\d M}\left(\frac{4(n-1)}{n-2}\frac{\d \bar{u}_z}{\d\eta_{g_0}}-2\cmz\bar{u}_z+2\cminfbar\bar{u}_z^{\frac{n}{n-2}}\right)\bar{u}_z\,\ds_{g_0}\notag
\\
&=-\sum_{b\in A}\int_{\d M}\left(\frac{4(n-1)}{n-2}\frac{\d \bar{u}_z}{\d\eta_{g_0}}-2\cmz\bar{u}_z+2\cminfbar\bar{u}_z^{\frac{n}{n-2}}\right)\psi_b\,\ds_{g_0}\notag
\\
&\hspace{2cm}\cdot\int_{\d M}\uinf^{\frac{2}{n-2}}\psi_b\,\bar{u}_z\,\ds_{g_0}\,.\notag
\end{align}

Substituting (\ref{Lemma6.5:2}) and (\ref{Lemma6.5:3}) in (\ref{Lemma6.5:1}), we obtain
\ba
&\left(\int_{\d M}\bar{u}_z^{\critbordo}\ds_{g_0}\right)^{\frac{n-2}{n-1}}\frac{\d}{\d z_a}E(\bar{u}_z)\notag
\\
&\hspace{1cm}=2\sum_{b\in A}\left(\int_{\d M}\bar{u}_z^{\critbordo}\ds_{g_0}\right)^{-1}
\cdot\int_{\d M}\bar{u}_z^{\frac{n}{n-2}}\tilde{\psi}_{a,z}\,\ds_{g_0}
\cdot\int_{\d M}\uinf^{\frac{2}{n-2}}\bar{u}_z\psi_b\,\ds_{g_0}\notag
\\
&\hspace{3cm}\cdot
\int_{\d M}\left(\frac{4(n-1)}{n-2}\frac{\d \bar{u}_z}{\d\eta_{g_0}}-2\cmz\bar{u}_z+2\cminfbar\bar{u}_z^{\frac{n}{n-2}}\right)\psi_b\,\ds_{g_0}\notag
\\
&\hspace{1cm}-2\int_{\d M}\left(\frac{4(n-1)}{n-2}\frac{\d \bar{u}_z}{\d\eta_{g_0}}-2\cmz\bar{u}_z+2\cminfbar\bar{u}_z^{\frac{n}{n-2}}\right)\psi_a\,\ds_{g_0}\,.\notag
\end{align}
Hence, there exists $C>0$ such that
$$
\left|\frac{\d}{\d z_a}E(\bar{u}_z)\right|
\leq C\sup_{a\in A}\left|\int_{\d M}\left(\frac{4(n-1)}{n-2}\frac{\d \bar{u}_z}{\d\eta_{g_0}}-2\cmz\bar{u}_z+2\cminfbar\bar{u}_z^{\frac{n}{n-2}}\right)\psi_a\,\ds_{g_0}\right|\,,
$$
from which the assertion follows.
\ep

We set 
\ba
\mathcal{A}_{\nu}=\Big{\{}(z, (x_k,\e_k,\a_k)_{k=1,...,m})\in &\,\R^A\times (\d M\times\R_+\times\R_+)^m\,,\:\text{such that}\notag
\\
&|z|\leq\zetaup\,,\: d_{g_0}(x_{k}, x^*_{k,\nu})\leq \e^*_{k,\nu}\,,\:\frac{1}{2}\leq\frac{\e_k}{\e^*_{k,\nu}}\leq 2
\,,\:\frac{1}{2}\leq\a_k\leq 2\Big{\}}\,.\notag
\end{align}
For each $\nu$, we can choose a pair $(z_{\nu}, (x_{k,\nu},\e_{k,\nu},\a_{k,\nu})_{k=1,...,m})\in\mathcal{A}_{\nu}$ such that
\ba
&\int_M\frac{2(n-1)}{n-2}\big| d(u_{\nu}-\bar{u}_{z_{\nu}}-\sum_{k=1}^{m}\a_{k,\nu}\uknu )\big|_{g_0}^2\dv_{g_0}\notag
\\
&\hspace{1cm}+\int_{\d M}\cmedia_{g_0}\big(u_{\nu}-\bar{u}_{z_{\nu}}-\sum_{k=1}^{m}\a_{k,\nu}\uknu \big)^2\ds_{g_0}\notag
\\
&\hspace{0.5cm}\leq\int_M\frac{2(n-1)}{n-2}\big{|}d(u_{\nu}-\bar{u}_{z}-\sum_{k=1}^{m}\a_{k}\uk)\big{|}_{g_0}^2\dv_{g_0}\notag
\\
&\hspace{1cm}+\int_{\d M}\cmedia_{g_0}\big(u_{\nu}-\bar{u}_{z}-\sum_{k=1}^{m}\a_{k}\uk\big)^2\ds_{g_0}\notag
\end{align}
for all $(z, (x_k,\e_k,\a_k)_{k=1,...,m})\in\mathcal{A}_{\nu}$.

\vspace{0.2cm}
The proofs of the next two propositions are the same of Propositions 6.6 and 6.7 in \cite{brendle-flow}:
\begin{proposition}\label{Propo6.6}
We have:

(i) For all $i\neq j$, 
$$
\lim_{\nu\to\infty}\left\{
\frac{\e_{i,\nu}}{\e_{j,\nu}}+\frac{\e_{j,\nu}}{\e_{i,\nu}}
+\frac{d_{g_0}(x_{i,\nu},x_{j,\nu})^2}{\e_{i,\nu}\e_{j,\nu}}\right\}=\infty\,.
$$

(ii) We have
$$
\lim_{\nu\to\infty}
\big{\|}\u-\bar{u}_{z_{\nu}}-\sum_{k=1}^{m}\a_{k,\nu}\uknu\big{\|}_{H^1(M)}=0\,.
$$
\end{proposition}
\begin{proposition}\label{Propo6.7}
We have $|z_{\nu}|=o(1)$,
and
$$
d_{g_0}(x_{k,\nu}, x^*_{k,\nu})\leq o(1)\e^*_{k,\nu}\,,\:\:\:
\frac{\e_{k,\nu}}{\e^*_{k,\nu}}=1+o(1)\,,\:\:\:
\text{and}\:\:\:\a_{k,\nu}=1+o(1)\,,
$$
for all $k=1,...,m$. In particular, $(z_{\nu},(x_{k, \nu},\e_{k, \nu},\a_{k, \nu})_{k=1,...,m})$ is an interior point of $\mathcal{A}_{\nu}$ for $\nu$ sufficiently large.
\end{proposition}
\begin{convention}
Assume that $\e_{i,\nu}\leq \e_{j,\nu}$ for all $i\leq j$, without loss of generality.
\end{convention}
\begin{notation}
We write $u_{\nu}=v_{\nu}+w_{\nu}$, where
\begin{equation}\label{def:v:w:2}
v_{\nu}=\bar{u}_{z_{\nu}}+\sum_{k=1}^{m}\a_{k,\nu}\uknu\:\:\:\:\text{and}\:\:\:\:
w_{\nu}=u_{\nu}-\bar{u}_{z_{\nu}}-\sum_{k=1}^{m}\a_{k,\nu}\uknu\,.
\end{equation}
\end{notation}
Observe that by Proposition \ref{Propo6.6} we have
\begin{equation}\label{propr:w:2}
\int_M\frac{2(n-1)}{n-2}|dw_{\nu}|_{g_0}^2\dv_{g_0}
+\int_{\d M}\cmedia_{g_0}w_{\nu}^2\ds_{g_0}=o(1)\,.
\end{equation}
\begin{proposition}\label{Propo6.8}
Let $\psi_{k,\nu}:B^+_{2\rho}(0)\to M$ be Fermi coordinates centered at $x_{k,\nu}$. If we set 
$$
C_{\nu}=\left(\int_{\d M}|\w|^{\critbordo}\ds_{g_0}\right)^{\frac{n-2}{2(n-1)}}
+\left(\int_{M}|\w|^{\crit}\dv_{g_0}\right)^{\frac{n-2}{2n}}\,,
$$ 
then for all $k=1,...,m$, and $a\in A$ we have:

(i) $\big{|} \int_{\d M} \uinf^{\frac{2}{n-2}}\psi_a\w\,\ds_{g_0}\big{|}
\leq o(1)\int_{\d M}|\w|\ds_{g_0}$\,.

(ii) $\displaystyle{\big{|} \int_{\d M} \uknu^{\frac{n}{n-2}}\,\w\,\ds_{g_0}\big{|}
\leq o(1)C_{\nu}}$\,.

(iii) $\displaystyle{\big{|} \int_{\psi_{k,\nu}(\d 'B^+_{2\rho}(0))} \uknu^{\frac{n}{n-2}}
\frac{\e_{k,\nu}^2-|\psi_{k,\nu}^{-1}(x)|^2}{\e_{k,\nu}^2+|\psi_{k,\nu}^{-1}(x)|^2}\,\w\,\ds_{g_0}\big{|}
\leq o(1)C_{\nu}}$\,.

(iv) $\displaystyle{\big{|} \int_{\psi_{k,\nu}(\d 'B^+_{2\rho}(0))} \uknu^{\frac{n}{n-2}}
\frac{\e_{k,\nu}\psi_{k,\nu}^{-1}(x)}{\e_{k,\nu}^2+|\psi_{k,\nu}^{-1}(x)|^2}\,\w\,\ds_{g_0}\big{|}
\leq o(1)C_{\nu}}$\,.
\end{proposition}
\bp
(i) Set $\tilde{\psi}_{a,z}=\frac{\d}{\d z_a}\bar{u}_z$. It follows from the identities (\ref{Lemma6.4:2}) and (\ref{Lemma6.5:2}) that $\tilde{\psi}_{a,0}=\psi_a$ for all $a\in A$.

By the definition of $(z_{\nu}, (x_{k,\nu}, \e_{k,\nu}, \a_{k,\nu})_{1\leq k\leq m})$, we have
$$
\int_M\critbordo<d\tilde{\psi}_{a,z_{\nu}},w_{\nu}>_{g_0}\dv_{g_0}
+\int_{\d M}\cmz\tilde{\psi}_{a,z_{\nu}}w_{\nu}\,\ds_{g_0}=0\,.
$$
Hence,
\ba
\l_a\int_{\d M}&\uinf^{\frac{2}{n-2}}\psi_aw_{\nu}\,\ds_{g_0}\notag
\\
&=-\int_{\d M}\left(\critbordo\frac{\d\psi_a}{\d\eta_{g_0}}-\cmz\psi_a\right)w_{\nu}\,\ds_{g_0}\notag
\\
&=\int_{\d M}\left(\critbordo\frac{\d}{\d\eta_{g_0}}(\tilde{\psi}_{a,z_{\nu}}-\psi_a)
-\cmz(\tilde{\psi}_{a,z_{\nu}}-\psi_a)\right)w_{\nu}\,\ds_{g_0}\,.\notag
\end{align}
Then, since $\l_a>0$ and $|z_{\nu}|\to 0$ as $\nu\to\infty$, we conclude that 
\begin{equation}\label{Propo6.8:1}
\left|\int_{\d M} \uinf^{\frac{2}{n-2}}\psi_a\w\,\ds_{g_0}\right|
\leq o(1)\|\w\|_{L^1(\d M)}\,,
\:\:\:\:\:\:\text{for all}\:a\in A\,,
\end{equation}
from which the assertion (i) follows.

The proofs of (ii), (iii), and (iv) are similar to Proposition \ref{Propo5.3}.
\ep
\begin{proposition}\label{Propo6.9}
There exists $c>0$ such that 
\ba
\frac{n}{n-2}&\cminfbar\int_{\d M}\left(\uinf^{\frac{2}{n-2}}+\sum_{k=1}^{m}\uknu^{\frac{2}{n-2}}\right)\w^2\,\ds_{g_0}\notag
\\
&\leq
(1-c)\left\{\int_{M}\frac{2(n-1)}{n-2}|d\w|_{g_0}^2\dv_{g_0}
+\int_{\d M}\cmedia_{g_0}\w^2\,\ds_{g_0}\right\}\notag
\end{align}
for all $\nu$ sufficiently large.
\end{proposition}
\bp
Suppose by contradiction this is not true. Upon rescaling, we can find a sequence $\{\tilde{w}_{\nu}\}$ satisfying
\begin{equation}\notag
\int_M \frac{2(n-1)}{n-2}|d\tilde{w}_{\nu}|_{g_0}^2\dv_{g_0}
+\int_{\d M}\cmedia_{g_0}\tilde{w}_{\nu}^2\,\ds_{g_0}=1
\end{equation}
and
\begin{equation}\notag
\lim_{\nu\to\infty}\frac{n}{n-2}\cminfbar\int_{\d M}\left(\uinf^{\frac{2}{n-2}}+\sum_{k=1}^{m}\uknu^{\frac{2}{n-2}}\right)\tilde{w}_{\nu}^2\,\ds_{g_0}
\geq 1\,.
\end{equation}
Proceeding as in the proof of Proposition \ref{Propo5.4} and using the same notations, we only have two possibilities:
\\\\
{\underline{Case 1.}} We can suppose that 
\begin{equation}\label{Propo6.9:1}
\lim_{\nu\to\infty}\int_{\d M}\uinf^{\frac{2}{n-2}}\tilde{w}_{\nu}^2\,\ds_{g_0}>0
\end{equation}
and
\ba\label{Propo6.9:2}
\lim_{\nu\to\infty}\Big\{\int_{M\backslash\bigcup_{k=1}^{m}\Omega_{k,\nu}} \frac{2(n-1)}{n-2}|d\tilde{w}_{\nu}|_{g_0}^2\dv_{g_0}
&+\int_{\d M\backslash\bigcup_{k=1}^{m}\Omega_{k,\nu}}\cmedia_{g_0}\tilde{w}_{\nu}^2\,\ds_{g_0}\Big\}
\\
&\leq 
\lim_{\nu\to\infty}\frac{n}{n-2}\cminfbar\int_{\d M}\uinf^{\frac{2}{n-2}}\tilde{w}_{\nu}^2\,\ds_{g_0}\,.\notag
\end{align}

In this case, we can assume that $\tilde{w}_{\nu}\rightharpoonup \tilde{w}$ in $H^1(M)$ and, in view of (\ref{Propo6.9:1}) and (\ref{Propo6.9:2}), we have
\begin{equation}\label{Propo6.9:3}
\int_{\d M}\uinf^{\frac{2}{n-2}}\tilde{w}^2\,\ds_{g_0}>0
\end{equation}
and
$$
\int_{M} \frac{2(n-1)}{n-2}|d\tilde{w}|_{g_0}^2\dv_{g_0}
+\int_{\d M}\cmedia_{g_0}\tilde{w}^2\,\ds_{g_0}
\leq 
\frac{n}{n-2}\cminfbar\int_{\d M}\uinf^{\frac{2}{n-2}}\tilde{w}^2\,\ds_{g_0}\,.
$$
Then it follows from the definition of $\{\psi_a\}_{a\in\mathbb{N}}$ that
\begin{equation}\label{Propo6.9:4}
\sum_{a\in\mathbb{N}}\l_a\left(\int_{\d M}\uinf^{\frac{2}{n-2}}\psi_a\tilde{w}\,\ds_{g_0}\right)^2
\leq
\sum_{a\in\mathbb{N}}\frac{n}{n-2}\cminfbar\left(\int_{\d M}\uinf^{\frac{2}{n-2}}\psi_a\tilde{w}\,\ds_{g_0}\right)^2\,.
\end{equation}
By Proposition \ref{Propo6.8}, we have
$$
\int_{\d M}\uinf^{\frac{2}{n-2}}\psi_a\tilde{w}\,\ds_{g_0}=0\,,
\:\:\:\:\:\text{for all}\:\:a\in A\,.
$$
This, together with (\ref{Propo6.9:4}), implies that $\tilde{w}\equiv 0$ on $\d M$ and contradicts the inequality (\ref{Propo6.9:3}).
\\\\
{\underline{Case 2.}} There exists $1\leq k\leq m$ such that 
$$
\lim_{\nu\to\infty}\int_{\d M}\uknu^{\frac{2}{n-2}}\tilde{w}_{\nu}^2\,\ds_{g_0}>0
$$
and
\ba
\lim_{\nu\to\infty}\Big\{\int_{\Omega_{k,\nu}} \frac{2(n-1)}{n-2}|d\tilde{w}_{\nu}|_{g_0}^2\dv_{g_0}
&+\int_{\Omega_{k,\nu}\cap\d M}\cmedia_{g_0}\tilde{w}_{\nu}^2\,\ds_{g_0}\Big\}\notag
\\
&\leq 
\lim_{\nu\to\infty}\frac{n}{n-2}\cminfbar\int_{\d M}\uknu^{\frac{2}{n-2}}\tilde{w}_{\nu}^2\,\ds_{g_0}\,.\notag
\end{align}

In this case, we proceed exactly as in the proof of Proposition \ref{Propo5.4} to reach a contradiction.

This finishes the proof.
\ep
\begin{corollary}\label{Corol6.10}
For all $\nu$ sufficiently large we have
$$
\frac{n}{n-2}\cminfbar\int_{\d M}v_{\nu}^{\frac{2}{n-2}}w_{\nu}^2\,\ds_{g_0}
\leq
(1-c)\left\{\int_{M}\frac{2(n-1)}{n-2}|d\w|_{g_0}^2\dv_{g_0}
+\int_{\d M}\cmedia_{g_0}\w^2\,\ds_{g_0}\right\}\,.
$$
\end{corollary}
\bp
By the definition of $v_{\nu}$ (see (\ref{def:v:w:2})), we have
$$
\lim_{\nu\to\infty}\int_{\d M}\big|v_{\nu}^{\frac{2}{n-2}}-\uinf^{\frac{2}{n-2}}-\sum_{k=1}^{m}\bar{u}_{(x_{k,\nu}, \e_{k,\nu})}^{\frac{2}{n-2}}\big|^{n-1}\ds_{g_0}=0\,.
$$
Hence, the assertion follows from Proposition \ref{Propo6.9}.
\ep

The next two propositions are similar to Propositions 6.14 and 6.15 of \cite{brendle-flow} and we will just outline their proofs.
\begin{proposition}\label{Propo6.14}
There exist $C>0$ and $0<\gamma<1$ such that
\ba
E(\bar{u}_{z_{\nu}})&-E(\uinf)\notag
\\
&\leq
C\left\{
\int_{\d M}u_{\nu}^{\critbordo}|\cmedia_{g_{\nu}}-\cminfbar|^{\frac{2(n-1)}{n}}
\right\}^{\frac{n}{2(n-1)}(1+\gamma)}
+C\sum_{k=1}^{m}\e_{k,\nu}^{\frac{n-2}{2}(1+\gamma)}\notag
\end{align}
if $\nu$ is sufficiently large.
\end{proposition}
\bp
As in \cite[Lemmas 6.11 and 6.12]{brendle-flow},
making use of estimates (\ref{Lemma6.3:i})  and (\ref{Lemma6.3:ii}), we can show that there exists $C>0$ such that
\begin{equation}\label{Lemma6.11}
\|u_{\nu}-\bar{u}_{z_{\nu}}\|_{L^{\frac{n}{n-2}}(\d M)}^{\frac{n}{n-2}}
\leq C\|u_{\nu}^{\frac{n}{n-2}}(\cmedia_{g_{\nu}}-\cminfbar)\|_{L^{\frac{2(n-1)}{n}}(\d M)}^{\frac{n}{n-2}}
+C\sum_{k=1}^{m}\e_{k,\nu}^{\frac{n-2}{2}}
\end{equation}
and
\begin{equation}\label{Lemma6.12}
\|u_{\nu}-\bar{u}_{z_{\nu}}\|_{L^1(\d M)}
\leq C\|u_{\nu}^{\frac{n}{n-2}}(\cmedia_{g_{\nu}}-\cminfbar)\|_{L^{\frac{2(n-1)}{n}}(\d M)}
+C\sum_{k=1}^{m}\e_{k,\nu}^{\frac{n-2}{2}}\,,
\end{equation}
for $\nu$ sufficiently large.

We will prove the estimate
\ba\label{Lemma6.13}
\sup_{a\in A}&\left|
\int_{\d M}\psi_a\left(
\critbordo\frac{\d \bar{u}_{z_{\nu}}}{\d \eta_{g_0}}-\cmz\bar{u}_{z_{\nu}}+\cminfbar\bar{u}_{z_{\nu}}^{\frac{n}{n-2}}
\right)\ds_{g_0}
\right|
\\
&\leq
C\left\{
\int_{\d M}u_{\nu}^{\critbordo}|\cmedia_{g_{\nu}}-\cminfbar|^{\frac{2(n-1)}{n}}
\right\}^{\frac{n}{2(n-1)}}
+C\sum_{k=1}^{m}\e_{k,\nu}^{\frac{n-2}{2}}\notag
\end{align}
for $\nu$ is sufficiently large.

Integrating by parts, we obtain
\ba
\int_{\d M}&\psi_a\left(
\critbordo\frac{\d \bar{u}_{z_{\nu}}}{\d \eta_{g_0}}-\cmz\bar{u}_{z_{\nu}}+\cminfbar\bar{u}_{z_{\nu}}^{\frac{n}{n-2}}
\right)\ds_{g_0}\notag
\\
&=\int_{\d M}\psi_a\left(
\critbordo\frac{\d u_{\nu}}{\d \eta_{g_0}}-\cmz u_{\nu}+\cminfbar u_{\nu}^{\frac{n}{n-2}}
\right)\ds_{g_0}\notag
\\
&\hspace{0.5cm}+\l_a\int_{\d M}\uinf^{\frac{2}{n-2}}\psi_a(u_{\nu}-\bar{u}_{z_{\nu}})\,\ds_{g_0}
-\cminfbar\int_{\d M}\psi_a(u_{\nu}^{\frac{n}{n-2}}-\bar{u}_{z_{\nu}}^{\frac{n}{n-2}})\,\ds_{g_0}\,.\notag
\end{align}
Using the fact that $\critbordo\frac{\d}{\d \eta_{g_0}}u_{\nu}-\cmz u_{\nu}-\cminfbar u_{\nu}^{\frac{n}{n-2}}=-(\cmedia_{g_{\nu}}-\cminfbar)u_{\nu}^{\frac{n}{n-2}}$ on $\d M$ and the pointwise estimate
$$
|u_{\nu}^{\frac{n}{n-2}}-\bar{u}_{z_{\nu}}^{\frac{n}{n-2}}|
\leq C\bar{u}_{z_{\nu}}^{\frac{2}{n-2}}|u_{\nu}-\bar{u}_{z_{\nu}}|
+C|u_{\nu}-\bar{u}_{z_{\nu}}|^{\frac{n}{n-2}}\,,
$$
we obtain
\ba
\sup_{a\in A}&\left|
\int_{\d M}\psi_a\left(
\critbordo\frac{\d \bar{u}_{z_{\nu}}}{\d \eta_{g_0}}-\cmz\bar{u}_{z_{\nu}}+\cminfbar\bar{u}_{z_{\nu}}^{\frac{n}{n-2}}
\right)\ds_{g_0}
\right|\notag
\\
&\leq
C\|u_{\nu}^{\frac{n}{n-2}}(\cmedia_{g_{\nu}}-\cminfbar)\|_{L^{\frac{2(n-1)}{n}}(\d M)}
+C\|u_{\nu}-\bar{u}_{z_{\nu}}\|_{L^1(\d M)}
+C\|u_{\nu}-\bar{u}_{z_{\nu}}\|_{L^{\frac{n}{n-2}}(\d M)}^{\frac{n}{n-2}}\,.\notag
\end{align}
Then it follows from Lemmas \ref{Lemma6.11} and \ref{Lemma6.12} that
\ba\label{Lemma6.13:1}
\sup_{a\in A}&\left|
\int_{\d M}\psi_a\left(
\critbordo\frac{\d \bar{u}_{z_{\nu}}}{\d \eta_{g_0}}-\cmz\bar{u}_{z_{\nu}}+\cminfbar\bar{u}_{z_{\nu}}^{\frac{n}{n-2}}
\right)\ds_{g_0}
\right|
\\
&\leq
C\|u_{\nu}^{\frac{n}{n-2}}(\cmedia_{g_{\nu}}-\cminfbar)\|_{L^{\frac{2(n-1)}{n}}(\d M)}^{\frac{n}{n-2}}
+C\|u_{\nu}^{\frac{n}{n-2}}(\cmedia_{g_{\nu}}-\cminfbar)\|_{L^{\frac{2(n-1)}{n}}(\d M)}
+C\sum_{k=1}^{m}\e_{k,\nu}^{\frac{n-2}{2}}\,.\notag
\end{align}

On the other hand, since Corollary \ref{Corol3.2} implies 
$$
\|u_{\nu}^{\frac{n}{n-2}}(\cmedia_{g_{\nu}}-\cminfbar)\|_{L^{\frac{2(n-1)}{n}}(\d M)}
=\left(\int_{\d M}|\cmedia_{g_{\nu}}-\cminfbar|^{\frac{2(n-1)}{n}}\ds_{g_{\nu}}\right)^{\frac{n}{2(n-1)}}
\to 0
$$
as $\nu\to\infty$, we can assume that 
\begin{equation}\label{Lemma6.12:4}
\|u_{\nu}^{\frac{n}{n-2}}(\cmedia_{g_{\nu}}-\cminfbar)\|_{L^{\frac{2(n-1)}{n}}(\d M)}<1\,.
\end{equation}

The estimate (\ref{Lemma6.13}) now follows using the inequality (\ref{Lemma6.12:4}) in (\ref{Lemma6.13:1}).

Proposition \ref{Propo6.14} is a consequence of  Lemma \ref{Lemma6.5} and the estimate (\ref{Lemma6.13}).
\ep
\begin{proposition}\label{Propo6.15}
There exists $c>0$ such that
$$
E(v_{\nu})\leq \left(E(\bar{u}_{z_{\nu}})^{n-1}+\sum_{k=1}^{m}E(\bar{u}_{x_k,\e_{k,\nu}})^{n-1}\right)^{\frac{1}{n-1}}
-c\sum_{k=1}^{m}\e_{k,\nu}^{\frac{n-2}{2}}
$$
if $\nu$ is sufficiently large.
\end{proposition}
\bp
Observe that the inequality
\ba
&\left(F(\bar{u}_{z_{\nu}})^{n-1}\bar{u}_{z_{\nu}}^{\critbordo}+F(\uknu)^{n-1}\uknu^{\critbordo}\right)^{\frac{1}{n-1}}
\bar{u}_{z_{\nu}}\uknu\notag
\\
&\hspace{0.5cm}\geq
F(\bar{u}_{z_{\nu}})\bar{u}_{z_{\nu}}^{\frac{n}{n-2}}\uknu
+c\e_{k,\nu}^{-\frac{n}{2}}1_{\{d_{g_0}(x,x_{k,\nu})\leq \e_{k,\nu}\}}\notag
\end{align}
holds on $\d M$ for any $1\leq k\leq m$.

As in \cite{brendle-flow} we obtain
\ba\label{Propo6.15:4}
\frac{1}{2}E(v_{\nu})&\left(\int_{\d M}v_{\nu}^{\critbordo}\ds_{g_0}\right)^{\frac{n-2}{n-1}}
\\
&\leq
\frac{1}{2}\left(E(\bar{u}_{z_{\nu}})^{n-1}+\sum_{k=1}^{m}E(\uknu)^{n-1}\right)^{\frac{1}{n-1}}\left(\int_{\d M}v_{\nu}^{\critbordo}\ds_{g_0}\right)^{\frac{n-2}{n-1}}\notag
\\
&\hspace{0.5cm}-\sum_{k=1}^{m}2\a_{k,\nu}\int_{\d M}\Big(\critbordo\frac{\d}{\d\eta_{g_0}}\bar{u}_{z_{\nu}}-\cmz\bar{u}_{z_{\nu}}\notag
\\
&\hspace{5cm}+\frac{1}{2}F(\bar{u}_{z_{\nu}})\bar{u}_{z_{\nu}}^{\frac{n}{n-2}}\Big)\uknu\ds_{g_0}\notag
\\
&\hspace{0.5cm}-\sum_{i<j}2\a_{i,\nu}\a_{j,\nu}\int_M\critbordo\Delta_{g_0}\ujnu\,\uinu\dv_{g_0}\notag
\\
&\hspace{0.5cm}-\sum_{i<j}2\a_{i,\nu}\a_{j,\nu}\int_{\d M}\Big(\critbordo\frac{\d}{\d\eta_{g_0}}\ujnu-\cmz\ujnu\notag
\\
&\hspace{5cm}+\frac{1}{2}F(\ujnu)\ujnu^{\frac{n}{n-2}}\Big)\uinu\ds_{g_0}\notag
\\
&\hspace{0.5cm}-c\e_{k,\nu}^{\frac{n-2}{2}}-c\sum_{i<j}\left(\frac{\e_{i,\nu}\e_{j,\nu}}{\e_{j,\nu}^2+d_{g_0}(x_{i,\nu},x_{j,\nu})^2}\right)^{\frac{n-2}{2}}\,.\notag
\end{align}

Since $\frac{1}{2}F(\bar{u}_{z_{\nu}})\to \frac{1}{2}F(\uinf)=\cminfbar$ as $\nu\to\infty$, we have the estimate
\begin{equation}\label{Propo6.15:5}
\int_{\d M}\left|\critbordo\frac{\d}{\d\eta_{g_0}}\bar{u}_{z_{\nu}}-\cmz\bar{u}_{z_{\nu}}
+\frac{1}{2}F(\bar{u}_{z_{\nu}})\bar{u}_{z_{\nu}}^{\frac{n}{n-2}}\right|\uknu\ds_{g_0}
\leq o(1)\e_{k,\nu}^{\frac{n-2}{2}}\,.
\end{equation}

Now the assertion follows from the estimates (\ref{Propo5.6:5}), (\ref{Propo6.15:4}), and (\ref{Propo6.15:5}), choosing $\rho$ small and $\nu$ large.
\ep
\begin{corollary}\label{Corol6.16}
Under the hypothesis of Theorem \ref{main:thm}, 
there exist $C>0$ and $0<\gamma<1$ such that
\ba
E(v_{\nu})\,\leq\, &\big(E(\uinf)^{n-1}+mQ(B^n,\d B)^{n-1}\big)^{\frac{1}{n-1}}\notag
\\
&+C\left(
\int_{\d M}u_{\nu}^{\critbordo}|\cmedia_{g_{\nu}}-\cminfbar|^{\frac{2(n-1)}{n}}\ds_{g_0}
\right)^{\frac{n}{2(n-1)}(1+\gamma)}\,,\notag
\end{align}
if $\nu$ is sufficiently large.
\end{corollary}
\bp
Using Proposition \ref{Propo:energy:test}, we obtain $E(\uknu)\leq \Q$ for all $k=1,...,m$. Then the result follows from  Propositions \ref{Propo6.14} and \ref{Propo6.15}.
\ep


\section{Proof of the main theorem}\label{sec:mainthm}
Let $u(t)$, $t\geq 0$, be the solution of (\ref{eq:evol:u}) obtained in Section \ref{sec:prelim}. The next proposition, which is analogous to Proposition 3.3 of \cite{brendle-flow}, is a crucial step in our argument. The blow-up analysis of Section \ref{sec:blowup} is used in its proof.
\begin{proposition}\label{Propo3.3}
Let $\{t_{\nu}\}_{\nu=1}^{\infty}$ be a sequence such that $\lim_{\nu\to\infty}t_{\nu}=\infty$.  Then we can choose $0<\gamma<1$ and $C>0$ such that, after passing to a subsequence, we have
$$
\overline{\cmedia}_{g(t_{\nu})}-\cminfbar 
\leq
C\left\{\int_{\d M} u(t_{\nu})^{\critbordo}|\cmedia_{g(t_{\nu})}-\cminfbar|
^{\conjbordo}\ds_{g_0}\right\}^{\frac{n}{2(n-1)}(1+\gamma)}
$$
for all $\nu$. 
\end{proposition}
\bp
Set $\u(x)=u(x,t_{\nu})$ and $\g=g(t_{\nu})=\u^{\frac{4}{n-2}}g_0$. 
We consider the non-negative smooth function $\uinf$ obtained in Proposition \ref{Propo4.1} and write $u_{\nu}=v_{\nu}+w_{\nu}$ as in the formula \eqref{def:v:w:1} if $\uinf\equiv 0$, or in the formula \eqref{def:v:w:2} if $\uinf>0$. Then, integrating by parts the equations (\ref{eq:R:H}), we obtain
\ba
\cmvbar
&=\int_M\frac{2(n-1)}{n-2}|dv_{\nu}|_{g_0}^2\dv_{g_0}+\int_{\d M}\cmz v_{\nu}^2\ds_{g_0}\notag
\\
&\hspace{0.5cm}+\int_M\frac{2(n-1)}{n-2}|dw_{\nu}|_{g_0}^2\dv_{g_0}+\int_{\d M}\cmz w_{\nu}^2\ds_{g_0}\notag
\\
&\hspace{0.5cm}+\int_M\frac{4(n-1)}{n-2}<dw_{\nu},dv_{\nu}>_{g_0}\dv_{g_0}+\int_{\d M}2\cmz v_{\nu}w_{\nu}\ds_{g_0}\,.\notag
\end{align}
On the other hand,
\ba
\int_M&\frac{4(n-1)}{n-2}<dw_{\nu},dv_{\nu}>_{g_0}\dv_{g_0}+\int_{\d M}2\cmz v_{\nu}w_{\nu}\ds_{g_0}\notag
\\
&=\int_{\d M}2\cmv u_{\nu}^{\frac{n}{n-2}}w_{\nu}\ds_{g_0}
-\int_M\frac{4(n-1)}{n-2}|w_{\nu}|_{g_0}^2\dv_{g_0}
-\int_{\d M}2\cmz w_{\nu}^2\ds_{g_0}\,.\notag
\end{align}
Hence,
\ba
\cmvbar
&=\int_M\frac{2(n-1)}{n-2}|dv_{\nu}|_{g_0}^2\dv_{g_0}+\int_{\d M}\cmz v_{\nu}^2\ds_{g_0}\notag
\\
&\hspace{0.5cm}-\int_M\frac{2(n-1)}{n-2}|dw_{\nu}|_{g_0}^2\dv_{g_0}-\int_{\d M}\cmz w_{\nu}^2\ds_{g_0}
+2\int_{\d M}\cmv u_{\nu}^{\frac{n}{n-2}}w_{\nu}\ds_{g_0}\,,\notag
\end{align}
which can be written as 
\ba\label{Propo3.3:1}
\cmvbar
&=\frac{1}{2}E(v_{\nu})\left\{\int_{\d M}v_{\nu}^{\critbordo}\ds_{g_0}\right\}^{\frac{n-2}{n-1}}
+2\int_{\d M}(\cmv-\cminfbar)u_{\nu}^{\frac{n}{n-2}}w_{\nu}\ds_{g_0}
\\
&\hspace{0.2cm}-\int_M\frac{2(n-1)}{n-2}|dw_{\nu}|_{g_0}^2\dv_{g_0}
-\int_{\d M}\left\{\cmz w_{\nu}^2-\frac{n}{n-2}\cminfbar v_{\nu}^{\frac{2}{n-2}}w_{\nu}^2\right\}\ds_{g_0}\notag
\\
&\hspace{0.2cm}+\cminfbar\int_{\d M}\left\{-\frac{n}{n-2}v_{\nu}^{\frac{2}{n-2}}w_{\nu}^2
+2(v_{\nu}+w_{\nu})^{\frac{n}{n-2}}w_{\nu}\right\}\ds_{g_0}\,.\notag
\end{align}
We can prove that
$$
\left\{\int_{\d M}v_{\nu}^{\critbordo}\ds_{g_0}\right\}^{\frac{n-2}{n-1}}-1
\leq \frac{n-2}{n-1}\int_{\d M}v_{\nu}^{\critbordo}\ds_{g_0}-\frac{n-2}{n-1}\,.
$$
Thus, it follows from the volume normalization $\int_{\d M}(v_{\nu}+w_{\nu})^{\critbordo}\ds_{g_0}=1$ that
\begin{equation}\label{Propo3.3:2}
\left\{\int_{\d M}v_{\nu}^{\critbordo}\ds_{g_0}\right\}^{\frac{n-2}{n-1}}-1
\leq \int_{\d M}\left\{\frac{n-2}{n-1}v_{\nu}^{\critbordo}
-\frac{n-2}{n-1}(v_{\nu}+w_{\nu})^{\critbordo}\right\}\ds_{g_0}\,.
\end{equation}
Using the inequality (\ref{Propo3.3:2}) in the equation (\ref{Propo3.3:1}), we obtain
\ba\label{Propo3.3:3}
\cmvbar
\leq \cminfbar&+(\frac{1}{2}E(v_{\nu})-\cminfbar)\left\{\int_{\d M}v_{\nu}^{\critbordo}\ds_{g_0}\right\}^{\frac{n-2}{n-1}}
+2\int_{\d M}(\cmv-\cminfbar)u_{\nu}^{\frac{n}{n-2}}w_{\nu}\ds_{g_0}\notag
\\
&+\cminfbar\int_{\d M}\Big\{
\frac{n-2}{n-1}v_{\nu}^{\critbordo}
-\frac{n}{n-2}v_{\nu}^{\frac{2}{n-2}}w_{\nu}^2
\\
&\hspace{2cm}+2(v_{\nu}+w_{\nu})^{\frac{n}{n-2}}w_{\nu}
-\frac{n-2}{n-1}(v_{\nu}+w_{\nu})^{\critbordo}
\Big\}\ds_{g_0}\notag
\\
&-\int_M\frac{2(n-1)}{n-2}|dw_{\nu}|_{g_0}^2\dv_{g_0}
-\int_{\d M}\left\{\cmz w_{\nu}^2-\frac{n}{n-2}\cminfbar v_{\nu}^{\frac{2}{n-2}}w_{\nu}^2\right\}\ds_{g_0}\,.\notag
\end{align}

Now we estimate some terms of the right-hand side of (\ref{Propo3.3:3}). 
By the H\"{o}lder's inequality,
\ba\label{Propo3.3:4}
\int_{\d M}&u_{\nu}^{\frac{n}{n-2}}(\cmv-\cminfbar)w_{\nu}\ds_{g_0}
\\
&\leq\left\{\int_{\d M}u_{\nu}^{\critbordo}|\cmv-\cminfbar|^{\conjbordo}\ds_{g_0}\right\}^{\frac{n}{2(n-1)}}
\cdot\left\{\int_{\d M}|w_{\nu}|^{\critbordo}\ds_{g_0}\right\}^{\frac{n-2}{2(n-1)}}\,.\notag
\end{align}

It follows from Corollaries \ref{Corol5.5} and \ref{Corol6.10} that 
\ba
\int_M&\frac{2(n-1)}{n-2}|dw_{\nu}|_{g_0}^2\dv_{g_0}
+\int_{\d M}\left\{\cmz w_{\nu}^2-\frac{n}{n-2}\cminfbar v_{\nu}^{\frac{n}{n-2}}w_{\nu}^2\right\}\ds_{g_0}\notag
\\
&\geq c\int_M\frac{2(n-1)}{n-2}|dw_{\nu}|_{g_0}^2\dv_{g_0}
+c\int_{\d M}\cmz w_{\nu}^2\ds_{g_0}\,.\notag
\end{align}
Since $Q(M,\d M)>0$, this implies
\ba\label{Propo3.3:5}
\int_M&\frac{2(n-1)}{n-2}|dw_{\nu}|_{g_0}^2\dv_{g_0}
+\int_{\d M}\left\{\cmz w_{\nu}^2-\frac{n}{n-2}\cminfbar v_{\nu}^{\frac{n}{n-2}}w_{\nu}^2\right\}\ds_{g_0}
\\
&\geq c'\left\{\int_{\d M} |w_{\nu}|^{\critbordo}\ds_{g_0}\right\}^{\frac{n-2}{n-1}}\,.\notag
\end{align}
By the pointwise estimate
\ba
&\left|\frac{n-2}{n-1}v_{\nu}^{\critbordo}-\frac{n}{n-2}v_{\nu}^{\frac{2}{n-2}}w_{\nu}^2
+2(v_{\nu}+w_{\nu})^{\frac{n}{n-2}}w_{\nu}-\frac{n-2}{n-1}(v_{\nu}+w_{\nu})^{\critbordo}\right|\notag
\\
&\hspace{1cm}\leq Cv_{\nu}^{\max\{0,\critbordo -3\}}|w_{\nu}|^{\min\{\critbordo,3\}}
+C|w_{\nu}|^{\critbordo}\,,\notag
\end{align}
we have
\ba\label{Propo3.3:6}
\int_{\d M}&\left|\frac{n-2}{n-1}v_{\nu}^{\critbordo}-\frac{n}{n-2}v_{\nu}^{\frac{2}{n-2}}w_{\nu}^2
+2(v_{\nu}+w_{\nu})^{\frac{n}{n-2}}w_{\nu}-\frac{n-2}{n-1}(v_{\nu}+w_{\nu})^{\critbordo}\right|\ds_{g_0}\notag
\\
&\hspace{2cm}\leq C\int_{\d M}v_{\nu}^{\max\{0,\critbordo -3\}}|w_{\nu}|^{\min\{\critbordo,3\}}\ds_{g_0}+C\int_{\d M}|w_{\nu}|^{\critbordo}\ds_{g_0}\notag
\\
&\hspace{2cm}\leq C\left\{\int_{\d M}|w_{\nu}|^{\critbordo}\ds_{g_0}
\right\}^{\frac{n-2}{n-1}\min\{\frac{n-1}{n-2},\frac{3}{2}\}}\,.
\end{align}
Recall that $\|w_{\nu}\|_{L^{\critbordo}(\d M)}=o(1)$ by (\ref{propr:w:1}) and (\ref{propr:w:2}).
Using the estimates (\ref{Propo3.3:4}), (\ref{Propo3.3:5}), and (\ref{Propo3.3:6}) in the inequality (\ref{Propo3.3:3}), we obtain
\ba
\cmvbar 
\leq \cminfbar
&+(\frac{1}{2}E(v_{\nu})-\cminfbar)\left\{\int_{\d M}v_{\nu}^{\critbordo}\ds_{g_0}\right\}^{\frac{n-2}{n-1}}\notag
\\
&+C\left\{\int_{\d M}u_{\nu}^{\critbordo}|\cmv-\cminfbar|^{\conjbordo}\ds_{g_0}\right\}^{\frac{n}{n-1}}\,.\notag
\end{align}
In view of equation (\ref{eq:cminbar}), Corollaries \ref{Corol5.7} and \ref{Corol6.16} imply that
$$
\frac{1}{2}E(v_{\nu})-\cminfbar 
\leq C\left\{\int_{\d M}u_{\nu}^{\critbordo}|\cmv-\cminfbar|^{\conjbordo}\ds_{g_0}
\right\}^{\frac{n}{2(n-1)}(1+\gamma)}\,.
$$
Hence,
$$
\cmvbar\leq\cminfbar 
+C\left\{\int_{\d M}u_{\nu}^{\critbordo}|\cmv-\cminfbar|^{\conjbordo}\ds_{g_0}
\right\}^{\frac{n}{2(n-1)}(1+\gamma)}\,,
$$
where we also used Corollary \ref{Corol3.2} with $p=\frac{2(n-1)}{n}$. 

This finishes the proof of Proposition \ref{Propo3.3}.
\ep
Once we have proved Corollary \ref{Corol3.2} and Proposition \ref{Propo3.3}, the proof of the following result is a simple argument by contradiction as in  \cite[Proposition 3.4]{brendle-flow}.
\begin{proposition}\label{Propo3.4}
There exist $0<\gamma<1$ and $t_0>0$ such that
$$
\cmbar-\cminfbar 
\leq\left\{\int_{\d M} u(t)^{\critbordo}|\cm-\cminfbar|
^{\conjbordo}\ds_{g_0}\right\}^{\frac{n}{2(n-1)}(1+\gamma)}
$$
for all $t\geq t_0$.
\end{proposition}
\begin{corollary}\label{Corol3.4}
There exist $0<\gamma<1$, $C>0$ and $t_1>0$ such that
$$
\cmbar-\cminfbar 
\leq C\left\{\int_{\d M} u(t)^{\critbordo}|\cm-\cmbar|
^{\conjbordo}\ds_{g_0}\right\}^{\frac{n}{2(n-1)}(1+\gamma)}
$$
for all $t\geq t_1$.
\end{corollary}
\bp
It follows from Proposition \ref{Propo3.4} that
\ba
\cmbar-\cminfbar
\leq\, &C\left\{\int_{\d M} u(t)^{\critbordo}|\cm-\cmbar|
^{\conjbordo}\ds
_{g_0}\right\}^{\frac{n}{2(n-1)}(1+\gamma)}\notag
\\
&+C(\cmbar-\cminfbar)^{1+\gamma}\,,\notag
\end{align} 
from which the result follows.
\ep
\begin{proposition}\label{Propo3.5}
There exists $C>0$ such that
$$
\int_{0}^{\infty}\left\{\int_{\d M}u(t)^{\critbordo}(\cm-\cmbar)^2\ds_{g_{0}}\right\}^{\frac{1}{2}}dt
\leq C
$$
for all $t\geq 0$.
\end{proposition}
\bp
By the evolution equation (\ref{eq:evol:Hbar}) and Corollary \ref{Corol3.4}, there exists $C>0$ such that 
\ba
\frac{d}{dt}(\cmbar-\cminfbar)
&=-(n-2)\int_{\d M}(\cm-\cmbar)^2\,u(t)^{\critbordo}\ds_{g_0}\notag
\\
&\leq -(n-2)\left\{
\int_{\d M}\big{|}\cm-\cmbar\big{|}^{\conjbordo}u(t)^{\critbordo}\ds_{g_0}
\right\}^{\frac{n}{n-1}}\notag
\\
&\leq -c(\cmbar-\cminfbar)^{\frac{2}{1+\gamma}}\notag
\end{align}
for $t>0$ sufficiently large. Hence, 
$\frac{d}{dt}(\cmbar-\cminfbar)^{-\frac{1-\gamma}{1+\gamma}}\geq c$, which implies
$$
\cmbar-\cminfbar\leq Ct^{-\frac{1+\gamma}{1-\gamma}}\,.
$$
Then using H\"{o}lder's inequality and the equation (\ref{eq:evol:Hbar}) we obtain
\ba
\int_{T}^{2T}&\left(\int_{\d M}(\cm-\cmbar)^2u(t)^{\critbordo}\ds_{g_0}\right)^{\frac{1}{2}}dt\notag
\\
&\leq
\left(\int_{T}^{2T}dt\right)^{\frac{1}{2}}
\left(\int_{T}^{2T}\int_{\d M}(\cm-\cmbar)^2u(t)^{\critbordo}\ds_{g_0}\,dt\right)^{\frac{1}{2}}\notag
\\
&=\left\{\frac{1}{n-2}T(\overline{\cmedia}_{g(T)}-\overline{\cmedia}_{g(2T)})\right\}^{\frac{1}{2}}
\leq CT^{-\frac{\gamma}{1-\gamma}}\notag
\end{align}
for $T$ sufficiently large. This implies
\ba
\int_{0}^{\infty}&\left(\int_{\d M}(\cm-\cmbar)^2u(t)^{\critbordo}\ds_{g_0}\right)^{\frac{1}{2}}dt\notag
\\
&=\int_{0}^{1}\left(\int_{\d M}(\cm-\cmbar)^2u(t)^{\critbordo}\ds_{g_0}\right)^{\frac{1}{2}}dt\notag
\\
&\hspace{1cm}+\sum_{k=0}^{\infty}\int_{2^k}^{2^{k+1}}
\left(\int_{\d M}(\cm-\cmbar)^2u(t)^{\critbordo}\ds_{g_0}\right)^{\frac{1}{2}}dt\notag
\\
&\leq C\sum_{k=0}^{\infty}2^{-\frac{\gamma}{1-\gamma}k}\leq C\,,\notag
\end{align}
which concludes the proof.
\ep
\begin{proposition}\label{Propo3.6}
Given $\gamma_0>0$, there exists $r>0$ such that
$$
\int_{D_r(x)}u(t)^{\critbordo}\ds_{g_{0}}\leq \gamma_0
$$
for all $x\in \d M$ and all $t\geq 0$.
\end{proposition}
\bp
Let $\gamma_0>0$. Using Proposition \ref{Propo3.5}, we can choose $T>0$ large such that
\begin{equation}\label{Propo3.6:1}
\int_{T}^{\infty}\left\{\int_{\d M}u(t)^{\critbordo}(\cm-\cmbar)^2\ds_{g_{0}}\right\}^{\frac{1}{2}}dt
\leq \frac{\gamma_0}{2(n-1)}\,.
\end{equation}
Then we choose $r>0$ small such that 
\begin{equation}\label{Propo3.6:2}
\int_{D_{r}(x)}u(t)^{\critbordo}\ds_{g_0}\leq \frac{\gamma_0}{2}
\end{equation}
for all $t\in [0,T]$ and all $x\in \d M$. By the second equation of (\ref{eq:evol:u}), we see that
\ba\label{Propo3.6:3}
\int_{D_{r}(x)}u(t)\ds_{g_0}
&-\int_{D_{r}(x)}u(T)^{\critbordo}\ds_{g_0}
\\
&=\int_{T}^{t}\frac{d}{dt}\left\{\int_{D_{r}(x)}u(t)^{\critbordo}\ds_{g_0}\right\}dt\notag
\\
&=-(n-1)\int_{T}^{t}\int_{D_{r}(x)}(\cm-\cmbar)\,u(t)^{\critbordo}\ds_{g_0}\,dt\notag
\\
&\leq 
(n-1)\int_{T}^{\infty}\left\{\int_{\d M}(\cm-\cmbar)^2\,u(t)^{\critbordo}\ds_{g_0}\right\}^{\frac{1}{2}}dt\notag
\end{align}
for all $t\geq T$ and all $x\in\d M$, where we have used the boundary area normalization. 
Now the result follows from the inequalities (\ref{Propo3.6:1}), (\ref{Propo3.6:2}), and (\ref{Propo3.6:3}).
\ep
\begin{proposition}\label{Propo3.7}
There exist $C,c>0$ such that 
\begin{equation}\label{Propo3.7:1}
\sup_{M} u(t)\leq C
\:\:\:\:\text{and}\:\:\:\:
\inf_{M} u(t)\geq c
\,,\:\:\:\:\:\text{for all}\:t\geq 0\,.
\end{equation}
\end{proposition}
\bp
By the estimate (\ref{Propo2.4:3}) and the Sobolev embedding theorems, we can choose $C_1>0$ such that
$$
\int_M u(t)^{\crit}dv_{g_0}\leq C_1\,,
\:\:\:\:\:\text{for all}\:t\geq 0\,.
$$

Fix $n-1<q<p<n$. According to Corollary \ref{Corol3.2} there is $C_2>0$ such that
$$
\int_{\d M}|\cm|^p\ds_{g(t)}\leq C_2\,,
\:\:\:\:\:\text{for all}\:t\geq 0\,.
$$
Set $\gamma_0=\gamma_1^{\frac{p}{p-q}}C_2^{-\frac{q}{p-q}}$, where $\gamma_1$ is the constant obtained in Proposition \ref{PropoA.1}. By Proposition \ref{Propo3.6}, there is $r>0$ such that
$$
\int_{D_r(x)}\ds_{g(t)}\leq \gamma_0\,,\:\:\:\:\:\:\text{for all}\:t\geq 0,\,x\in\d M\,.
$$
Then
$$
\int_{D_r(x)}|\cm|^q\ds_{g(t)}
\leq 
\left\{\int_{D_r(x)}\ds_{g(t)}\right\}^{\frac{p-q}{p}}
\left\{\int_{D_r(x)}|\cm|^p\ds_{g(t)}\right\}^{\frac{q}{p}}
\leq \gamma_1\,.
$$
Hence, the first assertion of (\ref{Propo3.7:1}) follows from Proposition \ref{PropoA.1} and the maximum principle. The second one follows exactly as in the proof of the second estimate of (\ref{Propo2.4:1}). 
\ep


\bp[Proof of Theorem \ref{main:thm}]
Once we have proved Proposition \ref{Propo3.7}, it follows as in \cite{brendle-boundary} p.642 that all higher order derivatives of $u$ are uniformly bounded. The uniqueness of the asymptotic limit of $\cm$ follows from Proposition \ref{Propo3.5}.
\ep




\appendix
\renewcommand{\theequation}{A-\arabic{equation}}
\setcounter{equation}{0}
\renewcommand{\thetheorem}{A-\arabic{theorem}}
\setcounter{theorem}{0}
\section{Some elliptic estimates}

The next proposition is a modification of the arguments in \cite[Theorems 8.17 and 8.18]{gilbarg-trudinger}. We refer the reader to \cite[ Lemma 3.2]{sincich} and \cite[Lemma 3.3]{sincich-thesis} for similar results under boundary conditions. (See also the proof of Lemma A.1 in \cite{han-li}.)
\begin{proposition}\label{Propo:estim:Lp}
Let $(M,g)$ be a Riemannian manifold with boundary $\d M$. Let $q>n-1$ and $h\in L^q(\d M)$ with $\|h\|_{L^q(\d M)}\leq \Lambda$. We fix $r_0>0$ small and, for each $x\in \d M$, consider Fermi coordinates $\psi_x:B^+_{4r_0}(0)\to M$  centered at $x$. Then
\\\\
(a) If $p>1$, there exists $C=C(n,g,p,q,\Lambda)$ such that
$$
\sup_{\psi_x(B^+_r(0))}u
\leq
Cr^{-\frac{n}{p}}\|u\|_{L^p(\psi_x(B^+_{2r}(0)))}
+Cr^{1-\frac{n-1}{q}}\|f\|_{L^q(D_{4r}(x))}
$$
for any $x\in\d M$ and $r<r_0$, and any $0< u\in H^1(M)$ and $f\in L^q(\d M)$ satisfying 
\begin{equation}\notag
\begin{cases}
\Delta_gu\geq 0\,,&\text{in}\: M\,,
\\
\frac{\d}{\d\eta_g}u+hu\geq f\,,&\text{on}\: \d M\,.
\end{cases}
\end{equation}

\bigskip
\noindent
(b)  If $1\leq p< \frac{n}{n-2}$, there exists $C=C(n,g,p,q,\Lambda)$ such that
$$
r^{-\frac{n}{p}}\|u\|_{L^p(\psi_x(B^+_{2r}(0)))}
\leq 
C\inf_{\psi_x(B^+_r(0))}u
+Cr^{1-\frac{n-1}{q}}\|f\|_{L^q(D_{4r}(x))}
$$
for any $x\in\d M$ and $r<r_0$, and any $0< u\in H^1(M)$ and $f\in L^q(\d M)$ satisfying
\begin{equation}\notag
\begin{cases}
\Delta_gu\leq 0\,,&\text{in}\: M\,,
\\
\frac{\d}{\d\eta_g}u+hu\leq f\,,&\text{on}\: \d M\,.
\end{cases}
\end{equation}
\end{proposition}
\begin{remark}
According to our notations, $D_r(x)=\psi_x(\d 'B_r(0))$ (see Section \ref{sec:prelim}).
\end{remark}
\begin{proposition}\label{PropoA.1}
Let $(M^n,g_0)$ be a compact Riemannian manifold with boundary $\d M$ and with dimension $n\geq 3$. Choose $\rho>0$ small such that, for all $x\in\d M$, we have Fermi coordinates $\psi_x:B^+_{2\rho}(0)\to M$ centered at $x$. For each $q>n-1$ and $C_1>0$, we can find constants $\gamma_1=\gamma_1(n,g_0,q,C_1)>0$ and $C=C(n,g_0,q)>0$ with the following significance: if $g=u^{\frac{4}{n-2}}g_0$ is a conformal metric satisfying 
$$
\int_{M}\dv_{g}\leq C_1
\:\:\:\:\:\:
\text{and}
\:\:\:\:\:\:
\int_{D_r(x)}|\cmedia_g|^qd\s_g\leq \gamma_1
$$
for $x\in\d M$ and $0<r<\rho$, then we have
$$
u(x)\leq Cr^{-\frac{n-2}{2}}\left(\int_{\psi_x(B^+_r(0))}\dv_g\right)^{\frac{n-2}{2n}}\,.
$$
\end{proposition}
\bp
Suppose that $C_1=1$. 
We can assume that, for any $x\in\d M$, using Fermi coordinates $\psi_x:B^+_{2\rho}(0)\to M$, we have
\begin{equation}\label{PropoA.1:0}
\frac{1}{\sqrt{2}}|z|\leq d_{g_0}(\psi_x(z),x)\leq \sqrt{2}|z|\,,
\:\:\:\:\:\:\text{for all}\:z\in B^+_{2\rho}(0)\,.
\end{equation}

Given $r\in(0,\rho)$ and $x\in \d M$, we define 
$f(s)=(r-s)^{\frac{n-2}{2}}\sup_{B^+_s(0)}u\circ \psi_x$ for $s\in(0,r]$,
and $f(0)=r^{\frac{n-2}{2}}u(x)$. Then we can choose $r_0\in [0,r)$ satisfying
$f(r_0)\geq f(s)$ for all $s\in [0,r)$,
and $x_0=(x_0^1,...,x_0^n)\in\Rn$, with $|x_0|\leq r_0$, such that 
$$
u\circ\psi_x(x_0)\geq u\circ\psi_x(z)\,,
\:\:\:\:\:\:\text{for all}\: z\in B^+_{r_0}(0)\,.
$$ 

Set $\bar{x}_0=(x_0^1,...,x_0^{n-1},0)$ and choose a $0<s\leq \frac{r-r_0}{2}$. We first assume $s> 8x_0^n$. 
It follows from Proposition \ref{Propo:estim:Lp} that there exists $C=C(n,g_0,q)$ such that
\ba
s^{\frac{n-2}{2}}\sup_{B^+_{s/8}(0)}u\circ\psi_{\psi_x(\bar{x}_0)}
&\leq
C\left\{\int_{\psi_{\psi_x(\bar{x}_0)}(B^+_{s/4}(0))}u^{\crit}\dv_{g_0}\right\}^{\frac{n-2}{2n}}\notag
\\
&+Cs^{\frac{n}{2}-\frac{n-1}{q}}
\left\{\int_{\psi_{\psi_x(\bar{x}_0)}(\d 'B^+_{s/2}(0))}\Big{|}\critbordo\frac{\d u}{\d\eta_{g_0}}-\cmz u\Big{|}^q\ds_{g_0}
\right\}^{\frac{1}{q}}\,,\notag
\end{align}
where $\psi_{\psi_x(\bar{x}_0)}:B^+_{2\rho}(0)\to M$ are Fermi coordinates centered at $\psi_x(\bar{x}_0)$. Then, 
by (\ref{PropoA.1:0}) and the fact that $s/8> x_0^n$, we have
\ba
s^{\frac{n-2}{2}}u\circ\psi_x(x_0)
&\leq
C\left\{\int_{\psi_x(B^+_s(\bar{x}_0))}u^{\crit}\dv_{g_0}\right\}^{\frac{n-2}{2n}}\notag
\\
&+Cs^{\frac{n}{2}-\frac{n-1}{q}}
\left\{\int_{\psi_x(\d 'B^+_s(\bar{x}_0))}\Big{|}\critbordo\frac{\d u}{\d\eta_{g_0}}-\cmz u\Big{|}^q\ds_{g_0}
\right\}^{\frac{1}{q}}\,.\notag
\end{align}
Using the second equation of (\ref{eq:R:H}), we conclude that
\ba
s^{\frac{n-2}{2}}u\circ\psi_x(x_0)
&\leq
C\left\{\int_{\psi_x(B^+_s(\bar{x}_0))}\dv_{g}\right\}^{\frac{n-2}{2n}}\notag
\\
&+Cs^{\frac{n}{2}-\frac{n-1}{q}}
\left\{\int_{\psi_x(\d 'B^+_s(\bar{x}_0))}u^{\frac{n}{n-2}q-\critbordo}|\cmedia_g|^q\ds_{g}
\right\}^{\frac{1}{q}}\notag
\end{align}
holds whenever $8x_0^n< s\leq \frac{r-r_0}{2}$.

On the other hand, by a standard interior estimate for linear elliptic equations (see \cite[Theorem 8.17]{gilbarg-trudinger}), if $s\leq 8x_0^n$ and $s\leq \frac{r-r_0}{2}$ then there exists $C=C(n,g_0)$ such that
$$
s^{\frac{n-2}{2}}u\circ\psi_x(x_0)
\leq
C\left\{\int_{\psi_x(B^+_{r}(0))}u^{\crit}\dv_{g_0}\right\}^{\frac{n-2}{2n}}\,.
$$

By the definitions of $r_0$ and $x_0$, we obtain
$$
\sup_{B^+_{\frac{r-r_0}{2}}(\bar{x}_0)}u\circ\psi_x
\leq \sup_{B^+_{\frac{r+r_0}{2}}(0)}u\circ\psi_x
\leq 2^{\frac{n-2}{2}}u\circ\psi_x(x_0)\,.
$$
Hence, there exists $K=K(n,g_0,q)>0$ such that
\ba\label{PropoA.1:1}
s^{\frac{n-2}{2}}u\circ\psi_x(x_0)
&\leq
K\left\{\int_{\psi_x(B^+_{r}(0))}\dv_g\right\}^{\frac{n-2}{2n}}
\\
&+K\big(s^{\frac{n-2}{2}}u\circ\psi_x(x_0)\big)^{\frac{n}{n-2}-\critbordo\frac{1}{q}}
\left\{\int_{D_{r}(x))}|\cmedia_g|^q\ds_g
\right\}^{\frac{1}{q}}\notag
\end{align}
for all $0<s\leq\frac{r-r_0}{2}$. 

Now we choose $\gamma_1=\gamma_1(n,g_0,q)>0$ such that 
\begin{equation}\label{PropoA.1:2}
(2K)^{\frac{n}{n-2}-\critbordo\frac{1}{q}}\gamma_1^{\frac{1}{q}}\leq \frac{1}{2}
\end{equation}
and claim that, if $\int_{M}\dv_{g}\leq 1$ and $\int_{D_r(x)}|\cmedia_g|^qd\s_g\leq \gamma_1$, then
\begin{equation}\label{PropoA.1:3}
\left(\frac{r-r_0}{2}\right)^{\frac{n-2}{2}}u\circ\psi_x(x_0)\leq 2K\,.
\end{equation}
Indeed, if $\left(\frac{r-r_0}{2}\right)^{\frac{n-2}{2}}u\circ\psi_x(x_0)> 2K$, then we can use  (\ref{PropoA.1:1}) with $s=\left(\frac{2K}{u\circ\psi_x(x_0)}\right)^{\frac{2}{n-2}}<\frac{r-r_0}{2}$ to conclude that $2K\leq K+K(2K)^{{\frac{n}{n-2}-\critbordo\frac{1}{q}}}\gamma_1^{\frac{1}{q}}$, which contradicts our choice of $\gamma_1$.

Using  (\ref{PropoA.1:1}) with $s=\frac{r-r_0}{2}$, and (\ref{PropoA.1:3}), we can see that
\ba
\left(\frac{r-r_0}{2}\right)^{\frac{n-2}{2}}u\circ\psi_x(x_0)
\leq &K\left(\int_{\psi_x(B^+_r(0))}\dv_g\right)^{\frac{n-2}{2n}}\notag
\\
&+\frac{1}{2}(2K)^{{\frac{n}{n-2}-\critbordo\frac{1}{q}}}\gamma_1^{\frac{1}{q}}
\left(\frac{r-r_0}{2}\right)^{\frac{n-2}{2}}u\circ\psi_x(x_0)\,.\notag
\end{align}
Hence, by (\ref{PropoA.1:2}),
$$
\left(\frac{r-r_0}{2}\right)^{\frac{n-2}{2}}u\circ\psi_x(x_0)
\leq 2K\left(\int_{\psi_x(B^+_r(0))}\dv_g\right)^{\frac{n-2}{2n}}\,.
$$
This implies
$$
r^{\frac{n-2}{2}}u(x)
\leq (r-r_0)^{\frac{n-2}{2}}u\circ\psi_x(x_0)
\leq 2^{\frac{n}{2}}K\left(\int_{\psi_x(B^+_r(0))}\dv_g\right)^{\frac{n-2}{2n}}\,,
$$
proving the case $C_1=1$.

Now we turn attention to the case $C_1>1$. Let $\gamma_1$ be the constant obtained above. If $\int_{M}\dv_g\leq C_1$, we choose $\l=C_1^{-\frac{n-2}{2n}}$ and set $\tilde{g}=\l^{\frac{4}{n-2}}g=(\l u)^{\frac{4}{n-2}}g_0$. Then $\int_M\dv_{\tilde{g}}\leq 1$, and we have $\int_{D_r(x)}|\cmedia_{\tilde{g}}|^q\ds_{\tilde{g}}\leq \gamma_1$ whenever  $\int_{D_r(x)}|\cmedia_g|^q\ds_{g}\leq \gamma_1 C_1^{\frac{n-1-q}{n}}$. In this case, we proved above that  
$$
\l u(x)\leq Cr^{-\frac{n-2}{2}}\left(\int_{\psi_x(B^+_r(0))}\dv_{\tilde{g}}\right)^{\frac{n-2}{2n}}\,,
$$
which is equivalent to
$$
u(x)\leq Cr^{-\frac{n-2}{2}}\left(\int_{\psi_x(B^+_r(0))}\dv_g\right)^{\frac{n-2}{2n}}
$$ 
by rescaling. This finishes the proof.
\ep

Using Proposition \ref{Propo:estim:Lp}(b) and interior Harnack estimates for elliptic linear equations (see \cite[Theorem 8.18]{gilbarg-trudinger}), one can prove the next proposition by adapting the arguments in \cite[Proposition A.2]{brendle-flow}.
\begin{proposition}\label{PropoA.2}
Suppose $u>0$ satisfies
\begin{equation}\notag
\begin{cases}
-\D_{g_0}u\geq 0\,,&\text{in}\:M\,,
\\
-\frac{\d}{\d\eta_{g_0}}u+Pu\geq 0\,,&\text{on}\:\d M\,,
\end{cases} 
\end{equation}
where $P\in C^{\infty}(\d M)$.
Then there exists $C=C(M,g_0,P)$ such that
$$\int_{M}u \dv_{g_0}\leq C\inf_{M}u\,.$$
In particular,
$$\int_{M}u^{\crit}\dv_{g_0}\leq C\inf_{M}u\left(\sup_{M}u\right)^{\frac{n+2}{n-2}}\,.$$
\end{proposition}


\renewcommand{\theequation}{B-\arabic{equation}}
\setcounter{equation}{0}
\renewcommand{\thetheorem}{B-\arabic{theorem}}
\setcounter{theorem}{0}
\section{Construction of the Green function}

In this section, we prove the existence of the Green function used in this paper and
some of its properties. The construction performed here is similar to the ones in \cite[p.106]{aubin} and \cite[p.201]{druet-hebey-robert}.
\begin{lemma}\label{lemma:holder}
Let $(M,g)$ be a compact Riemannian manifold of dimension $n\geq 2$ and fix $x\in M$ and $\a\in \R$.
Let $u:M\backslash\{x\}\to \R$ be a function satisfying
$$
|u(y)|\leq C_0d_g(x,y)^{\a}
\:\:\:\:\text{and}\:\:\:\:
|\nabla_gu(y)|_g\leq C_0d_g(x,y)^{\a-1}\,,
$$
for any $y\in M$, with $x\neq y$. Then, for any $0<\theta\leq 1$, there exists $C_1=C_1(M,n,g,C_0,\a)$ such that
$$
|u(y)-u(z)|\leq C_1d_g(y,z)^{\theta}(d_g(x,y)^{\a-\theta}+d_g(x,z)^{\a-\theta})
$$
for any $y,z\in M$, with $y\neq x\neq z$.
\end{lemma}
\bp
Let $y\neq x$ and $z\neq x$.\\\\
{\underline{1st case:}} $d_g(y,z)\leq \frac{1}{2}d_g(x,y)$.
Let $\gamma:[0,1]\to M$ be a smooth curve such that $\gamma(0)=y$, $\gamma(1)=z$, and
$\int_0^1|\gamma'(t)|_gdt\leq\frac{3}{2}d_g(y,z)$.

\vspace{0.2cm}\noindent
{\it{Claim.}} We have
$\frac{1}{4}d_g(x,y)\leq d_g(\gamma(t),x)\leq\frac{7}{4}d_g(x,y)$.

Indeed, since $d_g(y,\gamma(t))\leq\frac{3}{2} d_g(y,z)\leq\frac{3}{4}d_g(x,y)$, we have
$$
d_g(x,\gamma(t))\geq d_g(x,y)-d_g(\gamma(t),y)\geq
d_g(x,y)-\frac{3}{4}d_g(x,y)=\frac{1}{4}d_g(x,y)\,.
$$
Moreover,
$$
d_g(\gamma(t),x)\leq d_g(\gamma(t),y)+d_g(y,x)\leq
\frac{3}{4}d_g(x,y)+d_g(x,y)=\frac{7}{4}d_g(x,y)\,.
$$
This  proves the claim.

Observe that $u(z)-u(y)=\int_0^1g(\nabla_gu(\gamma(t)),\gamma'(t))\,dt$. Thus, 
\ba 
|u(y)-u(z)| 
&\leq \sup_{t\in
[0,1]}|\nabla_gu(\gamma(t))|_g\int_0^1|\gamma'(t)|_gdt\notag
\\
&\leq C\sup_{t\in
[0,1]}d_g(\gamma(t),x)^{\a-1}\frac{3}{2}d_g(y,z)\notag
\\
&\leq C(\a)d_g(x,y)^{\a-1}d_g(y,z)\notag
\\
&\leq C(\a)d_g(x,y)^{\a-\theta}d_g(y,z)^{\theta}\,.\notag
\end{align}
\noindent
{\underline{2nd case:}} $d_g(y,z)> \frac{1}{2}d_g(x,y)$.
In this case, we have
\ba
|u(y)-u(z)|
&\leq |u(y)|+|u(z)|\notag
\\
&\leq Cd_g(y,x)^{\a}+Cd_g(z,x)^{\a}\notag
\\
&\leq
Cd_g(y,x)^{\a-\theta}d_g(z,y)^{\theta}+Cd_g(z,x)^{\a-\theta}(d_g(x,y)+d_g(y,z))^{\theta}\notag
\\
&\leq Cd_g(y,z)^{\theta}(d_g(x,y)^{\a-\theta}+d_g(x,z)^{\a-\theta})\,.\notag
\end{align}
This proves the lemma.
\ep
\begin{notation}
In what follows, $(M,g)$ will denote a compact Riemannian manifold with boundary $\d M$, dimension
$n\geq 3$, and positive Sobolev quotient $Q(M,\d M)$. We denote by $L_g$ the conformal Laplacian
$\Delta_g-\frac{n-2}{4(n-1)}R_g$, and by $B_g$ the boundary conformal operator
$\frac{\d}{\d\eta_g}-\frac{n-2}{2(n-1)}H_g$, where $\eta_g$ is the inward unit normal
vector to $\d M$.
\end{notation}
\begin{proposition}\label{green:point}
Fix $x_0\in\d M$ and assume there exist $C=C(M,n,g)$ and $N>1$ such that
\begin{equation}\label{hyp:H}
H_g(y)\leq Cd_g(x_0,y)^{N}\,,\:\:\:\:\text{for all}\:y\in \d M\,.
\end{equation}
If $N$ is sufficiently large, then there exists a positive $G_{x_0}\in C^{\infty}(M\backslash \{x_0\})$ satisfying
\begin{equation}\label{green:formula:0}
\phi(x_0)=-\int_M G_{x_0}(y)L_g\phi(y)\dv_g(y)-\int_{\d M}G_{x_0}(y)B_g\phi(y)\ds_g(y)
\end{equation}
for any $\phi\in C^2(M)$. Moreover, the following properties hold:
\\\\
(P1) There exists $C=C(M,n,g)$ such that, for any $y\in M$ with $y\neq x_0$,
$$
|G_{x_0}(y)|\leq Cd_g(x_0,y)^{2-n}\:\:\:\:\text{and}\:\:\:\:
|\nabla_g G_{x_0}(y)|_g\leq Cd_g(x_0,y)^{1-n}\,.
$$
\\
(P2) Consider Fermi coordinates $y=(y_1,...,y_n)$ centered at $x_0$. In those coordinates, write $g_{ab}=\exp(h_{ab})$, where $h_{ab}$, $a,b=1,...,n$, is a symmetric 2-tensor of the form
$$
h_{ij}(y)=\sum_{|\a|=1}^{d}h_{ij,\a}y^{\a}+O(|y|^{d+1})\,,
$$
for $i,j=1,...,n-1$, and $h_{an}=0$ for $a=1,...,n$. Here, $d=\left[\frac{n-2}{2}\right]$. Then there exists $C=C(M,n,g)$ such that
\begin{align}
&\Big|G_{x_0}(y)-\frac{|y|^{2-n}}{(n-2)\sigma_{n-1}}\Big|\leq
C\sum_{a,b=1}^{n-1}\sum_{|\a|=1}^{d}|h_{ab,\a}|d_g(x_0,y)^{|\a|+2-n}+Cd_g(x_0,y)^{d+3-n}\,,\notag
\\
&\Big|\nabla_g\Big(
G_{x_0}(y)-\frac{|y|^{2-n}}{(n-2)\sigma_{n-1}}
\Big)\Big|_g\leq
C\sum_{a,b=1}^{n-1}\sum_{|\a|=1}^{d}|h_{ab,\a}|d_g(x_0,y)^{|\a|+1-n}+Cd_g(x_0,y)^{d+2-n}\,.\notag
\end{align}
\end{proposition}
\bp
Firstly, we define an appropriate coordinate system for points near the boundary.
Set $d_x=d(x,\d M)$ for $x\in M$, and $M_{\rho}=\{x\in M\,;\:d_x<\rho\}$ for $\rho>0$.

We choose $\rho_0>0$ small such that the function
\ba
M_{2\rho_0}&\to \d M\notag
\\
x&\mapsto \bar{x}\notag
\end{align}
is well defined and smooth, where $\bar{x}$ is defined by $d(x,\bar{x})=d(x,\d M)$. Then, for any $0<t<2\rho_0$, the set $\d_tM=\{x\in M\,;\:d_x=t\}$ is a smooth embedded $(n-1)$-submanifold of $M$. For each $x\in M_{\rho_0}$, define the function
\ba
M_{2\rho_0}&\to \d_{d_x} M\notag
\\
y&\mapsto y_x\,,\notag
\end{align}
where $y_x$ is defined by $d(y,y_x)=d(y,\d_{d_x} M)$.

For any $x\in M_{\rho_0}$, we define the local coordinate system $\psi_x(y)=(y_1,...,y_n)$ on $M_{2\rho_0}$. Here, $y_n=d_y$, and $(y_1,...,y_{n-1})$ are normal coordinates of $y_x$, centered at $x$, with respect to the submanifold $\d_{d_x}M$ . Then $(x,y)\mapsto \psi_x(y)$ is locally defined and smooth.

Observe that $\psi_x(x)=(0,...,0,d_x)$ for any $x\in M_{\rho_0}$, and that $\psi_x$ are Fermi coordinates for any $x\in\d M$.
Moreover, in the coordinates $\psi_x(y)=(y_1,...,y_n)$ we have $g_{an}\equiv \delta_{an}$ and $g_{ab}(x)=\delta_{ab}$, for $a,b=1,...,n$. It is also clear that $d\psi_x^{-1}(\d/\d y_n)$ is the normal unit vector to $\d M$. Choosing $\rho_0$ possibly smaller,  we can assume that, for any $x\in M_{\rho_0}$, $\psi_x(y)=(y_1,...,y_n)$ is defined for $0\leq y_n <2\rho_0$ and $|(y_1,...,y_{n-1})|<\rho_0$.

Let $\chiup:\R_+\to [0,1]$ be a smooth cutoff function satisfying $\chiup(t)=1$ for $t<\rho_0/2$, and $\chiup(t)=0$ for $t\geq \rho_0$. For each $x\in M_{\rho_0}$, set
\ba
K_1(x,y)=\chiup(y_n/2)&\chiup(|(y_1,...,y_{n-1})|)\notag
\\
&\cdot\left\{|(y_1,...,y_{n-1},y_n-d_x)|^{2-n}+|(y_1,...,y_{n-1},y_n+d_x)|^{2-n}\right\}\,,\notag
\end{align}
where we are using the coordinates $\psi_x(y)=(y_1,...,y_n)$. Observe that
$$
\sum_ {a=1}^n\frac{\d^2}{\d y_ a^2}K_1(x,y)=0 \,,\:\:\:
\text{for}\:\:|(y_1,...,y_{n-1})|<\rho_0/2\,,\:0\leq y_n<\rho_0\,,\:\text{and}\:x\neq y\,.
$$
Moreover,  $\d K_1/\d y_n (x,y)=0$ if $y\in \d M$ with $x\neq y$.

For each $x\in M\backslash M_{\rho_0/2}$, set
$$
K_2(x,y)=\chiup(4d_g(y,x))d_g(y,x)^{2-n}\,,\:\:\:\:\text{if}\:\:0<d_g(y,x)<\rho_0/4\,,
$$
and $0$ otherwise. We can assume that $\rho_0/4$ is smaller than the injectivity radius of $(M,g)$.  If we express $y\mapsto K_2(x,y)$ in normal coordinates $(y_1,...,y_n)$ centered at $x$, we have $K_2(x,y)=\chiup(4|(y_1,...,y_n)|)|(y_1,...,y_n)|^{2-n}$, and thus
$$
\sum_ {a=1}^n\frac{\d^2}{\d y_ a^2}K_2(x,y)=0 \,,\:\:\:\:
\text{for}\:\:0<d_g(y,x)<\rho_0/8\,.
$$

Define $K:M\times M\backslash D_M\to \R$ by the expression
$$
K(x,y)=\frac{1}{(n-2)\sigma_{n-1}}\chiup(d_x)K_1(x,y)+\frac{1}{(n-2)\sigma_{n-1}}(1-\chiup(d_x))K_2(x,y)\,.
$$
Here,  $D_M=\{(x,x)\in M\times M\,;\:x\in M\}$.
Thus, $K(x,y)=K_1(x,y)$ if $x\in M_{\rho_0/2}$, and $K(x,y)=K_2(x,y)$ if $x\in M\backslash M_{\rho_0}$.
 Observe that $\d K/\d\eta_{g,y}(x,y)=0$ if $y\in\d M$ with $y\neq x$.

Expressing  $y\mapsto K_1(x,y)$ and $y\mapsto K_2(x,y)$ in their respective coordinate systems (as
described above)  one can check that there exists  $C=C(M,g,n)$ such that
$$
|L_{g,y}K(x,y)|\leq Cd_g(x,y)^{1-n}\,.
$$

For any $\phi\in C^2(M)$ and $x\in M$, we have
\begin{align}\label{form:green:H}
\phi(x)&=\int_M\big( \Delta_{g,y}K(x,y)\phi(y)-K(x,y)\Delta_g\phi(y)\Big)\dv_g(y)\notag
\\
&\hspace{1cm}-\int_{\d M} K(x,y)\frac{\d}{\d\eta_g}\phi(y)\ds_g(y)\,.
\end{align}
Indeed, this expression holds for $\frac{1}{(n-2)\sigma_{n-1}}K_1$ when $x\in M_{\rho_0}$, and for $\frac{1}{(n-2)\sigma_{n-1}}K_2$ when $x\in M\backslash M_{\rho_0/2}$.
In particular, $\Delta_{distr, y}K(x,y)=\Delta_{g,y}K(x,y)-\delta_x$.

We define $\Gamma_{k}:M\times M\backslash D_M\to \R$ inductively by setting
$$
\Gamma_1(x,y)=L_{g,y}K(x,y)
$$
and
$$
\Gamma_{k+1}(x,y)=\int_M\Gamma_{k}(x,z)\Gamma_1(z,y)\dv_g(z)\,.
$$
According to  \cite[Proposition 4.12]{aubin}, which is a result due to Giraud (\cite[p.50]{giraud}), we have
\begin{equation}\label{estim:gamma}
|\Gamma_k(x,y)|\leq
\begin{cases}
Cd_g(x,y)^{k-n}\,,\:\:&\text{if}\:k<n\,,
\\
C(1+|\log d_g(x,y)|)\,,\:\:&\text{if}\:k=n\,,
\\
C\,,\:\:&\text{if}\:k>n\,,
\end{cases}
\end{equation}
for some $C=C(M,g,n)$.
Moreover, $\Gamma_{k}$ is continuous on $M\times M$ for $k>n$, and on $M\times M\backslash D_M$ for $k\leq n$.

Now we will refine the estimate (\ref{estim:gamma}) around the point $x_0\in \d M$, using the expansion $g_{ab}=\exp(h_{ab})$. 
Since $K(x,y)=K_1(x,y)$ for $x\in \d M$, one can see that
$$
|L_{g,y}K(x_0,y)|\leq C\sum_{a,b=1}^n\sum_{|\a|=1}^{d}|h_{ab,\a}|d_g(x_0,y)^{|\a|-n}+Cd_g(x_0,y)^{d+1-n}\,.
$$
Then Giraud's result implies
\begin{equation}\label{estim:gamma:h}
|\Gamma_k(x_0,y)|\leq  C
\sum_{a,b=1}^n\sum_{|\a|=1}^{d}|h_{ab,\a}|d_g(x_0,y)^{k-1+|\a|-n}+d_g(x_0,y)^{k+d-n}\,,\:\:\text{if}\:k<n-d\,.
\end{equation}

\vspace{0.2cm}\noindent
{\it{Claim 1.}} Given $0<\theta<1$, there exists $C=C(M,g,n,\theta)$ such that
\begin{equation}\label{estim:claim1}
|\Gamma_{n+1}(x,y)-\Gamma_{n+1}(x,y')|\leq Cd_g(y,y')^{\theta}\,, \:\:\:\: \text{for any}\:y\neq
x\neq y'\,.
\end{equation}
In particular, $\Gamma_{n+1}(x_0,\cdot)\in C^{0,\theta}(M )$.

Indeed, observe that
$
|\Gamma_1(x,y)-\Gamma_1(x,y')|\leq Cd_g(y,y')^{\theta}(d_g(x,y)^{1-\theta-n}+d_g(x,y')^{1-\theta-n})\,,
$
according to Lemma \ref{lemma:holder}. So, Claim 1 follows from the estimates (\ref{estim:gamma}) and Giraud's result.

Set
$$
F_{k}(x,y)=K(x,y)+\sum_{j=1}^{k}\int_M \Gamma_j(x,z)K(z,y)\dv_g(z)\,.
$$
{\it{Claim 2.}} 
For any $\phi\in C^2(M)$ and  $x\in M$, and for all $k=1,2,...$, we have
\ba\label{claim3}
\phi(x)=
&-\int_M F_k(x,y)L_g\phi(y)\dv_g(y)-\int_{\d M} F_k(x,y)B_g\phi(y)\ds_g(y)
\\
&+\int_M\Gamma_{k+1}(x,y)\phi(y)\dv_g(y)-\int_{\d M} \frac{n-2}{2(n-1)}H_g(y)F_k(x,y)\phi(y)\ds_g(y)\,.\notag
\end{align}

Claim 2 can be proved by induction on $k$.

\vspace{0.2cm}\noindent
{\it{Claim 3.}} 
For any $x\in M$ and $0<\theta<1$, the function $y\mapsto F_{n}(x,y)$ is in $C^{1,\theta}(M\backslash \{x\})$ and satisfies
\begin{equation}\label{claim2:1}
|F_{n}(x,y)|\leq Cd_g(x,y)^{2-n}\,,
\:\:\:\:
|\nabla_{g,y}F_{n}(x,y)|_g\leq Cd_g(x,y)^{1-n}\,,
\end{equation}
and
\begin{equation}\label{claim2:2}
\frac{|\nabla_{g,y}F_{n}(x,y)-\nabla_{g,y'}F_{n}(x,y')|_g}{d_g(y,y')^{\theta}}\leq C'd_g(x,y)^{1-\theta-n}+C'd_g(x,y')^{1-\theta-n}\,.
\end{equation}
Here, $C=C(M,g,n)$ and $C'=C'(M,g,n,\theta)$.
In particular, for any $x\in \d M$, $y\mapsto \d F_{n}/\d \eta_{g,y}(x,y)$ defines a continuous function on $\d M\backslash \{x\}$.

As a consequence of Claim 3, we can choose $N$ large enough such that $y\mapsto H_g(y) F_{n}(x_0,y)$ is in $C^{1,\theta}(\d M)$ for $0<\theta<1$ and satisfies
\begin{equation}\label{estim:HF}
\|H_g(\cdot) F_{n}(x_0,\cdot)\|_{C^{1,\theta}(\d M)}\leq C(M,g,n,\theta)\,.
\end{equation}

Let us prove Claim 3.  Choose a $y\neq x$, and let $y_t$ be a smooth curve such that $y_0=y$. Then, for any $r>0$,
$$
\frac{d}{dt}\int_{M\backslash B_r(y)}\Gamma_j(x,z)K(z,y_t)\dv_g(z)=\int_{M\backslash B_r(y)}\Gamma_j(x,z)\frac{d}{dt}K(z,y_t)\dv_g(z)
$$
Choose $t_0=t_0(x,y)>0$ such that $\displaystyle\frac{1}{2}\leq\frac{d_g(x,y_t)}{d_g(x,y)}\leq \frac{3}{2}$ for all $t\in [0,t_0)$.

For any $r>0$ such that $2r<d_g(x,y)$ and $t\in (0,t_0)$,  we have 
\ba
\int_{B_r(y)}&\Gamma_j(x,z)\Big|\frac{K(z,y_t)-K(z,y)}{t}\Big|\dv_g(z)\notag
\\
&\leq
C\int_{B_r(y)}d_g(x,z)^{1-n}(d_g(z,y_t)^{1-n}+d_g(z,y)^{1-n})\dv_g(z)\notag
\\
&\leq
C2^{n-1}d_g(x,y)^{1-n}\int_{B_r(y)}(d_g(z,y_t)^{1-n}+d_g(z,y)^{1-n})\dv_g(z)\notag
\\
&\leq
C2^{n-1}(2^{n-1}+1)\,d_g(x,y)^{1-n}\int_{B_r(y)}d_g(z,y)^{1-n}\dv_g(z)\,,\notag
\end{align}
and the right-hand side goes to $0$ as $r\to 0$.
Hence,
\begin{equation}\label{eq:deriv:H}
\frac{d}{dt}\int_{M}\Gamma_j(x,z)K(z,y_t)\dv_g(z)=\int_{M}\Gamma_j(x,z)\frac{d}{dt}K(z,y_t)\dv_g(z)
\end{equation}
and the estimates in (\ref{claim2:1}) follow from Giraud's result.

Now,
\ba
\frac{1}{d_g(y,y')^{\theta}}
&\Big|
\int_{M}\Gamma_j(x,z)\frac{\d}{\d y_i}K(z,y)\dv_g(z)-\int_{M}\Gamma_j(x,z)\frac{\d}{\d y_i}K(z,y')\dv_g(z)
\Big|\notag
\\
&\leq
\int_{M}\Gamma_j(x,z)
\Big|
\frac
{\frac{\d}{\d y_i}K(z,y)-\frac{\d}{\d y_i}K(z,y')}
{d_g(y,y')^{\theta}}
\Big|
\dv_g(z)\notag
\\
&\leq
C\int_{M}d_g(x,z)^{1-n}(d_g(z,y)^{1-\theta-n}+d_g(z,y')^{1-\theta-n})\dv_g(z)\notag
\\
&\leq
C(M,g,n,\theta)(d_g(x,y)^{2-\theta-n}+d_g(x,y')^{2-\theta-n})\,,\notag
\end{align}
where we used Lemma \ref{lemma:holder} in the second inequality, and Giraud's result in the last one.

This proves Claim 3.

Using the  hypothesis $Q(M,\d M)>0$, we define $u_{x_0}\in C^{2,\theta}(M)$ as the unique solution of
\ba\label{eq:ux}
\begin{cases}
L_{g}u_{x_0}(y)=-\Gamma_{n+1}(x_0,y)\,,&\text{in}\:M\,,
\\
B_{g}u_{x_0}(y)=\frac{n-2}{2(n-1)}H_g(y)F_{n}(x_0,y)\,,&\text{on}\:\d M\,.
\end{cases}
\end{align}
It satisfies
\ba\label{estim:C2:ux}
\|u_{x_0}\|_{C^{2,\theta}(M)}
\leq C\|u_{x_0}\|_{C^{0}(M)}
&+C\|\Gamma_{n+1}(x_0,\cdot)\|_{C^{0,\theta}(M)}
\\
&+C\|H_g(\cdot)F_{n}(x_0,\cdot)\|_{C^{1,\theta}(\d M)}\,\notag
\end{align}
where $C=C(M,g,n,\theta)$  (see  \cite[Theorems 6.30 and 6.31]{gilbarg-trudinger}; see also \cite[Theorem 7.3]{agmon-douglis-niremberg}).

\vspace{0.2cm}\noindent
{\it{Claim 4.}} There exists $C=C(M,g,n,\theta)$ such that $\|u_{x_0}\|_{C^{2,\theta}(M)}\leq C$.

Indeed,  using  (\ref{claim3}) with $k=n$ and any
$\phi\in C^{2}(M)$, one can see that
$$
\sup_M|\phi|\leq C\sup_M|L_g\phi|+C\sup_{\d M}|B_g\phi|+C\|\phi\|_{L^2(M)}+C\|\phi\|_{L^2(\d M)}\,.
$$
Since $Q(M,\d M)>0$, there exists $C=C(M,g,n)$ such that
$$
\int_M\phi^2\dv_g+\int_{\d M}\phi^2\ds_g
\leq
C\int_M |L_g(\phi)\phi|\dv_g+C\int_{\d M} |B_g(\phi)\phi|\ds_g\,.
$$
Thus, the Young's inequality implies
$$
\int_M\phi^2\dv_g+\int_{\d M}\phi^2\ds_g
\leq
C\int_M L_g(\phi)^2\dv_g+C\int_{\d M} B_g(\phi)^2\ds_g\,.
$$
Hence,
$
\|\phi\|_{C^0(M)}\leq
C\|L_g\phi\|_{C^0(M)}+C\|B_g\phi\|_{C^0(\d M)}\,.
$
Setting $\phi=u_{x_0}$ and using the equations (\ref{eq:ux}), we see that
\begin{equation}\label{estim:C0:ux}
\|u_{x_0}\|_{C^0(M)}\leq
C\|\Gamma_{n+1}(x_0,\cdot)\|_{C^0(M)}+C\|H_g(\cdot)F_{n}(x_0,\cdot)\|_{C^0(\d M)}\,.
\end{equation}
Claim 4 follows from the estimates (\ref{estim:gamma}), (\ref{estim:claim1}),
(\ref{estim:HF}), (\ref{estim:C2:ux}), and (\ref{estim:C0:ux}).

We define the function $G_{x_0}\in C^{1,\theta}(M\backslash \{x_0\}$ by
$$
G_{x_0}(y)=K(x_0,y)+\sum_{k=1}^{n}\int_M\Gamma_{i}(x_0,z)K(z,y)\dv_g(z)+u_{x_0}(y)\,.
$$
One can check that the formula (\ref{green:formula:0}) holds.

\vspace{0.2cm}\noindent
{\it{Claim 5.}} We have $G_{x_0}\in C^{\infty}(M\backslash\{x_0\})$ and 
\begin{equation}\label{eq:G}
\begin{cases}
L_{g}G_{x_0}=0\,,&\text{in}\:M\backslash \{x_0\}\,,
\\
B_{g}G_{x_0}=0\,,&\text{on}\:\d M\backslash \{x_0\}\,.
\end{cases}
\end{equation}

In order to prove Claim 5, we rewrite (\ref{form:green:H}) as
\begin{align}\label{form:green:H:L}
\int_M K(x,y)&L_g\phi(y)\dv_g(y)
+\int_{\d M}K(x,y)B_g\phi(y)\ds_g(y)
\\
&=\int_M L_{g,y}K(x,y)\phi(y)\dv_g(y)\notag
\\
&\hspace{0.5cm}-\phi(x)
-\int_{\d M}\frac{n-2}{2(n-1)}H_g(y)K(x,y)\phi(y)\ds_g(y)\,.\notag
\end{align}
Thus,
\begin{align}
&\int_M\left\{\int_M\Gamma_j(x,z)K(z,y)\dv_g(z)\right\}L_g\phi(y)\dv_g(y)
\\
&\hspace{1cm}+\int_{\d M}\left\{\int_M\Gamma_j(x,z)K(z,y)\dv_g(z)\right\}B_g\phi(y)\ds_g(y)\notag
\\
&=\int_M \Gamma_j(x,z)\left\{
\int_M K(z,y)L_g\phi(y)\dv_g(y)+\int_{\d M}K(z,y)B_g\phi(y)\ds_g(y)
\right\}\dv_g(z)\notag
\\
&=\int_M \Gamma_j(x,z)\int_M L_{g,y}K(z,y)\phi(y)\dv_g(y)\dv_g(z)\notag
\\
&\hspace{1cm}-\int_M \Gamma_j(x,z)\left\{\int_{\d M}\frac{n-2}{2(n-1)}H_g(y)K(z,y)\phi(y)\ds_g(y)+\phi(z)\right\}\dv_g(z)\notag
\\
&=\int_M \left\{\int_M \Gamma_j(x,z)L_{g,y}K(z,y)\dv_g(z)-\Gamma_j(x,y)\right\}\phi(y)\dv_g(y)\notag
\\
&\hspace{1cm}-\int_{\d M}\left\{\int_M \Gamma_j(x,z)K(z,y)\dv_g(z)\right\}\frac{n-2}{2(n-1)}H_g(y)\phi(y)\ds_g(y)\,,\notag
\end{align}
where we used (\ref{form:green:H:L}) in the second equality. Hence, we proved that the equations
\begin{equation}\notag
\begin{cases}
L_{g,y}\int_M\Gamma_j(x,z)K(z,y)\dv_g(z)=\Gamma_{j+1}(x,y)-\Gamma_{j}(x,y)\,,&\text{in}\:M\,,
\\
B_{g,y}\int_M\Gamma_j(x,z)K(z,y)\dv_g(z)=-\frac{n-2}{2(n-1)}H_g(y)\int_M\Gamma_j(x,z)K(z,y)\dv_g(z)\,,&\text{on}\:\d M\,,
\end{cases}
\end{equation}
hold in the sense of distributions. Then it is easy to check that the equations (\ref{eq:G}) hold in the sense of distributions.
Since $G_{x_0}\in C^{1,\theta}(M\backslash\{x_0\})$, elliptic regularity arguments  imply that $G_{x_0}\in C^{\infty}(M\backslash\{x_0\})$. 
This proves Claim 5.

The property (P1) follows from  (\ref{claim2:1}) and Claim 4.  In order to
prove (P2), we use (\ref{estim:gamma}), (\ref{estim:gamma:h}), (\ref{eq:deriv:H}) and Claim 4.

\vspace{0.2cm}\noindent
{\it{Claim 6.}} The function $G_{x_0}$ is positive on $M\backslash\{x_0\}$.

Let us prove Claim 6.
Let
$$
G_{x_0}^ -=
\begin{cases}
-G_{x_0}\,,&\text{if}\:G_{x_0}<0\,,
\\
0\,,&\text{if}\:G_{x_0}\geq 0\,.
\end{cases}
$$
Since $G_{x_0}$ has support in $M\backslash\{x_0\}$, one has
\ba
0&=-\int_MG_{x_0}^-L_gG_{x_0}\dv_g-\int_{\d M}G_{x_0}^-B_gG_{x_0}\ds_g\notag
\\
&=\int_M\left(|\nabla_gG_{x_0}^-|_g^2+\frac{n-2}{4(n-1)}R_g(G_{x_0}^-)^2\right)\dv_g+\int_{\d M}\frac{n-2}{2(n-1)}H_g(G_{x_0}^-)^2\ds_g\,.\notag
\end{align}
By the hypothesis $Q(M,\d M)>0$, we have $G_{x_0}^-\equiv 0$ which implies  $G_{x_0}\geq 0$.

We now change the metric by a conformal positive factor $u\in C^{\infty}(M)$ such that $\tilde{g}=u^{\frac{4}{n-2}}g$ satisfies $R_{\tilde{g}}>0$ in $M$ and $H_{\tilde{g}}\equiv0$ on $\d M$ (see \cite{escobar3}). Observing the conformal properties (\ref{propr:L}) and (\ref{propr:B}), we see that $\tilde{G}=u^{-1}G_{x_0}\geq 0$ satisfies $L_{\tilde{g}}\tilde{G}=0$ in $M\backslash \{x_0\}$ and $B_{\tilde{g}}\tilde{G}=0$ on $\d M\backslash \{x_0\}$. Then the strong maximum principle implies $\tilde{G}>0$, proving Claim 6.

This finishes the proof of Proposition \ref{green:point}.
\ep

The next proposition extends our Green function to the set $M\times M\backslash D_M$, where $D_M=\{(x,x)\in M\times M\,;\:x\in M\}$.
In order to define $G_{x_0}$ for all points $x_0\in M$, we change conformally the background metric in such a way that $H_g\equiv 0$ on $\d M$ and $R_g>0$ in $M$ (see \cite{escobar3}).
\begin{proposition}\label{green:global}
There exists a  continuous function $G: M\times M\backslash D_M\to \R$ satisfying
\begin{equation}\label{green:formula}
\phi(x)=-\int_M G(x,y)L_g\phi(y)\dv_g(y)-\int_{\d M}G(x,y)B_g\phi(y)\ds_g(y)
\end{equation}
for any $\phi\in C^2(M)$ and $x\in M$. Moreover, the following properties hold:
\\\\
(Q1) For any $x,y\in M$ with $x\neq y$, we have $G(x,y)=G(y,x)$ and $G(x,y)>0$.
\\\\
(Q2) For each $x\in M$, the function $y\mapsto G(x,y)$ is in $C^{\infty}(M\backslash\{x\})$ and there exists $C=C(M,g,n)$ such that 
$$
|G(x,y)|\leq Cd_g(x,y)^{2-n}\:\:\:\:\text{and}\:\:\:\: |\nabla_{g,y} G(x,y)|_g\leq Cd_g(x,y)^{1-n}\,,
$$
for any $x,y\in M$ with $x\neq y$.
\end{proposition}
\begin{remark}
A conformal change of the metric does not affect the result of this proposition. More precisely, if we obtain $G(x,y)$ as above, then $\tilde{G}(x,y)=v(x)^{-1}v(y)^{-1}G(x,y)$ satisfies the conclusions of Proposition \ref{green:global} when we replace the metric $g$ by $\tilde{g}=v^{\frac{4}{n-2}}g$. Here, $0<v\in C^{\infty}(M)$. In this case, the formula (\ref{green:formula}) is clear when we use the conformal properties  (\ref{propr:L}) and (\ref{propr:B}).
\end{remark}
\bp Since  $H_g\equiv 0$ on $\d M$, the hypothesis (\ref{hyp:H}) is satisfied for any point
$x=x_0\in \d M$. Moreover, the construction in Proposition \ref{green:point} can be performed for
any other $x\in M\backslash \d M$, since we can always solve the equations (\ref{eq:ux}) using the
fact that $R_g>0$ in $M$. Then we define $G(x,y)=G_x(y)$ and the formula (\ref{green:formula})
follows from (\ref{green:formula:0}).

Here, we follow the notations of the proof of Proposition \ref{green:point} and set $u(x,y)=u_x(y)$.

As in the estimate (\ref{estim:C0:ux}), we have
$$
\|u_x-u_{x'}\|_{C^0(M)}\leq
C\|\Gamma_{n+1}(x,\cdot)-\Gamma_{n+1}(x',\cdot)\|_{C^0(M)}\,.\notag
$$

Using Claim 4, we obtain \ba |u(x,y)-u(x',y')| &\leq
|u_x(y)-u_{x'}(y)|+|u_{x'}(y)-u_{x'}(y')|\notag
\\
&\leq \|u_x-u_{x'}\|_{C^0(M)}+C\sup_{y\in M}|\nabla_gu_{x'}(y)|_gd_g(y,y')\notag
\\
&\leq C\|\Gamma_{n+1}(x,\cdot)-\Gamma_{n+1}(x',\cdot)\|_{C^0(M)}+Cd_g(y,y')\,,\notag
\end{align}
where the right-hand side goes to zero as $(x,y)\to (x',y')$ because $(x,y)\mapsto \Gamma_{n+1}(x,y)$ is continuous.
Hence,  $(x,y)\mapsto u(x,y)$ is continuous. From this we conclude that $G$ is continuous on $M\times M\backslash D_M$.

\vspace{0.2cm}\noindent
{\it{Claim 7.}} For any $x\neq y$ we have $G(x,y)=G(y,x)$.

In fact, given $0\leq f_1,f_2\in C_0^{\infty}(M\backslash \d M)$, we choose $\phi_1$ and $\phi_2$ satisfying
$$
\begin{cases}
L_g\phi_1=f_1\,,&\text{in}\:M\,,
\\
B_g\phi_1=0\,,&\text{on}\:\d M\,,
\end{cases}
$$
and
$$
\begin{cases}
L_g\phi_2=f_2\,,&\text{in}\:M\,,
\\
B_g\phi_2=0\,,&\text{on}\:\d M\,.
\end{cases}
$$
Then, by  (\ref{green:formula}) and Tonneli's theorem, \ba \int_M\int_M&
G(x,y)L_g\phi_1(y)L_g\phi_2(x)\dv_g(y)\dv_g(x)\notag
\\
&=-\int_M \phi_1(x)L_g\phi_2(x)\dv_g(x)\notag =-\int_M L_g\phi_1(y)\phi_2(y)\dv_g(y)\notag
\\
&=\int_M\int_M G(y,x)L_g\phi_1(y)L_g\phi_2(x)\dv_g(x)\dv_g(y)\,.
\end{align}
Thus,
$$
\int_M\int_M (G(x,y)-G(y,x))f_1(y)f_2(x)\dv_g(y)\dv_g(x)=0\,.
$$
Then we see that $G(x,y)=G(y,x)$ for all $x,y\in M\backslash \d M$ with $x\neq y$. Since the function
$(x,y)\mapsto G(x,y)-G(y,x)$ is continuous on $M\times M\backslash D_M$ and vanishes on
$\{(M\backslash \d M)\times(M\backslash \d M)\}\backslash D_M$, we see that $G(x,y)=G(y,x)$ for
all $x,y\in M$ with $x\neq y$. This proves Claim 7.

Property (Q2) is the property (P1) of Proposition \ref{green:point}. This proves Proposition \ref{green:global}. 
\ep

\begin{corollary}\label{green:corol}
Let $G:M\times M\backslash D_M\to \R$ be the Green function obtained in Proposition \ref{green:global}.
If $0<\a<1$, then
$$
T(f)(x)=\int_{\d M}G(x,y)f(y)\dv(y)
$$
defines a continuous linear map $T:L^1(\d M)\to W^{\a,1}(\d M)$. Here, $ W^{\a,1}(\d M)$ denotes the fractional Sobolev space.
\end{corollary}
\bp
Since
$$
|T(f)(x)|\leq C\int_{\d M}d_g(x,y)^{2-n}|f(y)|\ds_g(y)\,,
$$
we have $\|T(f)\|_{L^1(\d M)}\leq C\|f\|_{L^1(\d M)}$ for some $C=C(M,n,g)$. Moreover,
$$
\frac{T(f)(x)-T(f)(x')}{d_g(x,x')^{n-1+\a}}
=\int_{\d M}\frac{G(x,y)-G(x',y)}{d_g(x,x')^{n-1+\a}}f(y)\ds_g(y)\,.
$$
Let $\theta\in(\a,1)$. By Lemma \ref{lemma:holder} and (Q2) of Proposition \ref{green:global}, we have
$$
\frac{|G(x,y)-G(x,y')|}{d_g(y,y')^{\theta}}\leq Cd_g(x,y)^{2-\theta-n}+Cd_g(x,y')^{2-\theta-n}\,.
$$
Then
$$
\frac{|G(x,y)-G(x',y)|}{d_g(x,x')^{\theta}}=\frac{|G(y,x)-G(y,x')|}{d_g(x,x')^{\theta}}\leq
Cd_g(y,x)^{2-\theta-n}+Cd_g(y,x')^{2-\theta-n}\,.
$$
Hence,
\ba
&\iint_{\d M}\frac{|T(f)(x)-T(f)(x')|}{d_g(x,x')^{n-1+\a}}\ds_g(x,x')\notag
\\
&\hspace{0.5cm}\leq C\iiint_{\d M}
d_g(x,x')^{\theta-\a-n+1}(d_g(x,y)^{2-\theta-n}+d_g(x',y)^{2-\theta-n})|f(y)|\ds_g(x,x',y)\notag
\\
&\hspace{0.5cm}\leq C\int_{\d M}\left\{\int_{\d M}d_g(x',y)^{2-\a-n}\dv_g(x')\right\}|f(y)|\ds_g(y)\notag
\\
&\hspace{1cm}+C\int_{\d M}\left\{\int_{\d M}d_g(x,y)^{2-\a-n}\dv_g(x)\right\}|f(y)|\ds_g(y)\notag
\\
&\hspace{0.5cm}\leq C(M,n,g,\a)\|f\|_{L^1(\d M)}\,.\notag
\end{align}
This finishes the proof of Corollary \ref{green:corol}.
\ep


\noindent
{\bf{Acknowledgements.}}
I would like to thank  Professor F. Marques for first showing me the notion of mass in (\ref{def:mass}) and for helpful comments on an earlier version of this paper. I would also like to thank the hospitality of Professor A. Neves at Imperial College London where part of this research was carried out. While at Imperial College, I was supported by CAPES/Brazil and CNPq/Brazil grants, and while in Brazil, I was partially supported by a FAPERJ/Brazil grant.



\noindent
\textsc{Instituto de Matem\'{a}tica, Universidade Federal Fluminense (UFF), \\Rua M\'{a}rio Santos Braga S/N, 24020-140,  Niter\'{o}i-RJ, Brazil}
\\{\bf{almaraz@vm.uff.br}}

\end{document}